%% file: RAC.tex
\documentclass[letter,10pt]{article}

\input{my_commands}
\title{Managing Randomization in the Multi-Block Alternating Direction Method of Multipliers for Quadratic Optimization}
\author{%
Kre\v{s}imir Mihi\'{c}\footnote{Kresimir Mihic is with the School of Mathematics, The University of Edinburgh, UK; and Oracle Labs, Redwood Shores, CA, USA. Email: K.Mihic@sms.ed.ac.uk, kresimir.mihic@oracle.com} %
\and Mingxi Zhu\footnote{Mingxi Zhu is with the Graduate School of Business, Stanford University, USA. Email: mingxiz@stanford.edu}%
\and Yinyu Ye\footnote{Yinyu Ye is with the Department of Management Science and Engineering, School of Engineering, Stanford University, USA. Email: yyye@stanford.edu.}%
}

\begin{document}

\maketitle

\begin{abstract}




The Alternating Direction Method of Multipliers (ADMM) has gained a lot of attention for solving large-scale and objective-separable constrained optimization. However, the two-block variable structure of the ADMM still limits the practical computational efficiency of the method, because one big matrix factorization is needed at least once even for linear and convex quadratic programming (e.g.,\cite{stellato:2018, Zarepisheh:2018, zhang:2018}). 
This drawback may be overcome by enforcing a multi-block structure of the decision variables in the original optimization problem. Unfortunately, the multi-block ADMM, with more than two blocks, is not guaranteed to be convergent \cite{chen:2016}. On the other hand, two positive developments have been made: first, if in each cyclic loop one randomly permutes the updating order of the multiple blocks, then the method converges in expectation for solving any system of linear equations with any number of blocks \cite{sunYe:2015, sunYe:2019}. Secondly, such a randomly permuted ADMM also works for equality-constrained convex quadratic programming even when the objective function is not separable \cite{chen:2015}. The goal of this paper is twofold. First, we add more randomness into the ADMM by developing a randomly assembled cyclic ADMM (RAC-ADMM) where the decision variables in each block are randomly assembled. We discuss the theoretical properties of RAC-ADMM and show when random assembling helps and when it
hurts, and develop a criterion to guarantee that it converges almost surely.
Secondly, using the theoretical guidance on RAC-ADMM, we conduct multiple numerical tests on solving both randomly generated and large-scale benchmark quadratic optimization problems, which include continuous, and binary graph-partition and quadratic assignment, and 
selected machine learning problems. Our numerical tests show that the RAC-ADMM, with a variable-grouping strategy, could significantly improve the computation efficiency on solving most quadratic optimization problems.
\end{abstract}

\input{Introduction}

\input{RAC_Algorithm}
\input{RAC_Convergence}
\input{RAC_Grouping}
\input{RACQP}

\input{Empirical_Analysis}

\input{Summary}

\bibliography{references}
\bibliographystyle{siam}  

\end{document}

%% file: my_commands.tex
\usepackage[top=3cm, bottom=3cm, left=3cm, right=3cm, columnsep=0.75cm]{geometry}
\usepackage{caption}
\usepackage{varwidth}
\usepackage{color}
\usepackage{array}
\usepackage{longtable}
\usepackage{algpseudocode}
\algdef{SE}{Begin}{End}{\textbf{begin}}{\textbf{end}}
\usepackage{algorithmicx}
\usepackage{algorithm}
\usepackage{multirow}
\usepackage{graphicx,subcaption}
\usepackage[flushleft]{threeparttable} 
\usepackage{booktabs,caption}
\usepackage{cite}

\usepackage{bmpsize}
\usepackage{multirow}
\usepackage{epsfig}
\usepackage{amssymb}
\usepackage{amsmath}
\usepackage{amsfonts}
\usepackage{amsthm}
\usepackage{mathtools}
\usepackage{multimedia}
\usepackage{tikz}
\usepackage{lipsum}
\usetikzlibrary{positioning}
\usepackage{array,booktabs}
\usepackage{bm}
\usepackage{mathtools}
\usepackage{empheq}
\usepackage{pgfplots}
\usepackage[parfill]{parskip}
\usepackage{txfonts}

\usepackage{fancybox}
\usepackage{anyfontsize} 

\usepackage{sectsty}
\usepackage{subcaption}
\usepackage[utf8]{inputenc}
\usepackage[thinlines,thiklines]{easybmat}
\usepackage{hhline}
\usepackage{arydshln}
\pgfplotsset{compat=1.8}

\makeatletter
\DeclareMathSizes{\f@size}{11}{6}{6}
\makeatother

\DeclareMathOperator{\R}{\varmathbb{R}}
\DeclareMathOperator{\X}{\mathcal{X}}

\DeclareMathOperator{\Z}{\varmathbb{Z}}
\DeclareMathOperator{\E}{\mathbb{E}}

\DeclareMathOperator*{\argmin}{arg\,min}

\DeclareMathOperator{\e}{\mathbf e}
\DeclareMathOperator{\q}{\mathbf q}
\DeclareMathOperator{\x}{\mathbf x}
\DeclareMathOperator{\xh}{\hat {\mathbf x}}

\DeclareMathOperator{\bh}{\hat {\mathbf b}}
\DeclareMathOperator{\Ah}{\hat {A}}
\DeclareMathOperator{\y}{\mathbf y}

\DeclareMathOperator{\s}{\mathbf s}
\DeclareMathOperator{\m}{\mathbf m}

\DeclareMathOperator{\z}{\mathbf z}
\DeclareMathOperator{\w}{\mathbf w}
\DeclareMathOperator{\cc}{\mathbf c}
\DeclareMathOperator{\bb}{\mathbf b}
\DeclareMathOperator{\lb}{\mathbf l}
\DeclareMathOperator{\ub}{\mathbf u}
\DeclareMathOperator{\bz}{\mathbf 0}
\DeclareMathOperator{\bo}{\mathbf 1}
\DeclareMathOperator{\eig}{\hbox{eig}}
\DeclareMathOperator{\prob}{\hbox{Prob}}
\DeclareMathOperator{\vect}{\hbox{vec}}

\DeclareMathOperator{\cond}{\!\ |\ \!}
\DeclareMathOperator{\eq}{\!\ =\!}

\newcolumntype{a}{>{\columncolor{BlueGreen}}c}
\newcolumntype{b}{>{\columncolor{Dandelion}}c}
\newcolumntype{d}{>{\columncolor{GreenYellow}}c}
\def\VRHDW#1#2#3{\vrule height #1 depth #2 width #3}
\newcommand{\up}{\VRHDW{1.5em}{0em}{0em}}
\newcommand{\uph}{\VRHDW{1em}{0em}{0em}}

\newcommand\ddfrac[2]{\frac{\displaystyle #1}{\displaystyle #2}}

\newcommand{\half}{\frac{1}{2}}
\newcommand{\bhalf}{\frac{\beta}{2}}

\def\sgn{\mbox{Sign}}

\newcommand{\myeqdl}[1]{\begin{equation} {\begin{array}{rccl} #1 \end{array}}\end{equation}}

\newcommand{\myeqmodel}[1]{\begin{equation} {\begin{array}{cl} #1 \end{array}}\end{equation}}

\newcommand{\myeqmodeln}[1]{\begin{equation*} {\begin{array}{cl} #1 \end{array}}\end{equation*}}

\newcommand{\myeq}[1]{$  {#1} $}
\newcommand{\myeql}[1]{\begin{equation}  {#1} \end{equation}}
\newcommand{\myeqln}[1]{\begin{equation*}  {#1} \end{equation*}}
\newcommand{\myalg}[2]{\begin{equation} {#1 \left\{\begin{array}{l} #2 \end{array}\right.}\end{equation}}

\newcommand{\myalgn}[2]{\begin{equation*} {#1 \left\{\begin{array}{l} #2 \end{array}\right.}\end{equation*}}

\newcommand{\mycasesn}[2]{\begin{equation*} {#1 \begin{cases}\begin{array}{ll} #2 \end{array}\end{cases}}\end{equation*}}

\def\algnote{\par\parindent1em\ignorespaces\vspace{2pt}}
\def\algbody{\vspace{3pt}}
\def\algcaption#1{\caption{\hspace{16pt} #1}}
\def\ALGORITHM#1#2#3{{\algcaption{#1}}{\algbody#2}{\algnote#3}} 

\theoremstyle{definition}
\newtheorem{theorem}{Theorem}[section]
\newtheorem{corollary}{Corollary}[theorem] 
\newtheorem{lemma}[theorem]{Lemma} 
\newtheorem{assumption}[theorem]{Assumption}
\newtheorem{proposition}[theorem]{Proposition}
\newtheorem{example}{Example}[section]

\newcommand{\myasmp}[1]{\begin{assumption} {\it #1}\end{assumption}}
\newcommand{\mycol}[1]{\begin{corollary} {\it #1}\end{corollary}}
\newcommand{\myth}[1]{\begin{theorem} {\it #1}\end{theorem}}
\newcommand{\mylemma}[1]{\begin{lemma} {\it #1}\end{lemma}}
\newcommand{\myprop}[1]{\begin{proposition} {\it #1}\end{proposition}}
\newcommand{\myproof}[1]{\begin{proof} {#1 \qedsymbol}\end{proof}}


%% file: Introduction.tex
\section{Introduction}

In this paper we consider the linearly constrained convex minimization model with an objective function that is the sum of multiple separable functions and a coupled quadratic function:

 \myeqmodel{\label{eq:problem_model_gen}
 \min\limits_{\x} & \sum\limits_{i=1}^{p}f_i(\x_i) + \half \x^T H \x + \cc^T\x\\[0.3cm]
 \mbox{s.t.}& \sum\limits_{i=1}^{p}A_i\x_i =\bb\\[0.2cm]
 &\x\in\X 
 }
where \myeq{f_i:\R^{d_i}\mapsto (-\infty,+\infty]} are closed proper convex functions, $H\in\R^{n\times n}$ is a symmetric positive semidefinite matrix, vector $\cc\in\R^n$
and the problem parameters are the matrix $A=[A_1,\dots,A_p]$, $A_i\in\R^{m\times d_i}$, $i = 1,2,\dots, p$ with $\sum_{i=1}^{p} d_i = n$  and the vector $\bb\in\R^m$. The constraint set $\mathcal X$ is the Cartesian product of possibly non-convex real, closed, nonempty sets, ${\mathcal X} = {\mathcal X_1} \times \dots \times{\mathcal X_p}$, 
where ${\x_i\in\mathcal X_i} \subseteq \R^{d_i}$.

Problem (\ref{eq:problem_model_gen}) naturally arises from applications such as machine and statistical learning, image processing, portfolio management, tensor decomposition, matrix completion or decomposition, manifold optimization, data clustering and many
other problems of practical importance.
To solve problem (\ref{eq:problem_model_gen}), we consider in particular a randomly assembled multi-block and cyclic alternating direction method of multipliers (RAC-ADMM), a novel algorithm with which we hope to mitigate the problem of slow convergence and divergence issues of the classical alternating direction method of multipliers (ADMM) when applied to problems with cross-block coupled variables. 

ADMM was originally proposed in 1970's ({\cite{glowinski:2014, gabay:1976}) and after a long period without too much attention it has recently gained in popularity 
for a broad spectrum of applications \cite{forero:2011, ohlsson:2012, lai:2014, sun:2014, jiang:2015}. Problems successfully solved by ADMM range from classical linear programming (LP), semidefinite programming (SDP) and quadratically constrained quadratic programming (QCQP) applied to partial differential equations, mechanics, image processing, statistical learning, computer vision and similar problems (for examples see \cite{boyd:2011,peng:2012, tao:2011, mohan:2014, huang:2016, zhang:2018}) to emerging areas such as deep learning \cite{taylor:2016}, medical treatment \cite{Zarepisheh:2018} and social networking \cite{baingana:2015}. ADMM is shown to be a good choice for problems where high accuracy is not a requirement but a “good enough” solution is needed to be found fast.

Cyclic multi-block ADMM is an iterative algorithm that embeds a Gaussian-Seidel decomposition into each iteration of the augmented Lagrangian method (ALM) (\cite{hestenes:1969,powell:1978}). It consists of a cyclic update of the blocks of primal variables, $x_i\in\X_{i}$, $x=(x_1,\dots,x_p)$, and a dual ascent type update of the variable $y\in\R^m$, i.e.,
\myalg{\label{ADMM}\mbox{Cyclic multi-block ADMM}:=}{
\x_1^{k+1}=\argmin_{\x_1} \{L_{\beta}(\x_1,\x_2^k,\x_3^k,\dots,\x_p^k;\y^k)\,|\, \x_1\in \X_1\}, \\
\vdots \\
\x_p^{k+1}=\argmin_{\x_p} \{L_{\beta} (\x_1^{k+1},\x_2^{k+1},\x_3^{k+1},\dots,\x_p;\y^k)\,|\, \x_p\in \X_p\},\\[0.3cm]
\y^{k+1}=\y^k -\beta(\sum_{i=1}^{p}A_i\x_i^{k+1}- \bb).
}

Where $\beta > 0$ is a penalty parameter of the Augmented Lagrangian function \myeq{L_{\beta}}, 
\myeql{ \label{eq:argL} L_{\beta} (\x_1,\dots,\x_p;\y^k) := \sum\limits_{i=1}^{p}f_i(\x_i) + \half \x^T H \x + \cc^T\x 
-\y^T\bigl(\sum\limits_{i=1}^{p}A_i\x_i -\bb\bigr)
+ \frac{\beta}{2}
\big\|\sum_{i=1}^p A_i\x_i -\bb\big\|^2.
}
Note that the classical ADMM \cite{glowinski:2014, gabay:1976} admits only optimization problems that are separable in blocks of variables and with $p=2$. 

Another variant of multi-block ADMM was suggested in \cite{bert:2015}, where the authors introduce the distributed multi-block ADMM (D-ADMM) for separable problems. 
The method creates a Dantzig-Wolfe-Benders decomposition structure and sequentially solves a "master" problem followed by solving distributed multi-block "slave" problems.
It converts the multi-block problem into an equivalent two-block problem via variable splitting \cite{bert:1989} and performs a separate augmented Lagrangian minimization over $\x_i$.

\myalg{\label{ADMM-D}\mbox{Distributed multi-block ADMM}:=}{
\hbox{Update }\x_i, i=1,\dots,p\\
\hspace{10pt} \x_i^{k+1}=\argmin\limits_{\x_i\in\X_i} f_i(\x_i)-(\y^k)^T(A_i\x_i-\lambda_i^k)
+ \bhalf\|A_i\x_i-\lambda_i^k\|^2 \\[0.3cm]
\hbox{Update }\lambda_i, i=1,\dots,p\\
\hspace{10pt} \lambda_i^{k+1}=A_i\x_i^{k+1}-\frac{1}{p}\big(\sum_{j=1}^pA_j\x_j^{k+1}-\bb\big)\\[0.3cm]
\y^{k+1}=\y^k -\frac{\beta}{p}(\sum_{i=1}^{p}A_i\x_i^{k+1}- \bb).
}
Because of the variable splitting, the distributed ADMM approach based on (\ref{ADMM-D}) increases the number of variables and constraints in the problem, which in turn makes the algorithm not very efficient for large $p$ in practice. 
 In addition, the method is not provably working for solving problems with non-separable objective functions.

The classical two-block ADMM (Eq. \ref{ADMM} with $p=2$) and its convergence have been extensively studied in the literature (e.g. \cite{gabay:1976, eckstein:1992, he:2012, monteiro:2013, deng:2016}. However, the two-block variable structure of the ADMM still limits the practical computational efficiency of the method, because one factorization of a large matrix is needed at least once even for linear and convex quadratic programming (e.g.,\cite{stellato:2018, zhang:2018}). This drawback may be overcome by enforcing a multi-block structure of the decision variables in the original optimization problem.
Indeed, due to the simplicity and practical implications of a direct extension of ADMM to the multi-block variant (\ref{ADMM}), an active research recently has been going on in developing ADMM variants with provable convergence and competitive numerical efficiency and iteration simplicity (e.g. \cite{chen:2017, he:2012, hong:2017, peng:2012}), and on proving global convergence under some special conditions (e.g. \cite{lin:2015,lin:2016,esser:2010,cai:2014}). Unfortunately, in general the Cyclic multi-block ADMM, with more than two blocks, is not guaranteed to be convergent even for solving a single system of linear equations, which settled a long-standing open question \cite{chen:2016}.

Moreover, in contrast to the work on separable convex problems, little work has been done on understanding properties of the multi-block ADMM for (\ref{eq:problem_model_gen}) with a non-separable convex quadratic or even non-convex objective function. One of the rare works that addresses coupled objectives is \cite{chen:2017} where authors describe convergence properties for non-separable convex minimization problems. A good description of the difficulties of obtaining a rigorous proof is given in \cite{eckstein:2012}.
For solving non-convex problems, a rigorous analysis of ADMM is by itself a very hard problem, with only a couple of works being done for generalized, but still limited (by an objective function), separable problems. For examples see \cite{wang:2017,hong:2016,zhang:2017,jiang:2018,wang:2015}.

Randomization is commonly used to reduce information and computation complexity for solving large-scale optimization problems. Typical examples include Q-Learning or Reinforced Learning, Stochastic Gradient Descent (SGD) for
Deep Learning, Randomized Block-Coordinate-Descent (BCD) for convex programming, and so on. Randomization of ADMM has recently become a matter of interest as well. In \cite{sunYe:2015} the authors devised randomly permuted multi-block ADMM (RP-ADMM) algorithm, in which on every cyclic loop the blocks are solved or updated in a randomly permuted order. Surprisingly the algorithm eliminated the divergence example constructed in \cite{chen:2016}, and RP-ADMM  was shown to converge linearly in expectation for solving any square system of linear equations with any number of blocks. Subsequently, in \cite{chen:2017} the authors focused on solving the linearly constrained convex optimization with coupled convex quadratic objective, and proved the convergence in expectation of RP-ADMM for the non separable multi-block convex quadratic programming, which is a much broader class of computational problems.
\myalg{\label{RP-ADMM}\mbox{RP-ADMM}:=}{
\hbox{Randomly permute } (1,2,...,p) \hbox{ into } (\sigma_1,\sigma_2,...,\sigma_p)\hbox{,}\\[-3pt]
\hbox{then solve }\\
\x_{\sigma_1}^{k+1}=\argmin_{\x_{\sigma_1}} \{L_{\beta}(\x_{\sigma_1},\x_{\sigma_2}^k,x_{\sigma_3}^k,\dots,\x_{\sigma_p}^k,\y^k)\,|\, \x_{\sigma_1}\in X_{\sigma_1}\}, \\
\vdots\\
\x_{\sigma_p}^{k+1}=\argmin_{\x_{\sigma_p}} \{L_{\beta} (\x_{\sigma_1}^{k+1},\x_{\sigma_2}^{k+1},x_{\sigma_3}^{k+1}\dots,\x_{\sigma_p},\y^k)\,|\, \x_{\sigma_p}\in X_{\sigma_p}\},\\[0.3cm]
\y^{k+1}=\y^k -\beta(A\x^{k+1}- \bb).
}

The main goal of the work proposed in this paper is twofold. First, we add more randomness into the ADMM by developing a randomly assembled cyclic ADMM (RAC-ADMM) where the decision variables in each block are randomly assembled. In contrast to RP-ADMM in which the variables in each block are fixed and unchanged, 
RAC-ADMM randomly assembles new blocks at each cyclic loop. It can be viewed as a decomposition-coordination procedure that decomposes the problem in a random fashion and combines the solutions to small local sub-problems to find the solution to the original large-scale problem. RAC-ADMM, in-line with RP-ADMM, admits multiple blocks with possibly cross-block coupled variables and updates the blocks in the cyclic order. 
The idea of re-constructing block variables at each cyclic loop was first mentioned in \cite{mihic:2018}, where the authors present a framework for solving discrete optimization problems which decomposes a problem into sub-problems by randomly (without replacement) grouping variables into subsets. Each subset is then used to construct a sub-problem by considering variables outside the subset as fixed, and the sub-problems are then solved in a cyclic fashion. Subsets are constructed once per iteration. The algorithm presented in that paper is a variant of the block coordinate descent (BCD) method with an addition of methodology to handle a small number of special constraints, which can be seen as a special case of RAC-ADMM. In the current paper we discuss the theoretical properties of RAC-ADMM and show when the additional random assembling helps and when it hurts. 

Secondly, using the theoretical guidance on RAC-ADMM, we conduct multiple numerical tests on solving both randomly generated and bench-mark quadratic optimization problems, which include continuous, and  binary graph-partitioning and quadratic assignment problems, and 
selected machine learning problems such as linear regression, LASSO, elastic-net, and support vector machine. Our numerical tests show the RAC-ADMM, with a systematic variable-grouping strategy (designate a set of variables always belonging to a same block), could significantly improve the computation efficiency on solving most quadratic optimization problems.

The current paper is organized as follows. In the next section we present RAC-ADMM algorithm and present theoretical results with respect to convergence. Next we discuss the notion of special grouping, thus selecting variables in less-random fashion by analyzing a problem structure, and the use of partial Lagrangian, approaches, which improve convergence speed of the algorithm. In Section \ref{sect:solver}, we present a solver , RACQP, we built that uses RAC-ADMM to address linearly constrained quadratic problems. The solver is implemented in Matlab \cite{matlab} and the source code available online \cite{RACQP:code}. The solver's performance is investigated in Section \ref{sect:num}, where we compare RACQP with commercial solvers, Gurobi \cite{gurobi} and Mosek \cite{mosek}, and the academic OSQP which is a ADMM-based solver developed by \cite{stellato:2018}. The summary of our contributions with concluding remarks is given in Section \ref{sec:summary}.

%% file: RAC_Algorithm.tex
\section{RAC-ADMM}
\label{sect:RAC-ADMM}

In this section we describe our randomly assembled cyclic alternating direction method of multipliers (RAC-ADMM). We start by presenting the algorithm, then analyze its convergence for linearly constrained quadratic problems, and finalize the section by introducing accelerated procedures that improve the convergence speed of RAC-ADMM by means of a grouping strategy of highly coupled variables and a partial Lagrangian approach. Note that although our analysis of convergence is restricted to quadratic and/or special classes of problems, it serves as a good indicator of the convergence of the algorithm in more general case.

\subsection{The algorithm}
\label{subsect:RAC_alg}
RAC-ADMM is an algorithm that is applied to solve convex problems (\ref{eq:problem_model_gen}). The algorithm addresses equality and inequality constraints separately, with the latter converted into equalities using slack variables, $\s$:
\myeqmodel{\label{eiq:problem_model_gen}
 \min\limits_{\x,\s} & f(\x) =\sum\limits_{i=1}^{p}f_i(\x_i) + \half \x^T H \x + \cc^T\x \\[0.3cm]
 \mbox{s.t.}& A_{eq}\x   =\bb_{eq} \\[0.2cm]
                 & A_{ineq}\x + \s= \bb_{ineq}\\[0.2cm]
                 &\x\in\X, \s\ge \bz
 }
where  matrix $A_{eq}\in\R^{m_e\times n}$ and vector $\bb_{eq}\in\R^{m_e}$ describe equality constraints and matrix $A_{ineq}\in\R^{m_i\times n}$ and the vector $\bb_{ineq}\in\R^{m_i}$ describe inequality constraints. Primal variables $\x\in\X$ are in constraint set $\mathcal X\subseteq \R^{n}$ which is the Cartesian product of possibly non-convex real, closed, nonempty sets, and slack variables $\s\in\R^{m_i}_{+}$. 
The augmented Lagrangian function used by RAC-ADMM is then defined by
\myeql{ \label{eq:argL_rac} 
\begin{array}{cl} 
L_{\beta} (\x;\s;\y_{eq};\y_{ineq}) := &f(\x)
-\y_{eq}^T\bigl(A_{eq}\x -\bb_{eq}\bigr)
-\y_{ineq}^T\bigl(A_{ineq}\x +\s -\bb_{ineq}\bigr)\\
&+ \frac{\beta}{2} \big(\big\|A_{eq}\x -\bb_{eq}\big\|^2
+ \big\|A_{ineq}\x +\s -\bb_{ineq}\big\|^2\big)
\end{array}
}
with dual variables $\y\in\R^{m_e}$ and $\z\in\R^{m_i}$, and penalty parameter $\beta > 0$. In (\ref{eiq:problem_model_gen}) we keep inequality and equality constraint matrices separate so to underline a separate slack variable update step of (\ref{RAC-ADMM}) which has a close form solution described in more details in Section \ref{sect:solver}.

RAC-ADMM is an iterative algorithm that embeds a Gaussian-Seidel decomposition into each iteration of the augmented Lagrangian method (ALM). It consists of a cyclic update of randomly constructed blocks$^\dagger$ of primal variables, $\x_i\in\X_{i}$, followed by the update of slack variables $\s$ and a dual ascent type update for Lagrange multipliers $\y_{eq}$ and $\y_{ineq}$:
\myalg{\label{RAC-ADMM}\mbox{RAC-ADMM}:=}{\mbox{Randomly (without replacement) assemble primal variables in $\x$ $^\dagger$ into $p$ blocks} \\
\mbox{$\x_i$, $i=1,\dots,p$, then solve}:\\[0.3cm]
\x_1^{k+1}=\argmin\limits_{\x_1} \{L_{\beta}(\x_1,\x_2^k,\dots,\x_p^k;\s^k;\y_{eq}^k;\z_{ineq}^k)\,|\, \x_1\in X_1\}, \\[-.3cm]
\vdots \\
\x_p^{k+1}=\argmin\limits_{\x_p} \{L_{\beta} (\x_1^{k+1},\x_2^{k+1},\dots,\x_p;\s^k;\y_{eq}^k;\z_{ineq}^k)\,|\, \x_p\in X_p\},\\[0.3cm]
\s^{k+1}=\argmin\limits_{\s} \{L_{\beta} (\x_1^{k+1},\x_2^{k+1},\dots,\x_p^{k+1};\s;\y_{eq}^k;\z_{ineq}^k)\,|\, \s\ge0\},\\[0.3cm]
\y_{eq}^{k+1}=\y_{eq}^k -\beta(A_{eq}\x^{k+1} -\bb_{eq}),\\[0.3cm]
\y_{ineq}^{k+1}=\y_{ineq}^k -\beta(A_{ineq}\x^{k+1} +\s^{k+1} -\bb_{ineq}). 
}

\vspace{-.3cm}\hspace{2.5cm}{\footnotesize $\dagger$ structure of a problem, if known, can be used to guide grouping  as described in Section \ref{sect:detect_structure}}

Randomly assembled cyclic alternating direction method of multipliers
(RAC-ADMM), can be seen as a generalization of cyclic ADMM, i.e. cyclic multi-block ADMM is a special case of RAC-ADMM in which the blocks are constructed at each iteration using a deterministic rule and optimized following a fixed block order. Using the same analogy, RP-ADMM can be seen as a special case of RAC-ADMM, in which blocks are constructed using some predetermined rule and kept fixed at each iteration, but sub-problems (i.e. blocks minimizing primal variables) are solved in a random order. 

The main advantage of RAC-ADMM over other multi-block ADMM variants is in its potential to significantly reduce primal and, especially, dual residuals, which is a common obstacle for applying multi-block ADMMs. To illustrate this feature we ran a simple experiment in which we fix the number of iterations and check the final residuals among the aforementioned multi-block ADMM variants. 

\begin{figure}[h!]
    \centering
    \begin{subfigure}[t]{0.5\textwidth}
        \centering
        \includegraphics[height=2.1in]{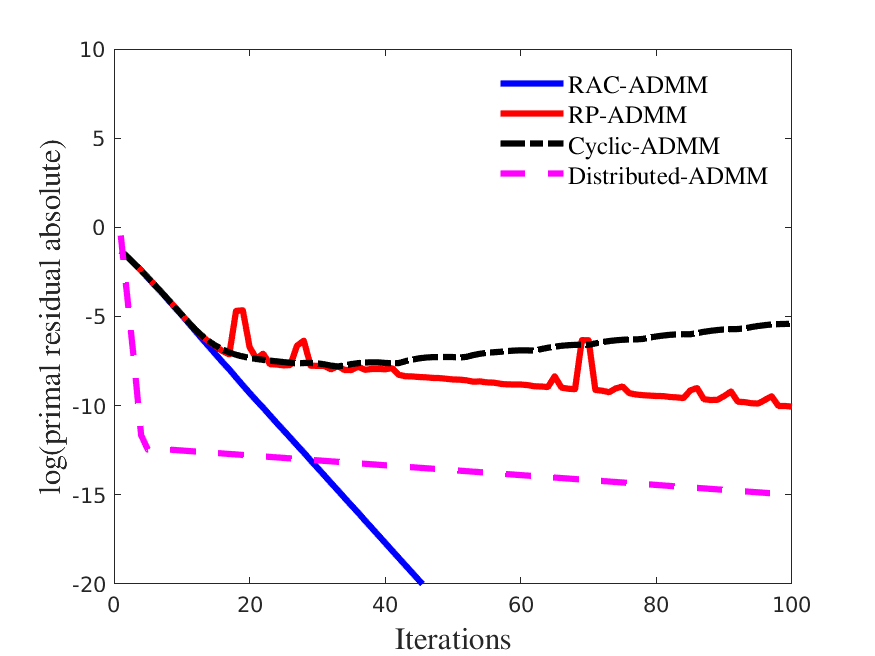}
        \caption{Primal residual}
    \end{subfigure}%
    \begin{subfigure}[t]{0.5\textwidth}
        \centering
        \includegraphics[height=2.1in]{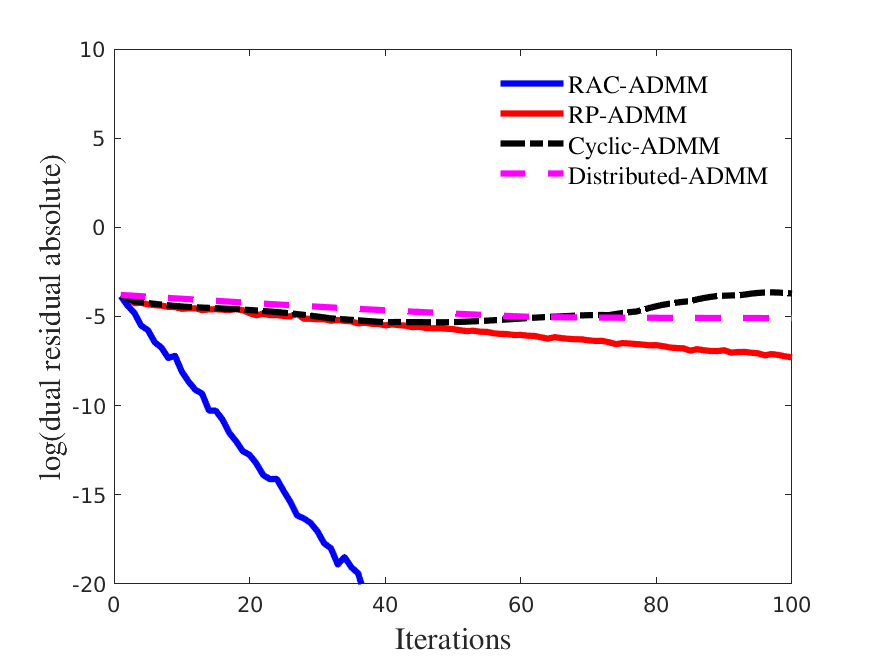}
        \caption{Dual residual}
    \end{subfigure}
    \caption{Iteration evolution of primal and dual residuals of ADMM variants}
\label{fig:admm_var}
\end{figure}

In Table \ref{tbl:RAC_Markov_comp} we show performance of the ADMMs when solving a simple quadratic problem with a single constraint, represented by a regularized Markowitz min-variance problem (defined in Section \ref{sect:rnd_lcqp}). Figure \ref{fig:admm_var} gives the insight in evolution of the both residuals with iterations. From the figure, it is noticeable that both D-ADMM (Eq. \ref{ADMM-D}) and RP-ADMM (Eq. \ref{RP-ADMM}) suffer from a very slow convergence speed, with the main difference that the latter gives a slightly lower error on dual residual. Multi-block Cyclic-ADMM (Eq. \ref{ADMM}) does not converge to a KKT point for any $k$, but oscillates  around a very weak solution. RAC-ADMM converges to the KKT solution very quickly with both residual errors below 10$^{-8}$ in less than 40 iterations.

\begin{table}[h!]
\footnotesize
\centering
\begin{tabular}{crccccc}  
\toprule 
\multirow{3}{*}{ADMM Variant}&\multicolumn{2}{c}{$k=$10 iterations} &\multicolumn{2}{c}{$k=$50 iterations} &\multicolumn{2}{c}{$k=$100 iterations}\\
\cmidrule(lr){2-7}
& primal & dual & primal & dual & primal & dual \\
\midrule
\up 
RAC-ADMM &7.2$\cdot$10$^{-3}$ & 3.1$\cdot$10$^{-4}$ & 3.0$\cdot$10$^{-10}$ & 4.6$\cdot$10$^{-12}$ & 1.2$\cdot$10$^{-14}$ & 4.4$\cdot$10$^{-16}$\\
RP-ADMM &7.4$\cdot$10$^{-3}$ & 1.0$\cdot$10$^{-2}$ & 2.0$\cdot$10$^{-4}$ & 3.3$\cdot$10$^{-3}$ & 4.3$\cdot$10$^{-5}$ & 6.8$\cdot$10$^{-4}$\\
Cyclic Multi-Block ADMM &7.4$\cdot$10$^{-3}$ & 1.2$\cdot$10$^{-2}$ & 6.8$\cdot$10$^{-4}$ & 4.9$\cdot$10$^{-3}$ & 4.5$\cdot$10$^{-3}$ & 2.5$\cdot$10$^{-2}$\\
Distributed Multi-block ADMM &3.7$\cdot$10$^{-6}$ & 1.8$\cdot$10$^{-2}$ & 1.2$\cdot$10$^{-6}$ & 8.0$\cdot$10$^{-3}$ & 3.1$\cdot$10$^{-7}$ & 6.2$\cdot$10$^{-3}$\\

\bottomrule
\end{tabular}
\caption{Primal and dual residuals of the result returned by ADMM variants after $k$ iterations for  a randomly generated Markowitz min-variance problem. Problem size $n=3000$, penalty parameter $\beta=1$.
}\label{tbl:RAC_Markov_comp}
\end{table}

%% file: RAC_Convergence.tex
\subsection{Convergence of RAC-ADMM}
\label{subsect:RAC_conv}
This section concerns with convergence properties of RAC-ADMM when applied to unbounded (i.e. $\x\in\R^n$) linearly-equality constrained quadratic optimization problems. 
To simplify the notation, we use $A=A_{eq}$ and $\bb=\bb_{eq}$.

\myeqmodel{\label{eq:problem_model_qp}
 \min\limits_{x} & \half \x^T H \x + \cc^T\x\\[0.3cm]
 \mbox{s.t.}& A\x =\bb 
 }

with $H\in\R^{n\times n}, H\succeq 0$, $\cc\in\R^n$, $A\in\R^{m\times n}$, $\bb\in\R^m$ and $\x\in\R^n$. 

Convergence analysis of problems that include inequalities (bounds on variables and/or inequality constraints) is still an open question and will be addressed in our subsequent work. 
\subsubsection{Preliminaries}
I) {\it Double Randomness Interpretation}

Let $\Gamma_{RAC(n,p)}$ denote all possible updating combinations for RAC with $n$ variables and $p$ blocks, and let $\sigma_{RAC}\in\Gamma_{RAC(n,p)}$ denote one specific updating combination for RAC-ADMM. Then the total  number of updating combinations for RAC-ADMM is given by
\myeqln{|\Gamma_{RAC(n,p)}|=\dfrac{n!}{(s!)^{p}}}
where $s\in\Z_+$ denotes size of each block with $p\cdot s=n$. 

RAC-ADMM could be viewed as a double-randomness procedure based on RP-ADMM with different block compositions. 
Let $\sigma_{RP}\in\Gamma_{RP(p)}$ denote an updating combinations of RP-ADMM with $p$ blocks where the variable composition in each block is fixed. Clearly, the total number of updating combinations for RP-ADMM is given by 
\myeqln{|\Gamma_{RP(p)}|=p!} 
the total number of possible updating orders of the $p$ blocks.
Then, one may consider RAC-ADMM first randomly chooses a block composition and then applies RP-ADMM.
Let $\upsilon_i\in \Upsilon(n,p)$ denote one specific block composition or partition of $n$ decision variables into $p$ blocks, where $\Upsilon(n,p)$ is the set of all possible block compositions. Then, the total number of all possible block compositions is given by
\myeqln{|\Upsilon(n,p)|=\dfrac{|\Gamma_{RAC(n,p)}|}{|\Gamma_{RP(p)}|}=\dfrac{n!}{(s!)^{p}p!}}
For convenience, in what follows let $\Gamma_{RP(p),\upsilon_i}$ denote all possible updating orders with a fixed block composition
$\upsilon_i$. 

To further illustrate the relations of RP-ADMM and RAC-ADMM, consider the following simple example.
\example{
Let $n=6$, $p=3$, so $|\Gamma_{RP(6,3)}|=3!=6$, and the total number of block compositions or partitions is $15$:
\myeqln{
\begin{array}{cl}
\upsilon_{i}\in \Upsilon(6,3)=&\big \{ 
\{[x_1,x_2],[x_3,x_4],[x_5,x_6]\},  \{[x_1,x_2],[x_3,x_5],[x_4,x_6]\}, \{[x_1,x_2],[x_3,x_6],[x_4,x_5]\}, \\ 
&\{[x_1,x_3],[x_2,x_4],[x_5,x_6]\}, \{[x_1,x_3],[x_2,x_5],[x_4,x_6]\}, \{[x_1,x_3],[x_2,x_6],[x_4,x_5]\}, \\ \null
&\{[x_1,x_4],[x_2,x_3],[x_5,x_6]\},  \{[x_1,x_4],[x_2,x_5],[x_3,x_6]\}, \{[x_1,x_4],[x_2,x_6],[x_3,x_5]\}, \\ \null
&\{[x_1,x_5],[x_2,x_3],[x_4,x_6]\},   \{[x_1,x_5],[x_2,x_4],[x_3,x_6]\}, \{[x_1,x_5],[x_2,x_6],[x_3,x_4]\}, \\ \null
&\{[x_1,x_6],[x_2,x_3],[x_4,x_5]\},   \{[x_1,x_6],[x_2,x_4],[x_3,x_5]\}, \{[x_1,x_6],[x_2,x_5],[x_3,x_4]\} \big\}
\end{array}
}
}
RAC-ADMM could be viewed as if, at each cyclic loop, the algorithm first selects a block composition $\upsilon_i$ uniformly random from all possible $15$ block compositions $\Upsilon(n,p)$, and then performs RP-ADMM with the chosen specific block composition $\upsilon_i$. In other words, RAC-ADMM then randomly selects $\sigma\in\Gamma_{RP(p),\upsilon_i}$, which leads to a total of $90$ possible updating combinations.

\vspace{.3cm}
II) {\it RAC-ADMM as a linear transformation}

Recall that the augmented Lagrangian function for (\ref{eq:problem_model_qp}) is given by
\myeqln{L_{\beta}(\x,\y)=\half\x^TH\x+\cc^T\x-\y^T(A\x-\bb)+\half{\beta}||A\x-\bb||^2.}
Consider one specific update order generated by RAC, $\sigma_{RAC}\in\Gamma_{RAC(n,p)}$. Note that we use $\sigma$ instead $\sigma_{RAC}$ when there is no confusion. One possible update combination generated by RAC, $\sigma = [\sigma_1,\dots,\sigma_p]$, where $\sigma_i$ is an index vector of size $s$, is as follows, 
\myalgn{\label{RAC-ADDM_one}\mbox{RAC-ADMM}_{k+1}=}{
\x_{\sigma 1}^{k+1}=\argmin\limits_{\x_{\sigma 1}} \{L_{\beta}(\x_{\sigma 1},\x_{\sigma 2}^k,\dots,\x_{\sigma p}^k;\y^k)\}, \\[-.3cm]
\vdots \\
\x_{\sigma p}^{k+1}=\argmin\limits_{\x_{\sigma p}} \{L_{\beta}(\x_{\sigma 1}^{k+1},\x_{\sigma 2}^{k+1},\dots,\x_{\sigma p};\y^k)\}, \\[0.4cm]
\y^{k+1}=\y^k -\beta(A\x^{k+1} -\bb). 
}

For convenience, we follow the notation in \cite{chen:2017} and  \cite{sunYe:2015, sunYe:2019} to describe the iterative scheme of RAC-ADMM in a matrix form. Let $L_{\sigma}\in\R^{n\times n}$ be $s\times s$ block matrix defined with respect to $\sigma_i$ rows and $\sigma_j$ columns as
\mycasesn{(L_{\sigma})_{\sigma_i,\sigma_j}:=}
{H_{\sigma_i,\sigma_j} + \beta A^T_{\sigma_i}A_{\sigma_j},  & \quad i\geq j\\
0,  & \quad \textup{otherwise}
}
and let $R _{\sigma}$ be defined as
\myeqln{R_{\sigma}:= L_\sigma - (H + \beta A^TA).} 
By setting $\z:=(\x;\y)$, RAC-ADMM could be viewed as a linear system mapping iteration
\myeqln{\z^{k+1}:= M_{\sigma} \z^{k}+\bar{L}^{-1}_{\sigma}\bar{\bb}}
where
\myeql{\label{eq:mapping_Mi}
M_{\sigma}:= \bar{L}^{-1}_{\sigma}\bar{R}_{\sigma}}
and
\myeqln{
\begin{array}{ccc}
\bar{L}_{\sigma}:=\begin{bmatrix}L_{\sigma} & 0\\
\beta A & I\end{bmatrix} & 
\bar{R}_{\sigma}:=\begin{bmatrix}R_{\sigma} & A^T\\
0 & I\end{bmatrix} &
\bar{\bb}:=\begin{bmatrix}-\cc+\beta A^T\bb \\
\beta \bb\end{bmatrix} 
\end{array}
}

Define the matrix $Q$ by
\myeqln{
Q:=\E_{\sigma}(L^{-1}_{\sigma})=\frac{1}{|\Gamma_{RAC(n,p)}|}\sum\nolimits_{\sigma\in\Gamma_{RAC(n,p)}}L_{\sigma}^{-1}=\frac{1}{|\Upsilon(n,p)|}\sum\nolimits_{\upsilon_i\in \Upsilon(n,p)}\left\{\frac{1}{p!}\sum\nolimits_{\sigma\in \Gamma_{RP(p),\upsilon_i}}L^{-1}_{\sigma}\right\}
}

Notice that for any block structure $\upsilon_{i}$ any update order within this fixed block structure $\sigma\in\Gamma_{RP(p),\upsilon_i}$, we have $L_{\sigma}^T=L_{\bar{\sigma}}$, where $\bar{\sigma}$ is a reverse permutation of $\sigma\in\Gamma_{RP(p),\upsilon_i}$. Specifically, let $\sigma=[\sigma_1,\dots,\sigma_p]$, we have $\bar{\sigma}=[\bar{\sigma}_1,\dots,\bar{\sigma}_p]$, and $\bar{\sigma}_i=\sigma_{p+1-i}$.
For a specific fixed block structure $\upsilon_i$, define matrix $Q_{\upsilon_i}$ as 
\myeqln{
Q_{\upsilon_i}:=\E(L_{\sigma}|\upsilon_i)=\frac{1}{p!}\sum\nolimits_{\sigma_i\in \Gamma_{RP(n,\upsilon_i)}}L^{-1}_{\sigma_i},}
and because $L_{\sigma}^T=L_{\bar{\sigma}}$, matrix $Q_{\upsilon_i}$ is symmetric for all $i$, and
\myeql{\label{eq:P}
Q:=\frac{1}{\Upsilon(n,p)}\sum\nolimits_{\upsilon_i\in \Upsilon(n,p)} Q_{\upsilon_i}}
Finally, the expected mapping matrix $M$ is given by
\myeqln{M:=\E_{\sigma}(M_\sigma)=\frac{1}{|\Gamma_{RAC(n,p)}|}\sum\nolimits_{\sigma\in\Gamma_{RAC(n,p)}}M_{\sigma}}
or, by direct computation, 
\myeqln{
M:=\begin{bmatrix}
I-QS & QA^T\\
-\beta A+\beta AQS & I-\beta AQA^T
\end{bmatrix}}
where $S=H + \beta A^TA$.

\subsubsection{Expected convergence of RAC-ADMM}
With the preliminaries defined, we are now ready to show that RAC-ADMM converges in expectation under the following assumption:\\
\myasmp{
\label{assumption1}
Assume that for any block of indices $\sigma_i$ that generated by RAC-ADMM
\myeqln{H_{\sigma_i,\sigma_i} + \beta A^T_{\sigma_i}A_{\sigma_i}\succ0}
where $\sigma_i$ is the index vector describing indices of primal variables of the block $i$.
}

\vspace{.3cm}
\myth{
\label{expected convergence}
Suppose that Assumption (\ref{assumption1}) holds, and that RAC-ADMM (\ref{RAC-ADMM}) is employed to solve problem (\ref{eq:problem_model_qp}). Then the expected output converges to some KKT point of (\ref{eq:problem_model_qp}).
}
Theorem \ref{expected convergence} suggests that the expected output converges to some KKT point of (\ref{eq:problem_model_qp}). Such convergence in expectation criteria has been widely used in many randomized algorithms, including convergence analysis for RP-BCD and RP-ADMM (e.g. \cite{chen:2015,sunYe:2015}), and stochastic quasi-newton methods (e.g. \cite{byrd2016stochastic}). 
It is worth mentioning that if the optimization problem is strictly convex (H>0), we are able to prove that the expected mapping matrix has specturm that is strictly less than 1, following corollary \ref{corollary1}. 

Although convergence in expectation is widely used in many literature, it is still a relatively weak convergence criteria. Thih is why in section \ref{subsect:variance} we propose a sufficient condition for almost surely convergence of RAC-ADMM. The section also provides an example showing a problem with $\rho(M)<1$ which does not converge. Rather it oscillates almost surely (Example \ref{example_diverge}). 
To the best of our knowledge, this is the first example showing that even if a randomized optimization algorithm has expected spectrum radius \textit{strictly less than 1}, the algorithm may still oscillate --
to construct an example with expected spectrum radius equals to 1 that does not converge is an easy task. Consider for example a  sequence $\{x_t,t\geq0\}$ with $x_t=-1$ and $x_t=1$, chosen with equal probabilities (prob=1/2). Then, the sequence does not converge with probability 1. However, under the such example, the expected spectrum of this mapping procedure $\rho(M)$ actually equals to 1, which implies that the sequence may not converge.

Despite the fact that such example exists for RAC-ADMM, in all the numerical tests provided in section \ref{sect:num}, RAC-ADMM converges to the KKT point of the optimization problem under few iterations. Such strong numerical evidences imply that in practice, our algorithm does not require taking expectation over many iterations  to converge.

The proof of Theorem \ref{expected convergence} follows the proof structure of \cite{sunYe:2015, chen:2017, sunYe:2019} to show that under Assumption \ref{assumption1}:
\begin{enumerate}
 \item[(1)] $\eig(QS)\in[0,\frac{4}{3})$ ;
 \item[(2)] $\forall \lambda\in\eig(M), \eig(QS)\in[0,\frac{4}{3})\implies \|\lambda\|<1$ or $\lambda=1$;
 \item[(3)] if $1\in\eig(M)$, then the eigenvalue 1 has a complete set of eigenvectors;
 \item[(4)] Steps (2) and (3) imply the convergence in expectation of the RAC-ADMM.
\end{enumerate}

The proof builds on Theorem 2 from \cite{chen:2017}, which describes RP-ADMM convergence in expectation under specific conditions put on matrices $H$ and $A$, and Weyl's inequality, which gives the upper bound on maximum eigenvalue  and the lower bound on minimum eigenvalue of a sum of Hermitian matrices. Proofs for items (2) and (3) are identical to proofs given in \cite{chen:2017}, Section 3.2, so here the  focus in on proving item (1).

The following lemma completes the proof of expected convergence of RAC.
\vspace{.3cm}
\mylemma{
\label{lemma2}
Under assumption $\ref{assumption1}$, the matrix $Q$ is positive definite, and
\myeqln{eig(QS)\subset[0,\frac{4}{3})} 
}

To prove Lemma \ref{lemma2}, we first show that for any block structure $\upsilon_i$, the following proposition holds:
\myprop{
\label{Proposition 1} 
$Q_{\upsilon_i}S$ is positive semi-definite and symmetric, and 
\myeqln{eig(Q_{\upsilon_i}S)\subseteq [0,\frac{4}{3})}
}
Intuitively, a different block structure of RAC-ADMM iteration could be viewed as relabeling variables and performing RP-ADMM procedure as described in \cite{chen:2017}. 

\myproof{
Define block structure $\{[x_1,\dots,x_{s}],[x_{s+1},\dots,x_{2s}],[x_{(p-1)s+1},\dots,x_{ps}]\}$ as $\upsilon_1$. For any block structure $\upsilon_i$, there exists $\tilde{S}$ and $\tilde{Q}_{\upsilon_1}$ s.t. 
\myeqln{
\eig(Q_{\upsilon_i}S)=\eig(\tilde{Q}_{\upsilon_1}\tilde{S})}
where $\tilde{Q}_{\upsilon_1}$ represents formulation of $\E_{\sigma}(L_{\sigma}^{-1})$ matrix with respect to block structure $\upsilon_1$ and matrix $\tilde{S}$. To prove this, we introduce permutation matrix $P_{\upsilon_1\to \upsilon_i}$ as follows. 

Given 
\myeqln{
\begin{array}{ll}
\upsilon_1=& \{[1,\dots, {s}],[{s+1},\dots, {2s}],[{(p-1)s+1},\dots,{ps}]\}\\
\upsilon_i=& \{[\pi(1),\dots, \pi({s})],[\pi({s+1}),\dots, \pi({2s})],[\pi({(p-1)s+1}),\dots,\pi({ps})]\}
\end{array}
}
define
\myeqln{
P_{\upsilon_i} =\begin{bmatrix}
\e_{\pi(1)}\\
\e_{\pi(2)}\\
\vdots\\
\e_{\pi(ps)}
\end{bmatrix}
}
Where $\e_{i}$ is the row vector with $i^{th}$ element equal to 1. Notice $P_{\upsilon_i}$ is orthogonal matrix for any $\upsilon_i$, i.e. $P_{\upsilon_i}P_{\upsilon_i}^T=I$. For any fixed block structure $\upsilon_i$, with an update order within $\sigma_{RP}\in \Gamma_{RP}(p)$, the following equality holds
\myeqln{
L_{\sigma_{RP},S,\upsilon_i}=P_{\upsilon_i}^T L_{\sigma_{RP},\tilde{S},\upsilon_1}P_{\upsilon_i}
}
where $L_{\sigma_{RP},S,\upsilon_i}$ is the construction of $L$ following update order $\sigma_{RP}\in\Gamma_{RP}(p)$ and block structure $\upsilon_i$ with respect to $S$, and $ L_{\sigma_{RP},\tilde{S},\upsilon_1}$ is the construction of $L$ following update order $\sigma_{RP}\in\Gamma_{RP}(p)$ and block structure $\upsilon_1$, with coefficient matrix $\tilde{S}$, and
\myeqln{
\tilde{S}=P_{\upsilon_i} S P_{\upsilon_i}^T
}
and
\myeqln{
L_{\sigma,S,\upsilon_i}^{-1}=(P_{\upsilon_i}^T L_{\sigma,\tilde{S},\upsilon_1}P_{\upsilon_i})^{-1}=P_{\upsilon_i}^T L^{-1}_{\sigma,\tilde{S},\upsilon_1}P_{\upsilon_i}.
}
Then by the definition of $Q$ matrix (Eq. \ref{eq:P}), we get
\myeqln{
Q_{\upsilon_i,S}=P_{\upsilon_i}^T\tilde{Q}_{\upsilon_1,\tilde{S}}P_{\upsilon_i}
}
so that
\myeqln{
Q_{\upsilon_i,S}S=P_{\upsilon_i}^T \tilde{Q}_{\upsilon_1,\tilde{S}} P_{\upsilon_i} P_{\upsilon_i}^{-1} \tilde{S} P_{\upsilon_i}=P_{\upsilon_i}^T \tilde{Q}_{\upsilon_1,\tilde{S}}  \tilde{S} P_{\upsilon_i}.
}
Considering the eigenvalues of $Q_{\upsilon_i,S}S$, 
\myeqln{
\eig(Q_{\upsilon_i,S}S)=\eig(P_{\upsilon_i}^T \tilde{Q}_{\upsilon_1,\tilde{S}}  \tilde{S} P_{\upsilon_i})=\eig(\tilde{Q}_{\upsilon_1,\tilde{S}}  \tilde{S} )
}
and from \cite{chen:2017}, under Assumption $(\ref{assumption1})$, $\tilde{Q}_{\upsilon_1,\tilde{S}}$ is positive definite, and
\myeqln{
eig(\tilde{Q}_{\upsilon_1,\tilde{S}}  \tilde{S} )\subset [0,\frac{4}{3})
}
which implies $Q_{\upsilon_i}$ is positive definite, and
\myeqln{
eig(Q_{\upsilon_i} S)\subset [0,\frac{4}{3}).
}
Notice that by definition of $Q$, we have 
\myeqln{
QS=\frac{1}{\Upsilon(d,n)}\sum_{\upsilon_i}Q_{\upsilon_i}S
} 
and $Q_{\upsilon_i}S$ is positive definite and symmetric. Let $\lambda_1(A)$ denote the maximum eigenvalue of matrix $A$, then as all $Q_{\upsilon_i}S$ are Hermitian matrices, by Weyl's theorem, we have
\myeqln{
\lambda_1(QS)=\lambda_1(\frac{1}{\Upsilon(d,n)} \sum_{i\in \Upsilon}Q_{\upsilon_i}S)\leq \frac{1}{\Upsilon(d,n)}\sum_{i\in \Upsilon}\lambda_1(Q_{\upsilon_i}S)
}
and as $\lambda_1(Q_{\upsilon_i}S)<\frac{4}{3}$ for each $i$, 
\myeqln{\eig(QS)\subseteq [0,\frac{4}{3})
}
which completes the proof of Lemma \ref{lemma2}, and thus establishes that RAC-ADMM is guaranteed to converge in expectation. }

When the problem is strongly convex ($H\succ0$), we introduce the following corollary.
\vspace{.3cm}
\mycol{
\label{corollary1}
Under assumption $\ref{assumption1}$, and $H\succ0$, 
\myeqln{\rho(M)<1}
}
\myproof{
When $H\succ0$, by definition $S=H+\beta A^TA\succ0$, and by Lemma \ref{lemma2}, $Q\succ0$, hence $\eig(QS)\subseteq (0,\dfrac{4}{3})$, and this implies $\rho(M)<1$.
}
Note that there are random sequences converging in expectation where their spectrum-radius equal to one. Therefore, for solving strongly non-separable convex quadratic optimization, the expected convergence rate of RAC-ADMM is proved to be linear, which result is stronger than just "convergence in expectation".

\subsubsection{Convergence speed of RAC-ADMM vs. RP-ADMM}
Following is a corollary to show that on average or in expectation, RAC-ADMM performs RP-ADMM with a fixed block composition in sense of spectral radius of mapping matrix.

\vspace{.4cm}
\mycol{\label{col:rac_speed}
Under Assumption \ref{assumption1}, with $H=0$ so that $S=\beta A^TA$, where $A\in\R^{n\times n}$ is a non-singular matrix, there exists some RP-ADMM (with specific block compositions), such that expected spectral radius of RAC-ADMM mapping matrix is (weakly) smaller than expected spectral radius of that of RP-ADMM.
}

\myproof{
We prove the corollary in solving linear system with $A$ non singular, with null objective function. In this setup, the expected output converges to the unique primal dual optimal solution to (\ref{eq:problem_model_qp}).

Notice in this setup, we have
\myeqln{
\begin{array}{l}
\lambda\in\eig(M)\Leftrightarrow \tau = \dfrac{(1-\lambda)^2}{1-2\lambda} \in \eig(QA^TA)\\
\lambda_{\upsilon_i}\in\eig(M_{RP,\upsilon_i})\Leftrightarrow \tau_{\upsilon_i} = \dfrac{(1-\lambda_{\upsilon_i})^2}{1-2\lambda_{\upsilon_i}} \in \eig(Q_{\upsilon_i}A^TA)
\end{array}
}
By calculation, we could characterize $\lambda$ as roots of quadratic polynomial \cite{sunYe:2019},
\myeqln{
\lambda_1 = 1-\tau+\sqrt{\tau(\tau-1)}, \quad \lambda_2 = 1-\tau-\sqrt{\tau(\tau-1)}.
}
Suppose corollary doesn't hold, $\rho(E(M_{RAC}))\geq \rho(\E(M_{RP,\upsilon_i}))$ for all possible block structure. Define $\underline {\tau}_{\upsilon_i}$ as the the smallest eigenvalue with respect to $Q_{\upsilon_i}S$, and $\bar{\tau}_{\upsilon_i}$ as the largest eigenvalue with respect to $Q_{\upsilon_i}S$. Similarly, $\underline{\tau}$ as the smallest eigenvalue with respect to $QS$, and $\bar{\tau}$ the largest eigenvalue of $QS$. Consider the following two cases. 

\underline {Case 1.} \hspace{10pt} $\lambda^*=\max_i|\lambda_i|\in \mathbb{C} \textup{ and }  \lambda^*\notin\R \Leftrightarrow \tau_{\lambda^*}<1$, where $\tau_{\lambda^*}\in\eig(QS)$ satisfies $\dfrac{(1-\lambda^*)^2}{1-2\lambda^*}=\tau_{\lambda^*}$.

We have, $\rho(E(M_{RAC}))\geq \rho(\E(M_{RP,\upsilon_i}))\ \forall \ i$, which implies that
\myeqln{
\sqrt{1-\tau_{\lambda^*}}>\max\{\sqrt{1-\tau_{\upsilon_i}},\tau_{\upsilon_i}-1+\sqrt{\tau_{\upsilon_i}(\tau_{\upsilon_i}-1)}\} \quad \forall i.
} 
Specifically
\myeqln{
\sqrt{1-\tau_{\lambda^*}}>\sqrt{1-\underline{\tau}_{\upsilon_i}} \quad \forall \upsilon_i,
} 
As $f(x)=\sqrt{1-x}$ is monotone decreasing with respect to $x$, the above implies that
\myeqln{
\tau_{\lambda^*}<\underline{\tau}_{\upsilon_i} \quad \forall \upsilon_i,
} 
and as $\tau_{\lambda^*}\geq\underline{\tau}$, the above equation implies
\myeqln{
\underline{\tau}<\underline{\tau}_{\upsilon_i} \quad \forall \upsilon_i,
} 
which is impossible, as by Weyl's theorem,
\myeqln{
\underline{\tau}\geq \dfrac{1}{|\Upsilon(d,b)|}\sum_i\underline{\tau}_{\upsilon_i}\geq \min_{i}\  \underline{\tau}_{\upsilon_i}. 
}

\underline{Case 2.} \hspace{10pt} $\lambda^*=\max_i|\lambda_i|\in \R \Leftrightarrow \tau_{\lambda^*}>1$.
 
We have $\rho(E(M_{RAC}))\geq \rho(\E(M_{RP,\upsilon_i}))\ \forall i$, what implies that
\myeqln{
\tau_{\lambda^*}-1+\sqrt{\tau_{\lambda^*}(\tau_{\lambda^*}-1)}>\max\{\sqrt{1-\tau_{\upsilon_i}},\tau_{\upsilon_i}-1+\sqrt{\tau_{\upsilon_i}(\tau_{\upsilon_i}-1)}\} \quad \forall i.
} 
Specifically,
\myeqln{
\tau_{\lambda^*}-1+\sqrt{\tau_{\lambda^*}(\tau_{\lambda^*}-1)}>\tau_{\upsilon_i}-1+\sqrt{\tau_{\upsilon_i}(\tau_{\upsilon_i}-1)} \quad \forall i,
 \quad \forall \upsilon_i.
} 
As $g(x)=x-1+\sqrt{x(x-1)}$ is a monotone increasing function for $x\in [1,\infty)$, the above implies
\myeqln{
\overline{\tau}\geq \tau_{\lambda^*}>\overline{\tau}_{\upsilon_i} \quad \forall \upsilon_i,
} 
which is impossible, as by Weyl's theorem,
\myeqln{
\bar{\tau} \leq \dfrac{1}{|\Upsilon(d,b)|}\sum_i\bar{\tau}_{\upsilon_i}\leq \max_{i}\  \underline{\tau}_{\upsilon_i}. 
}
}

\subsubsection{Variance of RAC-ADMM}
\label{subsect:variance}
Convergence in expectation may not be a good indicator of convergence for solving all problems, as there may exist a problem for which RAC-ADMM is not stable or possesses greater variance.
In order to give another probabilistic measure on performance of RAC-ADMM, this section introduces convergence almost surely (a.s.) as an indicator of the algorithm convergence. 
Convergence almost surely as a measure for stability has been used in linear control systems for quite some time, and is based on the mean-square stability criterion for stochastically varying systems \cite{costa:2006}. The criterion establishes conditions for asymptotic convergence of covariance of the system states (e.g. variables). 

This section builds on those results and establishes 
sufficient condition for RAC-ADMM to converge almost surely when applied to solve (\ref{eq:problem_model_qp}). 
The condition utilizes the Kronecker product of the mapping matrix, which captures the dynamics of the second moments of the
random sequences generated by RAC-ADMM algorithm, and the expectation over the products of mapping matrices that provides 
the bounds on the variance of the distance between the KKT point and the random sequence generated by our algorithm.

\vspace{0.4cm}
\myth{\label{th:almost_sure_conv}
Suppose that Assumption \ref{assumption1} holds, and that RAC-ADMM (\ref{RAC-ADMM}) is employed to solve problem (\ref{eq:problem_model_qp}). Then the output of RAC-ADMM converges almost surely to some KKT point of (\ref{eq:problem_model_qp}) if
\myeqln{\rho(\E(M_\sigma\otimes M_\sigma))<1}
where $M\otimes M$ is the Kronecker product of $M$ with itself.
}

\vspace{0.4cm}
\myproof{
Let $\overline{z}=[\overline{x};\overline{y}]\in\R^{N}$ denote the KKT point of (\ref{eq:problem_model_qp}), then, at $k+1^{th}$ iteration we have 
\myeqln{
(z_{k+1}-\overline{z})=M_{\sigma_k}(z_k-\overline{z}).
}
Define $d_k=z_{k}-\overline{z}$, and 
\myeqln{
P_k=\E(d_kd_k^T).
}
There exists a linear operator $\mathcal{T}$ s.t.
\myeql{
\label{eq:kronecker}
\vect(P_{k+1})=\mathcal{T}\vect(P_k)
}
where $\vect(\cdot)$ is vectorization of a matrix, and $\mathcal{T}=\E(M_\sigma \otimes M_\sigma)$, as 
\myeqln{
\vect(P_{k+1})=\vect(\E(d_{k+1}d_{k+1}^T))=\dfrac{1}{|\Upsilon(n,p)|}\sum^{|\Upsilon(n,p)|}_{i=1} \vect(M_i\E(d_kd_k^T)M_i^T)=\E(M_\sigma \otimes M_\sigma) \vect(P_k)
}
and $\rho(\E(M_\sigma \otimes M_\sigma))<1$ implies $d_k\overset{a.s.}\to\bm{0}$.

To prove this, let $||\cdot ||$ be the Frobenius norm of a matrix, $||A||=\sqrt{\sum^m_{i=1}\sum^{n}_{j=1}|a_{ij}|^2}$
\myeqln{
\E(||d_k||^2)=tr(P_k)\leq ||\vect(P_k)||^2
}
And by $(\ref{eq:kronecker})$,
\myeqln{
||\vect(P_k)||^2=||\mathcal{T}\vect(P_{k-1})||^2=||\mathcal{T}^{k}\vect(P_0)||^2\leq ||\mathcal{T}^k||^2\cdot ||\vect(P_0)||^2
}
If $\rho(\mathcal{T})<1$, we know that $\mathcal{T}$ is convergent, and there exists $\mu>0$, $0<\gamma<1$, s.t. 
\myeqln{||\mathcal{T}^k||^2\leq \mu \gamma^k,} 
thus there exists $M$ such that, 
\myeqln{
\sum^{\infty}_{k=0}\E(||d_k||^2)\leq M\sum^{\infty}_{k=0}\gamma^{k} \leq C<\infty
}
For any $\epsilon>0$, by Markov inequality we have
\myeqln{
\sum^{\infty}_{k=0}\E(||d_k||^2)\leq C  \Rightarrow \sum^{\infty}_{k=0}\prob(||d_k||^2>\epsilon)<\infty,
}
and as $\sum^{\infty}_{k=0}\prob(||d_k||^2<\epsilon)<\infty$, by Borel-Cantelli, and  $||d_k||^2\in m\mathcal{F}_{+}$,
\myeqln{
d_k\overset{a.s.}\to\bm{0} \quad \textup{as }k\to\infty
}
which then implies that randomized ADMM converges almost surely. 
}

To illustrate the stability issues with RAC-ADMM, consider the following example.

\vspace{.3cm}
\example{
\label{example_diverge}
Consider the following problem 
\myeqln{
\begin{array}{ll}
 \max &\ \mathbf{0}\cdot \mathbf{x}\\
\hbox{s.t.} &\ A\mathbf{x}=\mathbf{0}
\end{array}}
where
\myeqln{
A=\begin{bmatrix}
    1 & 1 & 1 & 1 &1 &1 \\
    1 & 1 &1 &1 &1 &1+\gamma  \\
    1 &1 & 1 & 1 & 1+\gamma & 1+\gamma\\
    1 & 1 &1 &1+\gamma &1+\gamma &1+\gamma\\
    1& 1 & 1+\gamma &1+\gamma &1+\gamma &1+\gamma\\
    1& 1+\gamma & 1+\gamma &1+\gamma &1+\gamma &1+\gamma\\
\end{bmatrix}
}
Let $[x_0,y_0]\sim N(0,5I)$, $\beta=1$, $\gamma=1$, and number of blocks $p=3$. Consider RP-ADMM with the fixed block composition  $[x_1,\ x_2],[x_3,\ x_4],[x_5,\ x_6]$.

Convergence in expectation for this particular block structure finds $\rho(\E(M_{RP,\upsilon_1}))=0.9887>\rho(\E(M_{RAC}))=0.8215$. In fact, for all block compositions for this example we have, $\rho(\E(M_{RAC}))>\rho(\E(M_{RP,\upsilon_i})$ . However, RAC-ADMM does not converge, as shown in Figure \ref{fig:div}, showing that that convergence in expectation may not be a sufficient indicator from this particular example.

Indeed, if we apply Theorem \ref{th:almost_sure_conv}, we find out that RAC-ADMM does not converge almost surely, while RP-ADMM does for solving this example. Namely, $\rho(\E_{RAC}(M_\sigma\otimes M_{\sigma}))=1.0948>1$, and  $\rho(\E_{RP}(M_\sigma\otimes M_{\sigma}))=0.9852<1$, what explains the results shown in Figure \ref{fig:div}. In fact, RP-ADMM converges almost surely for all $15$ block compositions of this example.
\begin{figure}[t!]
    \centering
    \begin{subfigure}[t]{0.5\textwidth}
        \centering
        \includegraphics[height=2.1in]{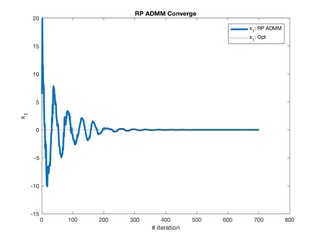}
        \caption{RP-ADMM}
    \end{subfigure}%
    \begin{subfigure}[t]{0.5\textwidth}
        \centering
        \includegraphics[height=2.1in]{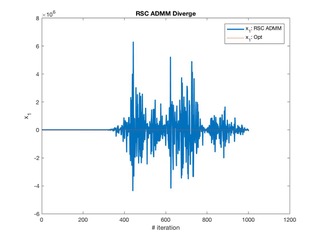}
        \caption{RAC-ADMM}
    \end{subfigure}
    \caption{Stability of RAC-ADMM and RP-ADMM}
\label{fig:div}
\end{figure}
}

%% file: RAC_Grouping.tex
\subsection{Variance Reduction in RAC-ADMM}
\label{subsec:grp}

The previous section described sufficient condition for the almost sure convergence of RAC-ADMM algorithm. This section address controlability of the algorithm. More precisely we ask,
given a linearly constrained quadratic problem (LCQP) (Eq. \ref{eq:problem_model_QP}), what means do we have at our disposal to control convergence of a LCQP-- how to bound the covariance and how to improve the convergence rate.

\subsubsection{Detecting and Utilizing a Structure in LCQP}
\label{sect:detect_structure}
Although some problem types inherit a known structure (e.g. network-flow problems), in general, the structure is not known. There are many sophisticated techniques used to detect a structure of a matrix   one can use and apply towards improving performance of RAC-ADMM. Although such elaborate methods have a potential of detecting hidden structure of Hessian and Jacobian matrices almost perfectly, using them or developing our own is beyond the scope of this paper. Instead, we adopt a simple matrix partitioning approach outlined in \cite{ferris:1998}.

In general, for RAC-ADMM we are interested in a structure of a constraint matrix, which can be detected using the following simple approach. Given a constraint matrix $A$ (describing equalities, inequalities or both), a desirable structure such as one shown in (\ref{block:A}) can be derived by applying a graph partitioning method.
\myeql{\label{block:A}
\begin{array}{ccc}
\underbrace{\begin{bmatrix}
    V_1 &  0 &\cdots &0 \\
    0 &  \ddots && \vdots\\
     \vdots& &  V_v&0 \\
    R_1&\cdots&R_v&R_{v+1}\\
\end{bmatrix}}_{\textstyle A} &
\underbrace{\begin{bmatrix}
\x_1\\
\vdots\\
\x_v\\
\x_{v+1}
\end{bmatrix}}_{\textstyle \x} &
=\underbrace{\begin{bmatrix}
\bb_1 \\
\vdots\\
\bb_v\\
\bb_{v+1}
\end{bmatrix}}_{\textstyle \bb}
\end{array}.
}
The outline of the process is as follows:
\begin{enumerate}
 \item Build a graph representation of matrix $A$: Each row $i$ and column $j$ is a vertex; vertices are connected with edges if $a_{i,j} \not= 0$.
 \item Partition the graph using a graph partitioning algorithm or solver, for example \cite{karypis_metis:1998}.
 \item Recreate $A$ as a block matrix from the graph partitions.
\end{enumerate}

\subsubsection{Smart Grouping}
\label{sect:smart_grouping}
Smart-grouping is a pre-processing method in which we use block structure of constraint matrix $A$ to pre-group certain variables as a single “super-variable” (a group of variables which always stay together in one block). Following the block structure shown in (\ref{block:A}), we make one super-variable $\xh_i$ for each group $\x_i$, $i=1,\dots,v$. Primal variables $\x_{v+1}$ stay shared and are randomly assigned to sub-problems to complement super-variables to which they are coupled with via block-matrices $R_i$, $i=1,\dots,v$. More than one super-variable can be assigned to a single sub-problem, dependent upon the maximum size of a subproblem, if defined.
Note that matrix partitioning based on $H+A^TA$ may result in a  better grouping, but is unpractical and thus not considered as a viable approach. 

\subsubsection{Partial Lagrangian}
\label{sect:partLan}
The idea of smart-grouping described in the previous section can be further extended by the means of the partial Lagrangian approach.
Consider a LCQP (\ref{eiq:problem_model_gen}) having the constraint matrix $A$ structure as shown in (\ref{block:A}). Now consider the scenario in which we split the matrix $A$ such that  the block $R\!=\![R_1,\dots,R_{v+1}]$ is admitted by the augmented Lagrangian while the rest of the constraints (blocks $V_i$) are solved exactly as a part of a sub-problem, i.e. a sub-problem $i$ is solved as
\myeqln{\x_i^{k+1}=\arg\min\{L_{\mathcal{P}} (\x_1^{k+1},\dots,\x_i,\dots,\x_p^k;\y^k)\cond V_j\xh_j=\bb_j, j\in\mathcal{J},\ \x_i\in\X_i\},}
where $\mathcal{J}$ is a set of indices of super-variables $\xh_j$ constituting sub-problem $i$ at any given iteration. The partial augmented Lagrangian is defined with
\myeqln{\label{lcqp:lp}L_{\mathcal{P}}(\x,\y)=\half\x^T H \x\ +\ c^T\x-\y^T(R\x-\bb_{v+1})+\frac{\beta}{2}\|R\x-\bb_{v+1}\|^2.}

 There are two advantages of the partial Lagrangian approach.
First, the rank of the constraint matrix used for the {\it global} constraints (matrix $R$) is lower than the rank of $A$, and the empirical results (Section \ref{sect:num}) suggest a strong correlation between a rank of a matrix and the stability of the algorithm and its rate of convergence .
Next, {\it local} constraints (matrices $V_i$) imply there is a feasibility region in which $\x_i$ exist, and that region may not be infinite. In other words, even when the variables themselves are unbounded (i.e.  $\x\in\R^n$), local constraints may put implicit bounds on maximum variation of values of $\x_i$.

Empirical results of the partial Lagrangian applied to mixed integer problems (Section \ref{subsect:num:bin}) show  the approach to be very useful. In such a scenario, local constraints are sets of rules that relate integer variables, while constraints between continuous variables are left global. In the case of a problems where such straight separation does not exist, or when problems are purely integer, a problem structure is let to guide the local/global constraints decision.

Although shown to be useful, the partial Lagrangian method suffers form being a mostly heuristic approach that depends on quality of solution methods applied to sub-problems -- in the case of continuous problems, a simple barrier based methodology can be applied, but for the mixed integers problems (MIP), sub-problems require a more complex solution (e.g. an external MIP solver). 

\vspace{10pt}
\example{To illustrate the usefulness of the smart grouping and partial Lagrange approaches, consider the following experiments done on selected instances taken from the Mittelmann LP test set \cite{bench:lp_test_set} and affiliated repositories augmented with diagonal Hessian $H$ to form a standard LCQP (\ref{eq:problem_model_QP}).

For each instance, a constraint matrix ($A_{eq}$, $A_{ineq}$ or $A=[A_{eq};A_{ineq}]$) was subjected to graph-partitioning procedure outlined in Section \ref{sect:detect_structure}, and then solved using the smart grouping (``s\_grp'') and the partial Lagrangian approach (``partial\_L''). 
Table \ref{tbl:pL} reports on the number of iterations required by RAC-ADMM algorithm to find a solution satisfying the primal/dual residual tolerance of $\epsilon=10^{-4}$. If the solution was not found the reason is noted (``time limit'' for exceeding sub-problem maximum run-time and ``iter. limit'' for exceeding maximum number of iterations). Fields showing ``divergence'' or ``oscillation'' mark experiments for which RAC-MBADMM algorithm experienced an unstable behavior. The baseline for the comparison is the {\it default} approach (sub-problems created at random) shown in column ``Default RAC''.   

\begin{table}[h!]\label{tbl:grouping}
\footnotesize
\centering
\begin{tabular}{lccccccccccc}  
\toprule
Instance & Num. & Num.  & Default&\multicolumn{2}{c}{Partitioning on $A_{eq}$} &\multicolumn{2}{c}{Partitioning on $A_{ineq}$} &\multicolumn{2}{c}{Partitioning on $A$} \\ 
 & col. & rows (A) & RAC & partial\_L & s\_grp & partial\_L & s\_grp & partial\_L & s\_grp\\ \midrule
qap10 & 4150 & 1820 & 182 & 102 & 550 &  &  &  & \\            
pds-02 & 7535 & 2953 & 2872 & 577 & 1530 & iter. limit & 798 & 1287 & 3878\\ 
n370a & 10000 & 5150 & 4797 & time limit & 852 & 211 & 1819 & 148 & 1298\\ 
supportcase10 & 14770 & 165684 & 5535 & 240 & 349 & 1370 & 2622 & 1480 & 2496\\ 
ex10 & 17680 & 69608 & 324 & 281 & 630 & oscillation & 410 & oscillation & 410\\ 
nug08-3rd & 20448 & 19728 & 963 & 437 & 414 & oscillation & 516 & 624 & 624\\   
brazil3 & 23968 & 14646 & 1083 & divergence & 1103 & oscillation & 2207 & divergence & 1126\\\bottomrule
\end{tabular}
\caption{Number of iterations until termination criteria is met for various benchmark instances\vspace{0pt}}\label{tbl:pL}
\end{table}

The partial Lagrangian approach has a potential to help stability and rate of convergence. However, before generalizing, one needs to consider the following: stability (i.e. convergence) of RAC-ADMM algorithm,
is a function, among other factors, of mapping operators (matrices $M_{\sigma}$, Eq. \ref{eq:mapping_Mi}) which are in turn functions, among other factors, of the constraint matrix of a problem being solved. In the case of partial Lagrangian methodology, this matrix is the matrix $R$,  meaning that if $R$ produces an unstable system (e.g. conditions set by Theorem \ref{th:almost_sure_conv} not met), no implicit bounding can help to stabilize it. 

Using smart grouping alone, on the other hand, does not make RAC-ADMM unstable, but in some cases increases the number of iterations needed to satisfy feasibility tolerance, a consequence of having less randomness as described by Corollary \ref{col:rac_speed}.

%% file: RACQP.tex
\section{RAC-ADMM Quadratic Programming Solver}
\label{sect:solver}

In this section we outline the implementation of the RAC-ADMM algorithm for linearly constrained quadratic problems as defined below:
 \myeqmodel{\label{eq:problem_model_QP}
 \min\limits_{\x} & f(\x)=\half\x^T H \x + \cc^T\x\\[0.3cm]
 \mbox{s.t.}& A_{eq}\x =\bb_{eq}\\[0.2cm]
    & A_{ineq}\x \le\bb_{ineq}\\[0.2cm]
    & \x\in\X
 }
where symmetric positive semidefinite matrix $H\in\R^{n\times n}$ and vector $\cc\in\R^n$ define the quadratic objective while matrix $A_{eq}\in\R^{m\times n}$ and the vector $\bb_{eq}\in\R^m$ describe equality constraints and matrix $A_{ineq}\in\R^{s\times n}$ and the vector $\bb_{ineq}\in\R^s$ describe inequality constraints. Primal variables $\x\in\X$, can be integral or continuous, thus the constraint set $\X$ is the Cartesian product of nonempty sets $\X_i\subseteq \R$ or $\X_i\subseteq \Z$, $i=1,\dots,n$. QP problems arise from many important applications themselves, and are also fundamental in general nonlinear optimization. 

Reformulate (\ref{eq:problem_model_QP}) as follows:
\myeqmodel{\label{eq:problem_model_QPref}
 \min\limits_{\x,\ {\tilde \x},\ \s} & \half \x^T H \x + \cc^T\x\\[0.3cm]
 \mbox{s.t.}& A_{eq}\x =\bb_{eq}\\[0.2cm]
    & A_{ineq}\x +\s = \bb_{ineq}\\[0.2cm]
    & \x-{\tilde \x} = \bz \\[0.2cm]
    & {\tilde \x}\in\X,\ \s\ge\bz, \x \hbox{free}
 }
where the augmented Lagrangian is given as
\myeql{ \label{eq:argL_racqp} 
\begin{array}{cl} 
L_{\beta} (\x\in \X;{\tilde\x};\s;\y_{eq};\y_{ineq};\z) := &\half\x^T H \x + \cc^T\x
-\y_{eq}^T(A_{eq}\x -\bb_{eq})
-\y_{ineq}^T(A_{ineq}\x +\s -\bb_{ineq})
-\z^T(\x-{\tilde \x})\\[0.2cm]
&+ \frac{\beta}{2} (\|A_{eq}\x -\bb_{eq}\|^2
+ \|A_{ineq}\x +\s -\bb_{ineq}\|^2+\|\x-{\tilde \x}\|^2)
\end{array}
}
RAC-ADMM, or simply RAC, quadratic programming (RACQP) solver admits continuous, binary and mixed integer problems. Algorithm \ref{alg:RACQP} outlines the solver: the solution vector is initialized to $-\infty$ at the beginning of the algorithm, and the main RAC-ADMM loop described (lines \ref{alg:RACQP:iter_start}-\ref{alg:RACQP:iter_end}). The main loop calls different procedures to optimize blocks of $\x$ (lines  \ref{alg:RACQP:block_start}-\ref{alg:RACQP:block_end}), followed by updates of slack and then dual variables. 
\input{Algorithm_RACQP}

Types of the block optimizing procedure being called to update the blocks depend on the structure of the problem being solved. The default, multi-block implementation for continuous problems is based on the Cholesky factorization, with a specialized one-block variant for very sparse problems that solves the iterates using the LDL factorization. 
Continuous problems that exhibit a structure (see Section \ref{sect:detect_structure}) can be addressed using the partial Lagrangian approach. In such a case, sub-problems are solved using either a simple interior point method based methodology, or, when sub-problems include only equality constraints, by employing Cholesky for solving KKT conditions.
In addition to the aforementioned methods, the solver supports calls to external solver(s) and specialized heuristic solution to handle hard sub-problem instances.

Binary and mixed integer problems require specialized optimization techniques (e.g. branch-and-bound), that require implementations which are beyond the scope of this paper, so we have decided to delegate optimizing of the blocks with mixed variables to an external solver. Mixed integer problems are addressed by using the partial Lagrangian to solve for primal variables and a simple procedure that helps to escape local optima, as described by Algorithm \ref{alg:RACQP:mip}.

Note that Algorithms given in this section are pseudo-algorithms which describe functionality of the solver rather than actual implementation. The implementation can be downloaded from \cite{RACQP:code}.

\subsection{Solving continuous problems}
\label{subsect:solver:cont}
For the continuous QP problems, we consider (\ref{eq:problem_model_QPref}) where $\X$ are possible simple lower and upper bounds on
each individual variable: 
\[l_i\le{\tilde x_i}\le u_i,\ i=1,\dots,n.\]
Continuous problems are solved as described by Algorithm \ref{alg:RACQP}, which repeats three steps until termination criteria is met: first update or optimize primal variables $\x$ in the RAC fashion, then update ${\tilde \x}$ and $\s$ in close forms and finally update dual variables $\y_{eq}$, $\y_{ineq}$ and $\z$.

\underline{Step 1}: Update primal variables $\x$

Let $\omega_i\in\Omega$ be a vector of indices of a block $i$, $i=1,\dots,p$, where $p$ is the number of blocks. The set of vectors $\Omega$ is randomly generated (with smart grouping when applicable as described in Section \ref{sect:smart_grouping}) at each iteration of the Algorithm \ref{alg:RACQP} (lines \ref{alg:RACQP:iter_start}-\ref{alg:RACQP:iter_end}). Let $\x_{\omega_i}$ be a sub-vector of $\x$ constructed of components of $\x$ with indices $\omega_i$, and let $\x_{-\omega_i}$ be the sub-vector of $\x$ with indices not chosen by $\omega_i$. 
Algorithm \ref{alg:RACQP} uses either Cholesky factorization or partial Lagrangian to solve each block of variables $\x_{\omega_i}$ while holding $\x_{-\omega_i}$ fixed. 

By rewriting (\ref{eq:argL_racqp}) to reflect the sub-vectors, we get

\myeql{ \label{eq:Qq}
\begin{array}{cl} 
L_\beta(\cdot)&=[\x_{\omega_i};\x_{-\omega_i}]^T(\half H+\bhalf A_{eq}^TA_{eq}+\bhalf A_{ineq}^TA_{ineq}+\bhalf I)[\x_{\omega_i};\x_{-\omega_i}]\\[0.2cm]
&\hspace{8pt}+(\cc-A_{eq}^T\y_{eq}-A_{ineq}^T\y_{ineq}-\z-\beta A_{eq}^T\bb_{eq}-\beta A_{ineq}^T\bb_{ineq} +\beta A_{ineq}^T\s-\beta{\tilde \x})^T[\x_{\omega_i};\x_{-\omega_i}]\\[0.2cm]
&=\frac{1}{2}[\x_{\omega_i};\x_{-\omega_i}]^TQ[\x_{\omega_i};\x_{-\omega_i}]+\q^T[\x_{\omega_i};\x_{-\omega_i}]
\end{array}
}
where \myeq{Q=(H+\beta A_{eq}^TA_{eq}+\beta A_{ineq}^TA_{ineq}+\beta I)}. Then
we can minimize in $\x_{\omega_i}$ by solving 
\myeq{\label{eq:racqp_multi}
Q_{\omega_i,\omega_i}\x_{\omega_i}=-(\q_{\omega_i}+{\hat \q})} 
using Cholesky factorization and back substitution. The linear term resulting from $Q$, ${\hat \q}$, is given as 

\vspace{-.4cm}
\myeql{\label{eq:q_hat}
{\hat\q}= (H{\hat \x})_{\omega_i}+\beta A_{\omega_i}^T(A{\hat \x})-
                          (H_{\omega_i,\omega_i}+\beta A_{\omega_i}^TA_{\omega_i})\x_{\omega_i}}

where \myeq{A=[A_{eq},{\mathit 0};A_{ineq},I]} and ${\hat \x}=[\x;\s]$. A square sub-matrix $H_{\omega_i,\omega_i}$ and column sub-matrix $A_{\omega_i}$ are constructed by extracting $\omega_i$ rows and columns from $H$ and $A$ respectively. 

When $p=1$, i.e. we are solving a problem using a single-block approach, then we solve the block utilizing LDL factorization to avoid calculating $A^TA$. Although the factorization can be relatively expensive if the problem size is large as we then factorize a large matrix, the factorization is done only once and re-used in each iteration of the algorithm. From (\ref{eq:argL_racqp}), we find minimizer $\x$ by solving 

\myeql{\label{eq:LDL}
Q\x=-(\cc-A_{eq}^T\y_{eq}-A_{ineq}^T\y_{ineq}-\z-\beta A_{eq}^T\bb_{eq}-\beta A_{ineq}^T\bb_{ineq} +\beta A_{ineq}^T\s-\beta{\tilde \x})=-\q}

With 
 \myeq{A=[A_{eq};A_{ineq}]} we can express the equivalent condition to (\ref{eq:LDL}) with

\myeql{\label{eq:single_p}
\begin{array}{ccc}
\begin{bmatrix}(H+\beta I) & \sqrt\beta A^T\\
\sqrt\beta A &  -I\end{bmatrix} & 
\begin{bmatrix}\x\\
\mu\end{bmatrix} &
=\begin{bmatrix}-\q \\
\bz\end{bmatrix} 
\end{array}.
}

We factorize the left hand side of the above expression and use the resulting matrices to find $\x$ by back substitution at each iteration of the algorithm. For single-block RACQP, LDL approach described above replaces lines \ref{alg:RACQP:m_start}-\ref{alg:RACQP:block_end} in Algorithm \ref{alg:RACQP}.

Furthermore, if \myeq{H} is diagonal, one can rewrite the system as
\myeql{\label{eq:single_p_diag}
\begin{array}{ccc}
\begin{bmatrix}I & (H+\beta I)^{-1}\sqrt\beta A^T\\
\sqrt\beta A & -I\end{bmatrix} & 
\begin{bmatrix}\x\\
\mu\end{bmatrix} &
=\begin{bmatrix}-\q \\
\bz\end{bmatrix} 
\end{array}.
}
Then we can factorize matrix \myeq{(I+\beta A (H+\beta I)^{-1}A^T)} to solve the system, which would be 
extremely effective when the number of constraints is very small and/or sparse, since \myeq{(H+\beta I)^{-1}} is diagonal and it does not change sparsity of \myeq{A}.

Partial Lagrangian approach to solving $\x$ blocks, described in Section \ref{sect:partLan}, uses the same implementation as Cholesky approach described above, with additional steps that build local constraints which reflect free and fixed components of $\x$, $\x_{\omega_i}$ and $\x_{-\omega_i}$ respectively. The optimization problem of partial Lagrangian is formulated as
\myeql{\label{eq:partL}
\begin{array}{lrl}
\x_{\omega_i}^* =& \argmin&\half\x_{\omega_i}^T\!{ Q}\x_{\omega_i}+(\q_{\omega_i}+{\hat \q}^T)\x_{\omega_i}\\[0.2cm]
& \hbox{s.t. } &\Ah_{eq,\ \omega_i}\x_{\omega_i}=\bh_{eq}-\Ah_{eq,\ -\omega_i}\x_{-\omega_i}\\[0.2cm]
&&\Ah_{ineq,\ \omega_i}\x_{\omega_i}\le\bh_{eq}-\Ah_{ineq,\ -\omega_i}\x_{-\omega_i}\\[0.2cm]
&&\lb_{\omega_i}\le\x_{\omega_i}\le\ub_{\omega_i}
\end{array}
}

with $\Ah_{eq},\ \bh_{eq}$ and $\Ah_{ineq},\ \bh_{ineq}$
describing local equality and inequality constraints, respectively.

Note that partial Lagrangian procedure is used by both continuous and mixed integer problems. In the case of the former we set $\X=\R^n$, while when we solve the latter we let $\X_i\subseteq\R$ and implicitly enforce the bounds. The blocks are solved by either an external solver (e.g.Gurobi) or by using Cholesky to solve KKT conditions when $\x_{\omega_i}$ is unbounded.

\underline{Step 2}: Update auxiliary variables ${\tilde \x}$

With all variables but ${\tilde \x}$ fixed, from augmented Lagrangian (\ref{eq:argL_racqp}) we find that the optimal vector $\lb\le {\tilde \x}\le \ub$ can be found by solving the optimization problem
\myeqln{\argmin\limits_{\lb\le {\tilde \x}\le \ub}\ \bhalf{\tilde \x}^T{\tilde \x}+(\z-\beta\x^T){\tilde \x}.}
The problem is separable and ${\tilde \x}$ has a closed form solution given by
\myeqln{{\tilde \x}=\min\big\{\max\{\lb,\x-\frac{1}{\beta}\z\},\ub\big\}}

\underline{Step 3}: Update slack variables $\s$

Similarly to the previous step, with all variables but $\s$ fixed, the optimal vector $\s$ is found by solving 
\myeqln{\argmin\limits_{\s\ge 0}\ \bhalf\s^T\s+(-\y_{ineq}+\beta(A_{ineq}\x-\bb_{ineq}))^T\s.}
The problem is separable and $\s$ has a closed form solution given by
\myeqln{\s=\max\big\{0,\frac{1}{\beta}\y_{ineq}+\bb_{ineq}-A_{ineq}\x\big\}.}

\subsubsection{Termination Criteria for Continuous Problems}
\label{subsect:solver:cont:term}

Termination criteria for continuous problems include maximum run-time limit settings, maximum number of iterations and primal-dual solution 
(found up to some tolerance). RACQP terminates when at least one criterion is met. 
For primal-dual solution criterion RACQP uses the optimality conditions of problem (\ref{eq:problem_model_QPref}) to define primal and dual relative residuals at iteration $k$,

\myeql{
\begin{array}{ll}
 r_{\hbox{prim}}^k&:= \max( r_{A_{eq}}^k, r_{A_{ineq}}^k, r_{bounds}^k)\\[0.2cm]
r_{\hbox{dual}}^k &:= \ddfrac{\|H\x^k+c-A_{eq}^T\y_{eq}^k-A_{ineq}^T\y_{ineq}^k-\z\|_\infty}
{1+\max(\|H\x^k\|_\infty, \|c\|_\infty,\ \|A_{eq}^T\y_{eq}^k\|_\infty,\ \|A_{ineq}^T\y_{ineq}^k\|_\infty,\ \|\z\|_\infty)}
\end{array}
}
where

\myeqln{
\begin{array}{ll}
r_{A_{eq}}^k&= \ddfrac{\|A_{eq}\x^k-\bb_{eq}\|_\infty}{1+\max(\|A_{eq}\x^k\|_\infty,\ \|\bb_{eq}\|_\infty)}\\ [0.4cm]
r_{A_{ineq}}^k&=\ddfrac{\|A_{ineq}\x^k+s-\bb_{ineq}\|_\infty}{1+\max(\|A_{ineq}\x^k+\s^k\|_\infty,\ \|\bb_{ineq}\|_\infty)}\\ [0.4cm]
r_{bounds}^k&=\ddfrac{\|\x^k-{\tilde\x^k}\|_\infty}{1+\max(\|\x^k\|_\infty,\ \|{\tilde\x^k}\|_\infty)}
\end{array}
}

and set RACQP to terminate when the residuals become smaller than some tolerance level $\epsilon>0$.
\myeql{\max(r_{\hbox{p}}^k,\ r_{\hbox{d}}^k)<\epsilon.}
Note that the aforementioned residuals are similar to those used in \cite{boyd:2011,stellato:2018} with relative and absolute residual tolerance ($\epsilon_{abs},\epsilon_{rel}$) set to be equal.


\subsection{Mixed Integer Problems}

\input{Algorithm_RACQP_MIP}

For mixed integer problems we tackle (\ref{eq:problem_model_QP}) without introducing $\tilde{\x}$, where augmented Lagrangian is given by
\myeqln{ 
\begin{array}{cl} 
L_{\beta} (\x;\s;\y_{eq};\y_{ineq}) := &\half\x^T H \x + \cc^T\x
-\y_{eq}^T(A_{eq}\x -\bb_{eq})
-\y_{ineq}^T(A_{ineq}\x +\s -\bb_{ineq})\\[0.2cm]
&+ \frac{\beta}{2} (\|A_{eq}\x -\bb_{eq}\|^2
+ \|A_{ineq}\x +\s -\bb_{ineq}\|^2)
\end{array}
}
where slack variables $\s\ge0$, and $x_i\in\X_i$, $\X_i\subseteq \R$ or $\X_i\subseteq \Z$, $i=1,\dots,n$. Mixed integer problems (MIP) are addressed by using the partial Lagrangian to solve for primal variables and a simple procedure that helps to escape local optima, as shown in Algorithm \ref{alg:RACQP:mip}. Note that MIP and continuous problems share the same main algorithm (Algorithm \ref{alg:RACQP}), but the former ignores the update to ${\tilde \x}$ as the bounds on $\x$ are explicitly set through $\X$, and thus ${\tilde \x} = \x$ always.

RACQP-MIP Solver, outlined in Algorithm \ref{alg:RACQP:mip}, consists of a sequence of steps that work on improving the current (or initial) solution which is then ``destroyed`` to be possibly improved again.  This solve-perturb-solve sequence (lines \ref{alg:RACQP:mip_lstart}-\ref{alg:RACQP:mip_lend}) is repeated until termination criteria is met. The  criteria for RACQP-MIP is usually set to be maximum run-time, maximum number of attempts to find a better solution, or a solution quality (assuming primal feasibility is met within some $\epsilon>0$). The algorithm can be seen as a variant of a neighborhood search technique usually associated with meta-heuristic algorithms for combinatorial optimization. 

After being stuck at some local optimum solution, the algorithm finds a new initial point $\x_0$ by perturbing the best known solution $\x_{best}$ and continues from there. The new initial point does not need to be feasible, but in some cases it may be beneficial to be constructed that way. To detect a local optimum we use a simple approach that counts number of times a ``feasible'' solution is found without improvement in objective value. A solution is considered to be feasible if \myeq{\max(\|A_{eq}\x-\bb_{eq}\|_\infty,\ \|A_{ineq}\x-\bb_{ineq}\|_\infty)\le\epsilon},  $\epsilon>0$. 
Perturbation (line \ref{alg:RACQP:mip:perturb}) can done, for example by choosing a random number (chosen from a truncated exponential distribution) of components of $\x_{best}$ and assigning them new values, or a more sophisticated approach can be used (see Section \ref{subsect:num:bin} for some implementation details). Parameters of permutation are encapsulated in a generic term $\kappa$.

%% file: Algorithm_RACQP.tex
\begin{algorithm}[t]
\ALGORITHM
{RACQP \label{alg:RACQP}}
{\begin{algorithmic}[1] 
\vspace{5pt}
\Require  Problem {\it model} (Eq. \ref{eq:problem_model_QPref}), {\it run-time parameters}$^\dagger$ with {\it termination criteria}$^{\dagger\dagger}$
\Ensure The optimal solution $\x^*$ or the best solution found before termination criteria met 
\vspace{5pt} 
    \State \myeq{\x\gets -\infty}
    \While  {{\it termination criteria} not met}\label{alg:RACQP:iter_start}    
    \State $\Omega\gets$ construct blocks at random; use smart grouping if applicable(Section \ref{sect:smart_grouping})\label{alg:RACQP:m_start}
    \ForAll {vectors ${\omega_i}\in\Omega$ of block indices}\label{alg:RACQP:block_start}\Comment{Solve $\x_{\omega_i}$ blocks}
    \State Prepare $Q_{\omega_i,\omega_i}$, $\q_{\omega_i}$ and ${\hat \q}$ following equations (\ref{eq:Qq}) and (\ref{eq:q_hat})
    \State $\x_{\omega_i}\gets$ solve \myeq{Q_{\omega_i,\omega_i}\x_{\omega_i}=-(\q_{\omega_i}+{\hat \q})} using:
    \State \hspace{10pt}  {\bf if} (sub-problem is mixed integer) {\bf then}
     \State \hspace{23pt} An external solver \label{alg:RACQP:external}
     \State \hspace{10pt}  {\bf else if} (partial Lagrangian (\ref{eq:partL}) used and sub-problem includes inequalities) {\bf then} 
     \State \hspace{23pt} Interior point method based procedure or an external solver  \label{alg:RACQP:IPM}
    \State \hspace{10pt}  {\bf else if} (partial Lagrangian used) {\bf then}
     \State \hspace{23pt} Cholesky factorization and back substitution solving KKT conditions
    \State \hspace{10pt}  {\bf else} 
     \State \hspace{23pt} Cholesky factorization and back substitution
    \State \hspace{10pt}  {\bf end if} 
    \EndFor\label{alg:RACQP:block_end}
    \State \myeq{\s^*\gets\max(0,\frac{1}{\beta}\y_{ineq}+\bb_{ineq}-A_{ineq}\x)} \Comment{Update slack variables}
    \If {(bounds on $\x$ not addressed by the partial Lagrangian)}
    \State \myeq{{\tilde \x} \gets \min(\max(\lb,\x-\frac{1}{\beta}\z),\ub)}\Comment{Update auxiliary split variables ${\tilde \x}$}
    \State \myeq{\z\gets\z-\beta(\x-{\tilde \x})} \Comment{Update dual variables for split variables}
    \EndIf
    \State \myeq{\y_{eq}\gets\y_{eq}-\beta(A_{eq}\x-\bb_{eq})}\Comment{Update dual variables for equality constraints}
    \State \myeq{\y_{ineq}\gets\y_{ineq}-\beta(A_{ineq}\x-\bb_{ineq})}\Comment{Update dual variables for inequality constraints}
    \EndWhile  \label{alg:RACQP:iter_end}
    \State \Return $\x$  
\end{algorithmic}}
{\footnotesize \hspace{-8pt} $\dagger$ Number of groups $p$, penalty parameter $\beta$, initial point $\x_0$, pre-grouped vars set ${\mathcal V}$.\\
$\dagger\dagger$ Termination criteria may include maximum run-time, number of attempts to find a better solution, solution quality and so on. }

\end{algorithm}

%% file: Algorithm_RACQP_MIP.tex
\begin{algorithm}[t!]
\ALGORITHM
{RACQP-MIP \label{alg:RACQP:mip}}
{Block Optimization with Integer Variables \begin{algorithmic}[1] 
\vspace{5pt}
\Require  Problem {\it model} (Eq. \ref{eq:problem_model_QPref}), {\it run-time parameters}$^\dagger$ , {\it termination criteria}$^{\dagger\dagger}$
\Statex \hspace{0.8cm} Perturbation parameters $\kappa$, number of trials before perturbing $nP$
\Ensure The optimal solution $\x^*$ or the best solution found before termination criteria met  
\State \myeq{\x_{best}\gets -\infty,\ k\gets 0}
\While  {{\it termination criteria} not met}\label{alg:RACQP:mip_lstart}
\State \myeq{\x^*\gets}\Call{RACQP}{{\it model}, {\it run-time parameters}}
\If {\myeq{f(\x^*)<f(\x_{best})}}
\State \myeq{\x_{best} \gets \x^*}
\Else
\State \myeq{k\gets k+1}
\EndIf
\If{$k=nP$}
\State \myeq{k\gets 0}
\State \myeq{\x_0\gets} perturb ($\x_{best}$, $\kappa$) \label{alg:RACQP:mip:perturb}
\EndIf
\EndWhile \label{alg:RACQP:mip_lend}
\State \Return $\x_{best}$
\end{algorithmic}}
{\footnotesize \hspace{-8pt} $\dagger$ RACQP run-time parameters (number of groups $p$, penalty parameter $\beta$, initial point $\x_0$, pre-grouped vars set ${\mathcal V}$, termination criteria).\\
$\dagger\dagger$ RACQP-MIP termination criteria (e.g. maximum run-time, number of attempts to find a better solution, solution quality and so on). }
\end{algorithm}

%% file: Empirical_Analysis.tex
\section{Computational Studies}
\label{sect:num}

The Alternating Direction Method of Multipliers (ADMM) has nowadays gained a lot of attention for solving many problems of practical importance (e.g. large-scale machine learning and signal processing, image processing, portfolio-management, to name a few). Unfortunately, the two most popular approaches, namely the two-block classical ADMM  and the variable-splitting multi-block \cite{boyd:2011}, both characterized by convergence speed and scaling issues  somehow hindered a wide acceptance of ADMM as the solution method of choice for ML problems. RAC-ADMM offers the multi-block solution that may help to overcome the problem of ADMM acceptance. 

The goal of this section is twofold: (1) to show that RAC-ADMM is a versatile algorithm that can be directly applied to a wide range of LCQP problems and compete with commercial solvers and (2) get an insight on  specific ML problems and devise a RAC-ADMM based solution that outperforms or matches the performance of the best tailored solution method(s) in both solution time and quality. To address the former, in Sections  \ref{subsect:num:cont} and \ref{subsect:num:bin} we compare RACQP with the state of the art commercial solvers, Gurobi \cite{gurobi} and Mosek \cite{mosek}, and the academic OSQP which is a ADMM-based solver developed by \cite{stellato:2018}. To address the latter, we focus on Linear
Regression (Elastic-Net) and Support Vector Machine (SVM),  machine learning algorithm used for classification and regression analysis, and in Section \ref{sect:ML} compare RACQP with glmnet \cite{friedman:2010, friedman:2011} and LIBSVM \cite{Chang:2011}.  

We conduct multiple numerical tests, solving randomly constructed problems and problems from benchmark test sets.
Data we collect include run-time, number of iterations until termination criteria is met  and quality of a solution, defined differently for continuous, mixed-integer and machine learning problems (described in corresponding subsections). Note that in some sections, due to space concerns we report on a subset of instances. Experiments using larger sets are  available together with RACQP solver code online \cite{RACQP:code} in ``demo'' directory.

The experiments were done on MacBook Pro with 2.8 GHZ Intel Core i7 and 16Gb memory running macOS High Sierra, v 10.13.2 (Section \ref{sect:lasso}) and 16-core Intel Xeon CPU E5-2650 machine with 96Gb memory running Debian linux 3.16.0-4-amd64 (all other sections). 

\input{Empirical_Analysis_Subsec_Cont}

\input{Empirical_Analysis_Subsec_Bin}

\input{Empirical_Analysis_Subsec_ML}

%% file: Empirical_Analysis_Subsec_Cont.tex
\subsection{Continuous problems}
\label{subsect:num:cont}
The section starts with the analysis of the regularized ($l_2$-norm) Markowitz mean-variance model applied to 2018 CSRP Quarterly Stock data \cite{WRDS} followed by randomly generated convex quadratic problems (QP) with coupled blocks. Next three sets of benchmark problems are addressed: relaxed QAPLIB \cite{bench:qaplib} (binary constraint on variables removed),  Maros and Meszaros
Convex QP \cite{bench:cute_r},  and the Mittelmann LP test set \cite{bench:lp_test_set} expanded to QP by adding  a diagonal Hessian to the problem model.

The goal of the section is to show that the multi-block ADMM approach adapted by RACQP can significantly reduce solution time compared to commercial solvers and two-block ADMM (used by OSQP) for most of the problems we addressed. Results obtained in this section are all done with a single RACQP run, using fixed random number generator seed. Performance of the solver when subjected to different seeds is described in Section \ref{sect:num:rnd_seed}.

The run-time settings applied to solvers to produce results reported in this section, unless noted otherwise, are shown in Table \ref{tbl:exp_defaults}. 

\begin{table}[h!]
\footnotesize
\centering
\setlength\tabcolsep{5pt}
\begin{tabular}{ll}  
\toprule
 Termination & Parameter \\
 criteria & value\\
\midrule
Max time & 3h \\
Max. num. iterations (OSQP, RACQP)& 4000\\
Primal residual (feasibility) tolerance & $\epsilon_{prim}=\epsilon\eq10^{-5}$ \\
Dual residual (optimality) tolerance &$\epsilon_{dual}=\epsilon\eq10^{-5}$ \\
Relative residual tolerance (OSQP) &  $\epsilon_{rel}=\epsilon\eq10^{-5}$ \\
Barrier convergence tolerance (Gurobi, Mosek) &$\epsilon_{comp}=\epsilon\eq10^{-5}$ \\
\bottomrule
\end{tabular}
\caption{Termination criteria used in this section by all solvers.}\label{tbl:exp_defaults}
\end{table}

 Authors are aware that either commercial solver can be tuned for maximum performance by adjusting run-time parameters to fit a specific problem structure, which is the same with RACQP and OSQP but to the much smaller extent. In addition, the latter do not have the access to a large number of real-world instances used by the former to fine-tune algorithms to exploit ``known'' problem structures nor manpower to build heuristics and/or preconditioners that boost solver performance. However, in order to create a more ``equal'' working conditions, we decided to let Mosek and Gurobi use their default settings, except for disabling multi-threading support and aforementioned optimality termination criteria (Table \ref{tbl:exp_defaults}).  Although allowing the solvers to execute presolve routines seems to be unfair to RACQP (which does not implement any presolving technique except for a very simple row scaling), disabling it would be even more unfair to the opposing solvers as their performance heavily depends on finesses of the presolve algorithm(s). Multi-threading is disabled for Mosek and Gurobi because both RACQP and OSQP are single-threaded, and leaving it on would be unfair. Finally, to make RACQP and OSQP comparison more fair, and because our target is to compare two ADMM variants, RAC-ADMM and operator splitting two-block ADMM, rather than solvers' implementations, the advanced option that OSQP uses to post-process results, ``Polish results'', was turned off. Note that such an option is relatively easy to implement and a variant of thereof will be added to a future RACQP version.

For continuous problems described in this section, performance is measured in terms of run-time, number of iterations and quality of solution, expressed via primal and dual residuals. Terminating a run after residual(s) have been met (Table \ref{tbl:exp_defaults}, rows 2-4) is one way of ensuring  quality of a solution. 
However, this criteria could be misleading. To start with, some solvers use absolute residuals as termination criteria (e.g. Gurobi), some depend on relative residuals (e.g. Mosek, RACQP), and some are adjustable like QSQP.

Next, solvers usually scale problems (e.g. row and column scaling of a constraint matrix) to avoid numerical problems and make matrices with favorable condition numbers. Residuals are then calculated and checked against these scaled models, meaning that a solver may prematurely terminate unless the results are periodically re-scaled and residuals recalculated on the actual model -- a ``good'' scaled solution can actually have a very bad ``actual'' residual. As each solver performs different scaling (and algorithms are not usually known as it is case with Gurobi and Mosek), direct comparison of residuals reported by the solvers is not possible.

To circumvent the issue, we re-calculate primal and dual residuals using the solutions (primal and dual variables), returned by the solvers as follows:

\myeqln{
\begin{array}{ll}
 r_{\hbox{prim}}^k&:= \max( r_{A_{eq}}^k, r_{A_{ineq}}^k, r_{bounds}^k)\\[0.2cm]
r_{\hbox{dual}} &:= \ddfrac{\|H\x^*+c-A^T\y^*-\y_{bounds}^*\|_\infty}
{1+\max(\|H\x^k\|_\infty, \|c\|_\infty,\ \|A^T\y^*\|_\infty,\ \|\y_{bounds}^*|_\infty)}
\end{array}
}
where $A=A_{eq};A_{ineq}]$, $y^*$ is a vector of dual variables related to equality and inequality constraints,  $y_{bounds}^*$ is a vector of dual variables related to primal variable bounds, and $\x^*$ is a vector of primal variables. 
Residuals due to equality and inequality constraints and bounds are defined with

\myeqln{
\begin{array}{ll}
r_{A_{eq}}&= \ddfrac{\|A_{eq}\x^*-\bb_{eq}\|_\infty}{1+\max(\|A_{eq}\x^*\|_\infty,\ \|\bb_{eq}\|_\infty)}\\ [0.4cm]
r_{A_{ineq}}&= \ddfrac{max(0, \|A_{ineq}\x^*-\bb_{ineq}\|_\infty)}{1+\max(\|A_{ineq}\x^*\|_\infty,\ \|\bb_{ineq}\|_\infty)}\\ [0.4cm]
r_{bounds}&=\max(\ddfrac{\|\max(0,\lb-\x^*)\|_\infty}{1+\max(\|\x\|_\infty,\ \|{\ \lb\ }\|_\infty)}, 
\ddfrac{\|\max(0,\x^*-\ub)\|_\infty}{1+\max(\|\x\|_\infty,\ \|\ {\ub}\ \|_\infty)})
\end{array}
}

Note that Gurobi does not provide dual variables for bounds ($\lb\le\x\le\ub$) directly. To get around we convert the bounds into inequality constraints, what makes Gurobi to produce the dual variables. This introduces negligible run-time cost as the additional constraints are discovered as bounds during presolve phase and consequently removed. The initial point $\x^0$ for all instances addressed by RACQP is $\max(\bz,lb)$.

\subsubsection{Choosing RACQP solver working mode}
\label{sect:racqp_rules}
To address differences in problem structure, the following simple rules are used to decide on the RACQP solver mode:
\begin{enumerate}
 \item If $H$ is non-diagonal and $A$ is non-structural or the problem is large, use multi-block mode (Eq. \ref{eq:racqp_multi}) \label{rule:1}
 
 \item  If $H$ is non-diagonal and $A$ is structural, which implies that $A$ has non-zero entries that follow some pattern  and problem structure is easy to detect, use multi-block mode with smart-grouping as described in Section \ref{sect:smart_grouping}.

\item If $H$ is diagonal, $m<<n$ or $H$ and $A$ are very sparse, and the problem is of moderate size, use single-block mode (group all primal variables $\x$ together in one block) with localized equality constraints for the sub-problem and apply (Eq. \ref{eq:single_p_diag}).

\item If $H$ is non-diagonal, both $H$ and $A$ are very sparse, and the problem is of moderate size, use single-block ADMM.  
If only a subset of primal variables is bounded, solve the block using an external solver (e.g. Gurobi or Mosek) with localized bounds. Otherwise, solve the block using (Eq. \ref{eq:single_p}).

\end{enumerate}

\subsubsection{Regularized Markowitz Mean-Variance Model}
\label{subsect:num:cont:markov}

The Markowitz mean-variance model describes $N$ assets characterized by a random vector of returns $R=(R_1,\dots,R_N)$ with known expected value $\m_i$ of each random variable $R_i$ and covariance $\sigma_{ij}$ for all pairs of random variables $R_i$ and $R_j$. Given some portfolio asset $\x=(x_1,\dots,x_N)$, where $x_i$ is the fraction of resources invested in asset $i$, an investor chooses  a portfolio $\x$, satisfying two objectives: expected value of the portfolio return $\m_{\x}=E(R_{\x})=\langle\m,\x\rangle$ is maximized and portfolio risk, measured by variance $\sigma_{\x}^2=\hbox{Var}(R_{\x})=\langle\x,V\!\x\rangle$, $V=(\sigma_{ij})$ is minimized \cite{evstigneev:2015}. 
The problem of finding the optimal portfolio can be formulated as a quadratic optimization problem,

\myeqmodel{\label{eq:mark:cont}
 \min\limits_{\x} & \x^T V \x - \tau\m^T\x+\kappa\|\x\|^2_2\\[0.3cm]
 \mbox{s.t.}& \e^T\x= \bo\\[0.2cm]
    & \x\in\R^n_+
 }
where $\tau \ge 0$ is {\it risk tolerance} parameter. The above problem formulation includes the regularization term with parameter $\kappa$.

The raw data was collected by the Center for Research in Security Price (CRSP), and provided through Wharton Research Data Services \cite{WRDS} covering daily prices of 4628 assets from Jan 01 to Dec 31, 2018, and monthly prices for 7958 stocks from Jan 31 to Dec 31, 2018. Missing data was filled using the yearly average price. The model uses risk tolerance parameter $\tau=1$, and is regularized with $\kappa=10^{-5}$. 
For the formulation $(\ref{eq:mark:cont})$, because Hessian ($V$) is dense and non-diagonal, the multi-block ADMM is used, following the rules on choosing the RACQP solver mode (rule 1, Section \ref{sect:racqp_rules}).  The number of groups $p$ is 50, and the augmented Lagrangian penalty parameter $\beta=1$. Default run settings (Table \ref{tbl:exp_defaults}) are used by all solvers, except for OSQP that had max iteration number set to 20000.

\begin{table}[h!]
\footnotesize
\centering
\setlength\tabcolsep{5pt}
\begin{tabular}{lrrccrrccrrcc}  
\toprule
\multirow{3}{*}{Solver}&\multicolumn{4}{c}{Quarterly Data  ($n=7958$)}&\multicolumn{4}{c}{Monthly Data ($n=7958$)}&\multicolumn{4}{c}{Daily Data ($n=4628$)}\\
\cmidrule(lr){2-5}\cmidrule(lr){6-9}\cmidrule(lr){10-13}
&\multirow{2}{*}{Run}&\multirow{2}{*}{Num.}& \multirow{2}{*}{Res.} &\multirow{2}{*}{Res} 
&\multirow{2}{*}{Run}&\multirow{2}{*}{Num.}& \multirow{2}{*}{Res.} &\multirow{2}{*}{Res} 
&\multirow{2}{*}{Run}&\multirow{2}{*}{Num.}& \multirow{2}{*}{Res.} &\multirow{2}{*}{Res} \\
\uph&\multicolumn{1}{c}{time} & \multicolumn{1}{c}{iter} & \multicolumn{1}{c}{primal} &\multicolumn{1}{c}{dual}
&\multicolumn{1}{c}{time} & \multicolumn{1}{c}{iter} & \multicolumn{1}{c}{primal} &\multicolumn{1}{c}{dual}
&\multicolumn{1}{c}{time} & \multicolumn{1}{c}{iter} & \multicolumn{1}{c}{primal} &\multicolumn{1}{c}{dual}\\
\midrule
\up Gurobi & 2305 & 11 & 6.3$\cdot$10$^{-14}$ & 1.1$\cdot$10$^{-8}$ & 2525 & 12 & 9.2$\cdot$10$^{-15}$ & 7.5$\cdot$10$^{-7}$ & 731 & 15 & 6.0$\cdot$10$^{-15}$ & 4.9$\cdot$10$^{-7}$\\ 
Mosek & 162 & 4 & 1.8$\cdot$10$^{-5}$ & 9.8$\cdot$10$^{-9}$ & 188 & 5 & 3.3$\cdot$10$^{-5}$ & 1.9$\cdot$10$^{-8}$ & 68 & 10 & 5.2$\cdot$10$^{-5}$ & 4.7$\cdot$10$^{-6}$\\ 
OSQP & 3599 & 20000 & 1.5$\cdot$10$^{-3}$ & 3.8$\cdot$10$^{-7}$ & 4359 & 20000 & 5.6$\cdot$10$^{-5}$ & 4.8$\cdot$10$^{-7}$ & 639 & 11475 & 1.0$\cdot$10$^{-5}$ & 1.1$\cdot$10$^{-8}$\\ 
RACQP & 97 & 356 & 7.3$\cdot$10$^{-6}$ & 8.6$\cdot$10$^{-6}$ & 314 & 1191 & 7.0$\cdot$10$^{-7}$ & 9.5$\cdot$10$^{-6}$ & 38 & 576 & 1.7$\cdot$10$^{-7}$ & 1.0$\cdot$10$^{-5}$\\
\bottomrule
\end{tabular}
\caption{Markowitz min-variance model (\ref{eq:mark:cont}). CRSP 2018 data \cite{WRDS}, run-time in seconds.}\label{tbl:markov:real}
\end{table}

The performance comparison between the solvers, given in Table \ref{tbl:markov:real}, shows that multi-block RAC finds the solution of high quality in a fraction of time needed by the commercial solvers. In addition, the results show that OSQP requires many iterations to converge to a solution meeting primal/dual tolerance criteria ($\epsilon = 10^{-5}$), confirming the slow convergence issue of a 2-block ADMM approach. 

{\it Low-rank re-formulation} 

Noting that the number of observations $k$ is not large and that the covariance matrix $V$ is of low rank and thus can be expressed as \myeq{V=B^TB}, where
\myeql{B = \dfrac{1}{\sqrt{k-1}}(R-\dfrac{1}{k}\e\e^T R)}
and $R\in\R^{k\times N}$, with rows corresponding to time series observations, and columns corresponding to different assets, we reformulate the problem as
\myeqmodel{\label{eq:mark:cont2}
 \min\limits_{\x} & \|\y\|^2_2 - \tau\m^T\x + \kappa\|\x\|^2_2 \\[0.3cm]
 \mbox{s.t.}& \e^T\x= \bo\\[0.2cm]
    & B\x - \y = \bm{0}\\[0.2cm]
    & \x\in\R^n_+
 }
Since the Hessian of $(\ref{eq:mark:cont2})$ is diagonal, and number of constraints is relatively small, the problem is solved using the single-block ADMM (rule 3, Section \ref{sect:racqp_rules}). Run-time settings are identical to those used for the regular model described previously, with the exception of the augmented Lagrangian penalty parameter which is set to $\beta=0.1$. 
The performance comparison between the solvers, given in Table \ref{tbl:markov:real_LR}, shows that RACQP is also competitive in low-rank formulation of the problem.

\begin{table}[h!]
\footnotesize
\centering
\setlength\tabcolsep{5pt}
\begin{tabular}{lrrccrrccrrcc}  
\toprule
\multirow{3}{*}{Solver}&\multicolumn{4}{c}{Quarterly Data  ($n=7958$)}&\multicolumn{4}{c}{Monthly Data ($n=7958$)}&\multicolumn{4}{c}{Daily Data ($n=4628$)}\\
\cmidrule(lr){2-5}\cmidrule(lr){6-9}\cmidrule(lr){10-13}
&\multirow{2}{*}{Run}&\multirow{2}{*}{Num.}& \multirow{2}{*}{Res.} &\multirow{2}{*}{Res} 
&\multirow{2}{*}{Run}&\multirow{2}{*}{Num.}& \multirow{2}{*}{Res.} &\multirow{2}{*}{Res} 
&\multirow{2}{*}{Run}&\multirow{2}{*}{Num.}& \multirow{2}{*}{Res.} &\multirow{2}{*}{Res} \\
\uph&\multicolumn{1}{c}{time} & \multicolumn{1}{c}{iter} & \multicolumn{1}{c}{primal} &\multicolumn{1}{c}{dual}
&\multicolumn{1}{c}{time} & \multicolumn{1}{c}{iter} & \multicolumn{1}{c}{primal} &\multicolumn{1}{c}{dual}
&\multicolumn{1}{c}{time} & \multicolumn{1}{c}{iter} & \multicolumn{1}{c}{primal} &\multicolumn{1}{c}{dual}\\
\midrule
\up Gurobi & 0.1 & 11 & 5.8$\cdot$10$^{-14}$ & 2.7$\cdot$10$^{-9}$ & 0.1 & 10 & 3.8$\cdot$10$^{-13}$ & 2.8$\cdot$10$^{-7}$ & 2.2 & 17 & 1.7$\cdot$10$^{-15}$ & 9.7$\cdot$10$^{-7}$\\ 
Mosek & 0.6 & 5 & 3.6$\cdot$10$^{-5}$ & 3.7$\cdot$10$^{-9}$ & 0.2 & 5 & 6.2$\cdot$10$^{-4}$ & 6.3$\cdot$10$^{-8}$ & 1.3 & 3 & 5.6$\cdot$10$^{-4}$ & 2.6$\cdot$10$^{-8}$\\ 
OSQP & 1.5 & 2600 & 2.5$\cdot$10$^{-9}$ & 3.1$\cdot$10$^{-6}$ & 1.0 & 1175 & 9.3$\cdot$10$^{-6}$ & 9.0$\cdot$10$^{-8}$ & 12.5 & 1900 & 9.6$\cdot$10$^{-6}$ & 9.6$\cdot$10$^{-9}$\\ 
RACQP & 0.6 & 350 & 2.0$\cdot$10$^{-12}$ & 8.9$\cdot$10$^{-6}$ & 0.8 & 767 & 1.3$\cdot$10$^{-12}$ & 9.3$\cdot$10$^{-6}$ & 5.8 & 561 & 1.5$\cdot$10$^{-13}$ & 1.0$\cdot$10$^{-5}$\\ 
\bottomrule
\end{tabular}
\caption{Low-rank reformulation Markowitz min-variance model (\ref{eq:mark:cont2}). CRSP 2018 data \cite{WRDS}, run-time in seconds.}\label{tbl:markov:real_LR}
\end{table}

\subsubsection{Randomly Generated Linearly Constrained Quadratic Problems (LCQP)}
\label{sect:rnd_lcqp}
In this section we analyze RACQP performance for different problem structures and run-time settings (number of blocks $p$, penalty parameter $\beta$, tolerance $\epsilon$). In order to have more control over problem structure we generate synthetic problem instances starting with a simple one row Markowitz-like problem to  multi-row problems of large sizes. Note that although we compare RACQP with Gurobi and Mosek on randomly generated instances, which may be considered to be unfair to the latter, our goal is not to diminish the importance of barrier type solution methods those solvers utilize, but to show that multi-block ADMM can be an approach to argument these methods when instances are large and/or dense. In this section we solve linearly constrained quadratic problems LCQP, described by (\ref{eq:problem_model_QP}), with $\x\in\R^n$.

Similarly to \cite{wright:2015} we construct a positive definite Hessian matrix $H$ from a random ($\sim U(0,1)$) matrix $U\in\R^{n\times n}$ and a normalized diagonal matrix $V\in\R_+^n$ whose elements were chosen from a log-uniform distribution to have a specific condition number:
\myeqmodel{\label{eq:randomGen}
U_{\eta}&=\eta U + (1-\eta)I\\[0.2cm]
H&=U_{\eta}V U_{\eta}^T + \zeta  \e\e^T
 }

where parameters $\eta\in(0,1)$ and $\zeta\ge 0$ induce different types of orientation bias. For convenience we normalize matrix $H$ and construct vector $\cc$ as a random vector ($\sim U(0,1)$).  Jacobian matrices $A_{eq}$ and $A_{ineq}$ are constructed in a way that the desired sparsity is met and $a_{i.j}\sim N(0,1)$ for both matrices. 
 Our analysis of LCQP is based on extensive experimentation using different problem structure embedded in the matrix $H$, by varying its orientation, condition number and the random seed used to construct $H$ (and vector $\cc$).  

{\it Markowitz-like Problem Instances}

RACQP implementation allows solving optimization problems by multi-block ADMM. A question that arises is the optimal number of blocks $p$ (i.e. sub-problems) to use. The optimal number, it turns out, is related to structure and density of both Hessian and Jacobian matrices. For any $H$ that is not a block matrix, and a dense $A$, as is the case with the Markowitz model, the number of blocks is related to the problem size -- having more blocks leads to having more iterations before the process meets the tolerance on residual error $\epsilon$ and more sub-problems to construct and solve. However, a sub-problem of a smaller size can be constructed and solved in less time than a larger sub-problem. Total time ($t_T$) is thus a function of opposing arguments. To show this interdependence, we solve simple Markowitz-like problem instances, with randomly generated $Q$ and $\cc$, and with $A_{eq}=\e^T$, $\bb = \bo$, and $\x\in\R^n_+$ (inequity constraints are not used). Following (\ref{eq:mark:cont}), we added a regularization term to the objective function with $\kappa=10^{-5}$. 
 
\begin{table}[h!]
\footnotesize
\centering
\begin{tabular}{rrccccccc}  
\toprule
\multirow{2}{*}{Num}&\multirow{2}{*}{Block}&\multicolumn{4}{c}{Number of iterations ($k$)}& \multicolumn{3}{c}{Cost per iteration [s]}\\
\cmidrule(lr){3-6}\cmidrule(lr){7-9}
 blocks & size & $\mu$ & $\sigma^2$ & min & max & $\mu$ & min & max  \\
\midrule
\up 50 & 180& 43.2 & 1.87 & 40 & 46 & 0.147 & 0.133 & 0.152\\
100 & 90 & 46.6 & 1.51 & 44 & 49 & 0.095 & 0.094 & 0.095\\
150 & 60 & 49.0 & 1.25 & 47 & 51 & 0.091 & 0.090 & 0.092\\
200 & 45 & 50.6 & 0.97 & 49 & 52 & 0.108 & 0.107 & 0.109\\
\bottomrule
\end{tabular}
\caption{RACQP performance with respect to number of blocks $p$ for randomly generated problems of type (\ref{eq:mark:cont}). Problem size $n=9000$, density$_Q$=0.05, $\epsilon = 10^{-5}$. \label{tbl:markov:block_time}}
\end{table}

Table \ref{tbl:markov:block_time} presents the aggregate results collected over a set of experiments (10 for each group size) using random problems constructed using (\ref{eq:randomGen}). The reason for constructing problems in such a way is to emulate a real-world situation when a problem model (Hessian, Jacobian, $\x$ upper and lower bounds) do not change, but coefficients do.  The results confirm that there exist a ``right'' number of blocks which minimizes overall run-time. For now, choosing that number is based on experience, but we are working on formalizing the procedure. 

In addition to run-time cost per iteration, Table \ref{tbl:markov:block_time} reports number of iterations until convergence ($k$) for different number of blocks. It is interesting to observe is that $k$ is very mildly affected by the choice of $p$, if tolerance $\epsilon$ is kept the same. This leads to another interesting question on how much a change in $\epsilon$ affect run-time. Table \ref{tbl:markov:epsilon} gives an answer to this question. The table lists RACQP performance over the same problem set, but with  different residual tolerances. As expected, results show that the number of iterations increases as the
tolerance gets tighter.

\begin{table}[h!]
\footnotesize
\centering
\begin{tabular}{cccccccc}  
\toprule
\multirow{2}{*}{$\epsilon$}&\multicolumn{4}{c}{Number of iterations ($k$)}& \multicolumn{2}{c}{Residuals (mean values)}\\
\cmidrule(lr){2-5}\cmidrule(lr){6-7}
 & $\mu$ & $\sigma^2$ & min & max& primal & dual \\
\midrule
\up 10$^{-4}$ & 30.4 & 1.43 & 28 & 33 &4.7$\cdot$10$^{-7}$ & 9.3$\cdot$10$^{-5}$\\
10$^{-5}$ & 46.6 & 1.51 & 44 & 49 & 4.7$\cdot$10$^{-8}$ & 9.5$\cdot$10$^{-6}$ \\
10$^{-6}$ & 63.4 & 1.65 & 60 & 65 &3.6$\cdot$10$^{-9}$ & 9.2$\cdot$10$^{-7}$ \\
10$^{-7}$ & 79.9 & 2.02 & 76 & 83 &3.5$\cdot$10$^{-10}$ & 9.5$\cdot$10$^{-8}$ \\
\bottomrule
\end{tabular}
\caption{A typical RACQP performance with respect to primal/dual residual tolerance $\epsilon$ for a randomly generated problems of type (\ref{eq:mark:cont}). Problem size $n=9000$, density$_Q$=0.05. Number of blocks $p=100$. Run-time shown is seconds.}\label{tbl:markov:epsilon}
\end{table}

\vspace{10pt}
{\it General LCQP}

Building on the results from the previous section, we expand the QP model to include general equality and inequality constraintswith unbounded variables $\x$. We analyze RACQP when solving sparse problems (dense problems are covered in the next section where we address relaxed QAP) for problems of size $n=6000$ and $n=9000$. The number of rows in both constraint matrices is equal ($m=m_{eq}=m_{ineq}$), and set to be a function of a problem size, $m=r\cdot n$, with $r=\{0.1,0.5\}$. The number of blocks used by RACQP is related to size of a block, $p_n=n/b_{\hbox{size}}$, with the optimal block size $b_{\hbox{size}}$ empirically determined to be 60. The penalty parameter $\beta=1$ was found to produce the best results.

Table \ref{tbl:rand:cont} presents comparative analysis of performance of the solvers with respect to run-time and primal/dual residuals. Although both OSQP and RACQP did well in terms of primal and dual residuals, the results show that multi-block RACQP converges to solutions much faster (4-10x) then OSQP. Both solvers outperform Gurobi and  Mosek in run-time, even though the tolerance on residual error is set to the same value ($\epsilon=10^{-5}$). 
Another observation is that Mosek produces solutions of inferior quality to all aforementioned solvers -- dual residuals are of $10^{-3}$ and $10^{-4}$ levels, far below the requested $\epsilon$ threshold. Investigation of the log files produced by Mosek reveled two problems: (1) Mosek terminates as soon as primal or dual or complementary gap residual criteria is met (unlike the other solvers which terminate when {\it all} the residual criteria  are met); (2) residuals are not periodically checked on a re-scaled model, resulting in a large discrepancy between internally evaluated residuals (scaled data) and the actual one. 

\begin{table}[h!]
\footnotesize
\setlength\tabcolsep{3pt}
\centering
\begin{tabular}{crccccccccccccccccc}  
\toprule  
&&&&&\multicolumn{1}{c}{}&\multicolumn{8}{c}{Residuals}\\
\cmidrule(lr){7-14}
\multirow{2}{*}{Problem}&\multirow{2}{*}{Num}&\multicolumn{4}{c}{Run-time [s]}
& \multicolumn{2}{c}{Gurobi} & \multicolumn{2}{c}{Mosek} & \multicolumn{2}{c}{OSQP} & \multicolumn{2}{c}{RACQP}\\
\cmidrule(lr){3-6}\cmidrule(lr){7-8}\cmidrule(lr){9-10}\cmidrule(lr){11-12}\cmidrule(lr){13-14}
 size & rows & Gurobi & Mosek &OSQP &RACQP  & primal & dual & primal & dual & primal & dual & primal & dual\\
\midrule
\up 
6000 & 600 & 1082 & 208 & 84 & 9 & 2.8$\cdot$10$^{-13}$ & 1.5$\cdot$10$^{-10}$ & 1.4$\cdot$10$^{-6}$ & 5.1$\cdot$10$^{-5}$ & 1.5$\cdot$10$^{-9}$ & 2.3$\cdot$10$^{-8}$ & 3.4$\cdot$10$^{-7}$ & 9.8$\cdot$10$^{-6}$\\ 
 & 3000 & 1861 & 143 & 98 & 26 & 1.2$\cdot$10$^{-12}$ & 1.4$\cdot$10$^{-10}$ & 9.3$\cdot$10$^{-6}$ & 8.2$\cdot$10$^{-3}$ & 1.4$\cdot$10$^{-9}$ & 4.0$\cdot$10$^{-8}$ & 1.4$\cdot$10$^{-6}$ & 9.6$\cdot$10$^{-6}$\\ 
 & 900 & 4222 & 365 & 293 & 22 & 5.7$\cdot$10$^{-13}$ & 1.9$\cdot$10$^{-9}$ & 7.0$\cdot$10$^{-6}$ & 6.9$\cdot$10$^{-4}$ & 4.9$\cdot$10$^{-9}$ & 4.3$\cdot$10$^{-7}$ & 1.6$\cdot$10$^{-7}$ & 8.7$\cdot$10$^{-6}$\\ 
9000 & 4500 & 6308 & 408 & 304 & 65 & 1.9$\cdot$10$^{-12}$ & 7.8$\cdot$10$^{-9}$ & 8.8$\cdot$10$^{-6}$ & 3.0$\cdot$10$^{-5}$ & 2.8$\cdot$10$^{-9}$ & 5.5$\cdot$10$^{-8}$ & 9.5$\cdot$10$^{-7}$ & 9.4$\cdot$10$^{-6}$\\
\bottomrule
\end{tabular}
\caption{Performance comparison between solvers for LCQP. Density=0.05, $p_{n=6000}\eq100$, $p_{n=9000}\eq150$}\label{tbl:rand:cont}
\end{table}

\subsubsection{Relaxed QAP}
\label{subsect:qap_relaxed}

As of this section we continue the study of RACQP but, instead of randomly generating problems, we use benchmark test sets compiled by other authors which reflect real-world problems. We start by addressing large scale instances from the QAPLIB benchmark library \cite{bench:qaplib} compiled by \cite{Bur:97} and hard problems of large size, described in \cite{Dre:05}. The quadratic assignment problem (QAP) problem is a binary problem, but for the purpose of more realistic comparison between the solvers, we relax it to a continuous problem. The numerical tests solving the binary problem formulation will be given later in Section \ref{subsect:num:bin:qap}. 

The quadratic assignment problem belongs to a class of combinatorial optimization problems that arise from problems of practical interest. The QAP objective is to assign $n$ facilities to $n$ locations in such a way that the assignment cost is minimized. The assignment cost is the sum, over all pairs, of a weight or flow between a pair of facilities multiplied by the distance between their assigned locations. Mathematically, the QAP can be presented as follows: 

\myeqmodel{\label{eq:qaplib:bin}
 \min\limits_{X} & \vect(X)^TH\vect(X)\\[0.3cm]
 \mbox{s.t.}& \sum_{i=1}^r x_{ij}=1,\ \forall j=1,\dots r\hspace{40pt}\hbox{(a)}\\[0.2cm]
    & \sum_{j=1}^r x_{ij}=1,\ \forall i=1,\dots r\hfill\hbox{(b)}\\[0.2cm]
    & 0\le x_{ij},\ \forall i,j=1,\dots r \hfill\hbox{(c)}
 }
where  $x_{ij}$ is the entry of the permutation matrix $X\in\R^{r\times r}$. To make the problem convex and be admitted by Cholesky factorization, we make $H\in\R^{n\times n}$ strict diagonally dominant, $H=\hat H+d\cdot I$, where $\hat H = (A\otimes B)$ and $d=\max(\sum_{i=1,i\not=j}^{n}\hat h_{i,j})+\delta$, with $\delta$ being some small positive number and $n=r^2$. The ``flow'' matrix $A\in\R^{r\times r}$ and the ``distance'' matrix $B\in\R^{r\times r}$.

For QAP we apply a method for variance reduction as described in Section \ref{subsec:grp} since the assignment constraints are highly structured and observable. We group variables following a simple reasoning -- given that the permutation matrix $X$ is doubly stochastic, each row (or column) can be seen as a single super-variable, an integer representing a permutation order. Thus, it makes sense to make one super-variable, $\x_i$  for each row $i$ of $X$, so that each super-variable is of size $r$. For each of the experiments whown we set number of groups $p=r$ (thus we solve for one super-variable per block), and penalty parameter $\beta$ to the best we found by running multiple experiments with different parameter values. We found that $\beta=r$ offered the best run-time.

\begin{table}[h!]
\footnotesize
\setlength\tabcolsep{5pt}
\centering
\begin{tabular}{crrrrrrrrrr}  
\toprule  
\multirow{2}{*}{Instance}&\multirow{2}{*}{Problem}&\multirow{2}{*}{Density}&\multicolumn{4}{c}{Run-time [s]}&\multicolumn{4}{c}{Num. iterations}\\
\cmidrule(lr){4-7}\cmidrule(lr){8-11}
 name & size ($n$)  &\multicolumn{1}{c}{(H)} &  Gurobi& Mosek & OSQP& RACQP& Gurobi &Mosek & OSQP& RACQP  \\
\midrule
\up dre110 & 12100 & 0.03 & 375 & 587 & 1259 & 12 & 10 & 6 & 50 & 45\\
sko100a & 10000 & 0.68 & 4305 & 408 & 401 & 13 & 8 & 7 & 50 & 22\\
sko100f & 10000 & 0.67 & 4694 & 396 & 405 & 14 & 9 & 7 & 50 & 23\\
tai100a & 10000 & 0.96 & 4214 & 419 & 416 & 12 & 8 & 7 & 50 & 20\\
tai125e01 & 15625 & 0.29 & limit & 1820 & 1544 & 16 & 5 & 8 & 50 & 22\\
tho150 & 22500 & 0.42 & limit & 4088 & 4586 & 91 & 1 & 7 & 50 & 26\\
wil100 & 10000 & 0.88 & 4529 & 497 & 409 & 13 & 8 & 9 & 50 & 19\\
\bottomrule
\end{tabular}
\caption{Relaxed QAP \cite{Bur:97, Dre:05} instances. Run-time and iteration count comparison between the solvers. }\label{tbl:qap:cont_time}
\end{table}

\begin{table}[h!]
\footnotesize
\centering
\begin{tabular}{ccccccccc}  
\toprule  
\multirow{2}{*}{Instance}&\multicolumn{2}{c}{Gurobi}&\multicolumn{2}{c}{Mosek}&\multicolumn{2}{c}{OSQP}&\multicolumn{2}{c}{RACQP}\\
\cmidrule(lr){2-3}\cmidrule(lr){4-5}\cmidrule(lr){6-7}\cmidrule(lr){8-9}
 name & primal & dual & primal & dual & primal & dual & primal & dual\\
\midrule
\up 
dre110 & 9.0$\cdot$10$^{-11}$ & 2.6$\cdot$10$^{-8}$ & 2.1$\cdot$10$^{-6}$ & 9.7$\cdot$10$^{-2}$ & 8.7$\cdot$10$^{-10}$ & 2.1$\cdot$10$^{-7}$ & 7.0$\cdot$10$^{-7}$ & 9.0$\cdot$10$^{-6}$ \\
sko100a & 1.7$\cdot$10$^{-13}$ & 4.5$\cdot$10$^{-8}$ & 1.5$\cdot$10$^{-6}$ & 1.8$\cdot$10$^{-2}$ & 8.6$\cdot$10$^{-10}$ & 1.6$\cdot$10$^{-7}$ & 1.3$\cdot$10$^{-6}$ & 6.0$\cdot$10$^{-6}$ \\
sko100f & 2.5$\cdot$10$^{-12}$ & 4.8$\cdot$10$^{-8}$ & 5.4$\cdot$10$^{-6}$ & 4.9$\cdot$10$^{-2}$ & 8.6$\cdot$10$^{-10}$ & 1.6$\cdot$10$^{-7}$ & 1.6$\cdot$10$^{-6}$ & 8.9$\cdot$10$^{-6}$ \\
tai100a & 7.5$\cdot$10$^{-13}$ & 1.2$\cdot$10$^{-9}$ & 1.5$\cdot$10$^{-7}$ & 3.0$\cdot$10$^{-3}$ & 8.6$\cdot$10$^{-10}$ & 1.6$\cdot$10$^{-7}$ & 6.4$\cdot$10$^{-6}$ & 4.7$\cdot$10$^{-6}$ \\
tai125e01 & NA & NA & 2.1$\cdot$10$^{-7}$ & 8.4$\cdot$10$^{-5}$ & 8.7$\cdot$10$^{-10}$ & 2.1$\cdot$10$^{-7}$ & 1.9$\cdot$10$^{-6}$ & 9.1$\cdot$10$^{-6}$  \\
tho150 & NA & NA & 2.3$\cdot$10$^{-6}$ & 2.3$\cdot$10$^{-2}$ & 8.7$\cdot$10$^{-10}$ & 2.5$\cdot$10$^{-7}$ & 1.2$\cdot$10$^{-6}$ & 7.4$\cdot$10$^{-6}$  \\
wil100 & 3.2$\cdot$10$^{-13}$ & 4.5$\cdot$10$^{-9}$ & 9.3$\cdot$10$^{-6}$ & 1.5$\cdot$10$^{-1}$ & 8.6$\cdot$10$^{-10}$ & 1.6$\cdot$10$^{-7}$ & 9.8$\cdot$10$^{-6}$ & 8.1$\cdot$10$^{-6}$  \\
\bottomrule
\end{tabular}
\caption{Relaxed QAP \cite{Bur:97, Dre:05}  instances. Primal and dual residuals comparison between the solvers.}\label{tbl:qap:cont_res}
\end{table}

The results showing performance of solvers on a selected set of large QAP instances are summarized in Tables \ref{tbl:qap:cont_time} and \ref{tbl:qap:cont_res}. The instances were chosen in such a way to cover a  variety of problem densities (Hessian) and sizes.  Table \ref{tbl:qap:cont_time} shows run-time and number of iterations. Note that any comparison between barrier based solvers (Gurobi and Mosek) and ADMM solvers (RACQP, OSQP) is not possible, as the solution methods are completely different, but giving the number of iterations allow us to compare performances within each class of the solvers.

Similarly to results presented previously, RACQP is the fastest solver. Solution quality (primal and dual residual tolerance) is achieved in a fraction of time required by the other solvers.  The average speedup is 214x, 86x and 83x with respect to Gurobi, Mosek and OSQP respectively.
OSQP, although performing a similar number of iterations as RACQP does, is much slower -- splitting a large problem into two parts (OSQP executes 2-block ADMM) still leaves two large matrices to solve! On the positive side, OSQP finds better solutions (primal residual smaller by the order of magnitude). Mosek is the worst performing solver -- run-time-wise it is close to OSQP, only one returned solution satisfies the dual residual (tai125e01). The other instances report the dual to be as low as $10^{-1}$. Gurobi found the best solutions, except for tai125e01 and tho150 instances, when max run-time limit (3h) was reached. 

\subsubsection{Maros and Meszaros Convex QP}
\label{sect:cuteR}

The Maros and Meszaros test set \cite{bench:cute_r} is a collection of convex quadratic programming examples from a variety of sources \cite{maros:1999} of the following form
\myeqmodel{\label{eq:cute}
 \min\limits_{\x} & \half \x^T H \x + \cc^T\x+c_0\\[5pt]
 \mbox{s.t.}&  A\x = \bb\\[5pt]
    & \lb\le\x\le\ub
 }
with $H\in\R^{n\times n}$ symmetric positive definite, $A\in\R^{m\times n}$, $\bb\in\R^m$ and $\lb, \ub\in\R^n$, meaning that some of components of $\lb$ and $\ub$ may be $-\infty$ and $+\infty$ respectively. Constant $c_0$ is assumed to be $|c_0|<\infty$. 

As in the previous section, only a subset of instances is used in experiments. The instances were chosen in such a way to cover a  variety of problem models (density, size) but also to point to strengths and weaknesses of ADMM-based algorithms. Problem sizes $n$ range from $4\cdot10^3$ to almost $10^5$ with the number of constraints $m$  up to $10^5$. The Hessians are diagonal matrices, with number of non-zero diagonal elements less or equal to $n$. The constraint matrices $A\in\R^{m\times n}$ are very sparse across the problems; for most of the instances density is below $10^{-3}$. In addition to being sparse, the Jacobian matrices are not block separable. 

RACQP mode was set to a single-block mode according to the rules 3 and 4 of Section \ref{sect:racqp_rules}, with $\beta=1$ for all instances except for CONT* and UBH1 which use $\beta=350$ and $\beta=12000$ respectively. 
Residual tolerance of  $\epsilon=10^{-4}$ was used in producing the results, reported in Tables \ref{tbl:cuter_time} and \ref{tbl:cuter_res}. The tolerance is lower than the default one ($10^{-5}$) because ADMM methods had hard time converging on CONT* and CVXQP* instances for tighter residuals (max number of iterations limit is 4000).

\begin{table}[h!]
\footnotesize
\setlength\tabcolsep{5pt}
\centering
\begin{tabular}{crrrrrrrrrrrr}  
\toprule  
\multirow{2}{*}{Instance}&\multirow{2}{*}{Problem}&\multirow{2}{*}{Num.}&\multirow{2}{*}{Density}&\multicolumn{4}{c}{Run-time [s]}&\multicolumn{4}{c}{Num. iterations}\\
\cmidrule(lr){5-8}\cmidrule(lr){9-12}
 name & size ($n$) & rows & \multicolumn{1}{c}{(A)}  &  Gurobi& Mosek & OSQP& RACQP& Gurobi &Mosek & OSQP& RACQP  \\
\midrule
\up 
AUG2DC & 20200 & 10000 & 2.0$\cdot10^{-4}$ & 0.1 & 0.8 & 0.1 & 0.6 & 1 & 4 & 50 & 1\\
AUG2DQP & 20200 & 10000 & 2.0$\cdot10^{-4}$ & 0.3 & 0.6 & 1.6 & 3.8 & 15 & 12 & 800 & 238\\
AUG3DC & 3873 & 1000 & 1.7$\cdot10^{-3}$ & 0.0 & 0.1 & 0.0 & 0.0 & 1 & 4 & 50 & 1\\
AUG3DQP & 3873 & 1000 & 1.7$\cdot10^{-3}$ & 0.0 & 0.1 & 0.0 & 0.1 & 13 & 10 & 100 & 154\\
BOYD1 & 93261 & 18 & 3.3$\cdot10^{-1}$ & 0.8 & 2.0 & 31.6 & 10.3 & 21 & 19 & 3325 & 826\\
CONT-050 & 2597 & 2401 & 1.9$\cdot10^{-3}$ & 0.1 & 0.2 & 1.2 & 2.4 & 10 & 10 & 2100 & 2058\\
CONT-100 & 10197 & 9801 & 5.0$\cdot10^{-4}$ & 0.4 & 1.0 & 15.4 & 10.5 & 10 & 13 & limit & 852\\
CONT-101 & 10197 & 10098 & 5.0$\cdot10^{-4}$ & 0.3 & 0.8 & 15.8 & 40.0 & 9 & 11 & limit & 2839\\
CONT-300 & 90597 & 90298 & 1.0$\cdot10^{-4}$ & 5.5 & 12.0 & 278.0 & 877.8 & 10 & 13 & limit & 3405\\
CVXQP1\_L & 10000 & 5000 & 3.0$\cdot10^{-4}$ & 23.5 & 15.3 & 77.6 & 67.9 & 10 & 8 & limit & limit\\
CVXQP2\_L & 10000 & 2500 & 3.0$\cdot10^{-4}$ & 6.4 & 15.4 & 23.0 & 4.3 & 9 & 15 & 1475 & 248\\
DTOC3 & 14999 & 9998 & 2.0$\cdot10^{-4}$ & 0.0 & 0.2 & 0.2 & 0.1 & 1 & 4 & 275 & 65\\
HUES-MOD & 10000 & 2 & 0.99 & 0.0 & 0.1 & 0.1 & 0.0 & 10 & 6 & 200 & 34\\
HUESTIS & 10000 & 2 & 0.99 & 0.0 & 0.1 & 0.0 & 0.0 & 11 & 8 & 75 & 34\\
UBH1 & 18009 & 12000 & 2.0$\cdot10^{-4}$ & 0.1 & 0.1 & 0.1 & 0.5 & 5 & 4 & 75 & 1\\
\bottomrule
\end{tabular}
\caption{Large Maros and Meszaros \cite{bench:cute_r} instances. Run-time and iteration count comparison between the solvers.}\label{tbl:cuter_time}
\end{table}

\begin{table}[h!]
\footnotesize
\centering
\begin{threeparttable}
\begin{tabular}{ccccccccc}  
\toprule  
\multirow{2}{*}{Instance}&\multicolumn{2}{c}{Gurobi}&\multicolumn{2}{c}{Mosek}&\multicolumn{2}{c}{OSQP}&\multicolumn{2}{c}{RACQP}\\
\cmidrule(lr){2-3}\cmidrule(lr){4-5}\cmidrule(lr){6-7}\cmidrule(lr){8-9}
 name & primal & dual & primal & dual & primal & dual & primal & dual\\
\midrule
AUG2DC & 8.0$\cdot$10$^{-13}$ & 9.7$\cdot$10$^{-7}$ & 9.4$\cdot$10$^{-8}$ & 6.8$\cdot$10$^{-8}$ & 3.6$\cdot$10$^{-10}$ & 3.9$\cdot$10$^{-11}$ & 1.6$\cdot$10$^{-12}$ & 1.0$\cdot$10$^{-16}$\\
AUG2DQP & 1.1$\cdot$10$^{-14}$ & 1.5$\cdot$10$^{-8}$ & 7.5$\cdot$10$^{-5}$ & 5.5$\cdot$10$^{-2}$ & 2.3$\cdot$10$^{-6}$ & 5.1$\cdot$10$^{-5}$ & 8.9$\cdot$10$^{-8}$ & 9.8$\cdot$10$^{-5}$\\  
AUG3DC & 1.4$\cdot$10$^{-14}$ & 9.5$\cdot$10$^{-7}$ & 4.8$\cdot$10$^{-10}$ & 8.3$\cdot$10$^{-11}$ & 7.9$\cdot$10$^{-10}$ & 4.4$\cdot$10$^{-10}$ & 1.2$\cdot$10$^{-14}$ & 8.8$\cdot$10$^{-17}$\\  
AUG3DQP & 8.9$\cdot$10$^{-16}$ & 3.0$\cdot$10$^{-7}$ & 3.1$\cdot$10$^{-4}$ & 9.2$\cdot$10$^{-3}$ & 4.2$\cdot$10$^{-5}$ & 5.5$\cdot$10$^{-6}$ & 9.6$\cdot$10$^{-5}$ & 2.6$\cdot$10$^{-6}$\\  
BOYD1 & 9.4$\cdot$10$^{-15}$ & 1.5$\cdot$10$^{-11}$ & 1.4$\cdot$10$^{-6}$ & 8.9$\cdot$10$^{-4}$ & 1.7$\cdot$10$^{-8}$ & 9.9$\cdot$10$^{-5}$ & 6.7$\cdot$10$^{-5}$ & 9.9$\cdot$10$^{-5}$\\  
CONT-050 & 2.2$\cdot$10$^{-15}$ & 2.5$\cdot$10$^{-7}$ & 4.0$\cdot$10$^{-6}$ & 9.6$\cdot$10$^{-11}$ & 4.1$\cdot$10$^{-6}$ & 2.2$\cdot$10$^{-6}$ & 9.9$\cdot$10$^{-5}$ & 1.1$\cdot$10$^{-5}$\\  
CONT-100 & 2.6$\cdot$10$^{-14}$ & 5.3$\cdot$10$^{-7}$ & 1.4$\cdot$10$^{-6}$ & 2.3$\cdot$10$^{-9}$ & 2.6$\cdot$10$^{-4}$ & 5.0$\cdot$10$^{-8}$ & 9.8$\cdot$10$^{-5}$ & 8.3$\cdot$10$^{-5}$\\  
CONT-101 & 4.3$\cdot$10$^{-10}$ & 3.6$\cdot$10$^{-7}$ & 9.8$\cdot$10$^{-6}$ & 4.6$\cdot$10$^{-7}$ & 1.8$\cdot$10$^{-3}$ & 6.2$\cdot$10$^{-7}$ & 9.9$\cdot$10$^{-5}$ & 9.8$\cdot$10$^{-5}$\\  
CONT-300 & 9.4$\cdot$10$^{-9}$ & 9.6$\cdot$10$^{-7}$ & 4.9$\cdot$10$^{-8}$ & 3.0$\cdot$10$^{-6}$ & 8.8$\cdot$10$^{-3}$ & 1.1$\cdot$10$^{-5}$ & 9.9$\cdot$10$^{-5}$ & 9.1$\cdot$10$^{-5}$\\  
CVXQP1\_L & 4.5$\cdot$10$^{-8}$ & 3.0$\cdot$10$^{-8}$ & 4.7$\cdot$10$^{-5}$ & 1.5$\cdot$10$^{-3}$ & 1.2$\cdot$10$^{-4}$ & 1.1$\cdot$10$^{-5}$ & 7.6$\cdot$10$^{-3}$ & 2.5$\cdot$10$^{-5}$\\  
CVXQP2\_L & 9.2$\cdot$10$^{-12}$ & 2.8$\cdot$10$^{-11}$ & 9.7$\cdot$10$^{-9}$ & 6.7$\cdot$10$^{-5}$ & 6.2$\cdot$10$^{-5}$ & 5.3$\cdot$10$^{-8}$ & 4.0$\cdot$10$^{-6}$ & 9.9$\cdot$10$^{-5}$\\  
DTOC3 & 5.2$\cdot$10$^{-11}$ & 6.3$\cdot$10$^{-9}$ & 6.7$\cdot$10$^{-10}$ & 2.4$\cdot$10$^{-13}$ & 9.2$\cdot$10$^{-5}$ & 2.6$\cdot$10$^{-6}$ & 7.2$\cdot$10$^{-13}$ & 9.9$\cdot$10$^{-5}$\\  
HUES-MOD & 2.8$\cdot$10$^{-15}$ & 3.4$\cdot$10$^{-7}$ & 7.2$\cdot$10$^{-5}$ & 1.0$\cdot$10$^{-1}$ & NA$^*$ & NA$^*$ & 8.5$\cdot$10$^{-5}$ & 8.6$\cdot$10$^{-6}$\\  
HUESTIS & 9.8$\cdot$10$^{-15}$ & 4.8$\cdot$10$^{-9}$ & 3.9$\cdot$10$^{-6}$ & 1.1$\cdot$10$^{-2}$ & NA$^*$ & NA$^*$ & 8.5$\cdot$10$^{-5}$ & 8.6$\cdot$10$^{-6}$\\  
UBH1 & 1.6$\cdot$10$^{-10}$ & 8.2$\cdot$10$^{-9}$ & 2.8$\cdot$10$^{-4}$ & 1.3$\cdot$10$^{-4}$ & 9.8$\cdot$10$^{-5}$ & 7.8$\cdot$10$^{-6}$ & 1.5$\cdot$10$^{-5}$ & 6.0$\cdot$10$^{-8}$\\
\bottomrule
\end{tabular}
\begin{tablenotes}
  \item[*] No feasible solution found.
  \end{tablenotes}
\end{threeparttable}
\caption{Large Maros and Meszaros \cite{bench:cute_r}. Primal and dual residuals comparison between the solvers.}\label{tbl:cuter_res}
\end{table}

Overall, for solving sparse and Hessian-diagonal problems, both Gurobi and Mosek seem more robust than OSQP and RACQP, probably due to the linear programming structure. The latter two are of the comparable performance. The results, in terms of the gap are of similar quality, and run-time is approximately the same, except for a couple of instances, where self-adjusting methodology used by OSQP for penalty parameter estimation, gives OSQP speed advantage. Also, some of the run-time variation can also be contributed to different languages used to implement solvers; OSQP is implemented in c/c++ while RACQP uses Matlab. 

RACQP solved more instances than OSQP, which in addition to not being able to meet primal/dual residuals for 25\% of instances, it also could not find a fesible solution for HUES-MOD and HUETIS instances. 
 Mosek residual issue reported in the previous section continues to persists on these problem instances. For example AUG2DQP instance solution has dual residual of $5.5\cdot10^{-2}$, the value that does not meet the requested tolerance. 

\subsubsection{Convex QP based on the Mittelmann LP test set}

In this section we report on the performance of solvers when applied to very large quadratic problems. Instances are taken from the Mittelmann LP test set \cite{bench:lp_test_set} augmented with a diagonal Hessian $H$ to form a standard LCQP (\ref{eq:problem_model_QP}).  The results are shown in Tables \ref{tbl:lp_test_set_time} and \ref{tbl:lp_test_set_res}. Residual tolerance was set to $10^{-4}$ (OSQP could not solve any instance but i\_n13 when default tolerance of $10^{-5}$ was used, and RACQP had hard time with nug30). Other default termination criteria apply (Table \ref{tbl:exp_defaults}). For all instances the number of blocks was set to $p=200$ and penalty parameter, to $\beta=5$ except for nug30 that used $\beta=50$.

\begin{table}[h!]
\footnotesize
\setlength\tabcolsep{5pt}
\centering
\begin{tabular}{crrrrrrrrrrrr}  
\toprule  
\multirow{2}{*}{Instance}&\multirow{2}{*}{Problem}&\multirow{2}{*}{Num.}&\multirow{2}{*}{Density}&
\multicolumn{4}{c}{Run-time [s]}&\multicolumn{4}{c}{Num. iterations}\\
\cmidrule(lr){5-8}\cmidrule(lr){9-12}
 name & size ($n$) & rows & \multicolumn{1}{c}{(A)}  &  Gurobi& Mosek & OSQP& RACQP& Gurobi &Mosek & OSQP& RACQP  \\
\midrule
\up 
nug30 & 753687 & 32769 & 7.91$\cdot$10$^{-5}$ & 9109 & 6738 & limit & 3976 & 13 & 7 & 1101 & 1057\\
wide15 & 753526 & 32762 & 6.10$\cdot$10$^{-5}$ & 9194 & 5267 & limit & 345 & 18 & 13 & 1407 & 136\\
square15 & 753690 & 32769 & 6.10$\cdot$10$^{-5}$ & 9158 & 3733 & limit & 363 & 18 & 9 & 1573 & 142\\
long15 & 379350 & 52260 & 6.10$\cdot$10$^{-5}$ & 8959 & 5238 & limit & 332 & 18 & 13 & 1534 & 136\\
i\_n13 & 741455 & 8192 & 2.44$\cdot$10$^{-4}$ & 71 & 40 & 156 & limit & 28 & 10 & 275 & 2171\\
16\_n14 & 262144 & 16384 & 1.22$\cdot$10$^{-4}$ & 24 & 13 & 267 & 2847 & 38 & 10 & limit & limit\\
\bottomrule
\end{tabular}
\caption{Convex QP based on the Mittelmann LP test set \cite{bench:lp_test_set}. Run-time and iteration count comparison between the solvers.}\label{tbl:lp_test_set_time}
\end{table}

\begin{table}[h!]
\footnotesize
\centering
\begin{tabular}{ccccccccc}  
\toprule  
\multirow{2}{*}{Instance}&\multicolumn{2}{c}{Gurobi}&\multicolumn{2}{c}{Mosek}&\multicolumn{2}{c}{OSQP}&\multicolumn{2}{c}{RACQP}\\
\cmidrule(lr){2-3}\cmidrule(lr){4-5}\cmidrule(lr){6-7}\cmidrule(lr){8-9}
name & primal & dual & primal & dual & primal & dual & primal & dual\\
\midrule
\up 
nug30 & 6.3$\cdot$10$^{-15}$ & 5.5$\cdot$10$^{-9}$ & 5.1$\cdot$10$^{-2}$ & 1.7$\cdot$10$^{-4}$ & 1.7$\cdot$10$^{-4}$ & 9.6$\cdot$10$^{-6}$ & 9.5$\cdot$10$^{-6}$ & 9.3$\cdot$10$^{-5}$\\  
wide15 & 1.9$\cdot$10$^{-13}$ & 5.5$\cdot$10$^{-13}$ & 3.1$\cdot$10$^{-5}$ & 2.2$\cdot$10$^{-3}$ & 4.1$\cdot$10$^{-9}$ & 2.4$\cdot$10$^{-3}$ & 1.5$\cdot$10$^{-5}$ & 7.1$\cdot$10$^{-5}$\\  
square15 & 8.3$\cdot$10$^{-16}$ & 6.4$\cdot$10$^{-13}$ & 4.6$\cdot$10$^{-5}$ & 1.4$\cdot$10$^{-3}$ & 2.1$\cdot$10$^{-9}$ & 7.0$\cdot$10$^{-4}$ & 1.0$\cdot$10$^{-5}$ & 9.4$\cdot$10$^{-5}$\\  
long15 & 1.9$\cdot$10$^{-13}$ & 5.5$\cdot$10$^{-13}$ & 3.1$\cdot$10$^{-5}$ & 2.2$\cdot$10$^{-3}$ & 4.0$\cdot$10$^{-9}$ & 1.1$\cdot$10$^{-3}$ & 1.5$\cdot$10$^{-5}$ & 7.1$\cdot$10$^{-5}$\\  
i\_n13 & 1.5$\cdot$10$^{-13}$ & 2.1$\cdot$10$^{-12}$ & 3.1$\cdot$10$^{-4}$ & 1.5$\cdot$10$^{-1}$ & 6.7$\cdot$10$^{-6}$ & 1.0$\cdot$10$^{-5}$ & 2.3$\cdot$10$^{-3}$ & 8.3$\cdot$10$^{-5}$\\  
16\_n14 & 1.1$\cdot$10$^{-14}$ & 4.4$\cdot$10$^{-12}$ & 2.0$\cdot$10$^{-4}$ & 2.2$\cdot$10$^{-1}$ & 1.0$\cdot$10$^{-4}$ & 7.4$\cdot$10$^{-5}$ & 1.4$\cdot$10$^{-3}$ & 4.2$\cdot$10$^{-5}$\\
\bottomrule
\end{tabular}
\caption{Convex QP based on the Mittelmann LP test set \cite{bench:lp_test_set}. Primal and dual residuals comparison between the solvers.}\label{tbl:lp_test_set_res}
\end{table}

RACQP solved very large ($n>750000$) quadratic problems to the required accuracy ($\epsilon=10^{-5}$) very fast. The results were obtained using different solution strategies: multi-block Cholesky factorization approach for wide15, square15 and long15 instances, and the partial Lagrangian approach for nug30 (localized lower and upper bound of sub-problem primal variables). The best set of parameters were found by a brute-force approach, which implies that additional research work needs to be done to identify algebraic methods to characterize instances so that run-time parameters can be chosen automatically. RACQP was unable to find a solution satisfying both primal and dual residual tolerances for two instances (i\_n13 and 16\_n14), no matter of what run-time settings we used.

OSQP solved only one instance (i\_n13) within given run-time and number of iterations limitations, while Gurobi solved all the instances to a high precison, regardless of having termination criteria, Table \ref{tbl:exp_defaults}, set to $\epsilon=10^{-5}$. Mosek did not find a single solution meeting the residual criteria, due to the aforementioned scaling and termination criteria issue.

\subsubsection{Changing Random Seed for RACQP}
\label{sect:num:rnd_seed}

When it comes to algorithms that are stochastic in nature, as RAC-ADMM is, the question that always comes onto mind is about robustness of the algorithm. More precisely, how much is RAC-ADMM sensitive to variations in problem data for a given problem model, and to variations arising from differences in sub-problems due to from randomness of block building procedure (Algorithm \ref{alg:RACQP}, line \ref{alg:RACQP:m_start}). The answer to the former question has been provided in Section   \ref{sect:rnd_lcqp}, and this section tackles the latter.   

To answer the question on RACQP sensitivity to sub-problem structure we subject RACQP to different {\it random seeds} -- each sub-problem is solving minimization problem defined by Lagrangian (Eq. \ref{eq:Qq}), which is, in turn, a function of blocks of primal variables constructed using a stochastic process, following the procedure outlined in Section \ref{subsect:solver:cont}, Step 1. This stochastic process is guided by values drawn from a pseudo random number generator, which is initialized using a {\it random seed} number. For different seeds the generator produces different sequences of numbers, what in turn produces different sub-problems addressed by RACQP. 

Table \ref{tbl:rnd_seed} shows results over a selected set of instances chosen to represent each problem type addressed so far. The table aggregates statistical data collected by solving each instances using ten different seeds per primal/dual tolerance $\epsilon$. Note that CuteR instances (Section \{sect:cuteR}) are not included in the analysis as all the instances are solved using a single-block approach.

\begin{table}[h!]
\footnotesize
\centering
\begin{tabular}{ccccccccccccc}  
\toprule  
\multirow{2}{*}{Instance}&\multicolumn{4}{c}{$\epsilon=$10$^{-4}$}
&\multicolumn{4}{c}{$\epsilon=$10$^{-5}$}
&\multicolumn{4}{c}{$\epsilon=$10$^{-6}$}\\
\cmidrule(lr){2-5}\cmidrule(lr){6-9}\cmidrule(lr){10-13}
& $\mu$ &$\sigma$& min & max & $\mu$ &$\sigma$& min & max & $\mu$ &$\sigma$& min & max\\
\midrule
\up regular\_monthly & 575.7 & 0.8 & 574 & 577 & 1181.8 & 30.9 & 1094 & 1192 & 2591.1 & 1.1 & 2589 & 2593\\
regular\_daily & 184.6 & 1.6 & 183 & 188 & 579.0 & 3.9 & 576 & 587 & 1242.7 & 2.7 & 1240 & 1248\\
sko100a & 17.5 & 0.5 & 17 & 18 & 21.3 & 0.5 & 21 & 22 & 25.2 & 0.6 & 24 & 26\\
tai125e01 & 18.6 & 0.5 & 18 & 19 & 22.5 & 0.5 & 22 & 23 & 26.6 & 0.5 & 26 & 27\\
tai150b & 27.8 & 0.6 & 27 & 29 & 34.8 & 0.6 & 34 & 36 & 41.4 & 0.8 & 40 & 42\\
tho150 & 21.8 & 0.8 & 21 & 23 & 26.4 & 0.7 & 25 & 27 & 30.6 & 0.8 & 29 & 32\\
wil100 & 15.0 & 0.0 & 15 & 15 & 19.4 & 0.5 & 19 & 20 & 24.0 & 0.0 & 24 & 24\\
square15 & 145.0 & 5.0 & 140 & 157 & 195.3 & 13.1 & 178 & 207 & 241.3 & 9.4 & 215 & 247\\
\bottomrule
\end{tabular}
\caption{RACQP performance -- number of iterations over different random seeds. Ten experiments per instance per $\epsilon$.}\label{tbl:rnd_seed}
\end{table}

The results show that RACQP is a robust algorithm and that using a single run was a correct choice to make, at least when it comes to problem instances reported in this section. To generalize the claim about RACQP robustness with respect to randomness of block building scheme would require much more experiments and theoretical analysis, what we delegate to our future work.

%% file: Empirical_Analysis_Subsec_Bin.tex
\subsection{Binary and Mixed Integer problems}
\label{subsect:num:bin}

The RAC-ADMM multi-block approach can be applied directly to binary (and mixed integer) problems without any adaptation. However, when dealing with combinatorial problems, a divide-and-conquer approach does not necessary lead to a good solution, because solver may get stuck in some local optima. To mitigate this problem RACQP, we introduce additional randomness into the implementation: a simple perturbation scheme shown in Algorithm \ref{alg:RACQP:mip} that helps the solver to ``escape'' the local optimum and to continue search for another one (and possibly find the global optimum). Thus, in addition to the run-time parameters used for continuous problems, for MIP we need to specify perturbation parameters such as probability distribution to use when choosing how many variables are perturbed ($N_p$) and the parameters thereof. As a default, RACQP implements truncated the exponential distribution,  $N_p\sim\hbox{Exp}(\lambda)$ with parameter $\lambda=0.4n$, minimum number of variables $N_{p,\hbox{min}}=2$, and maximum number of variables $N_{p,\hbox{max}}=n$, based on the observation that for most of the problems ``good'' solutions tend to be grouped. Variables are chosen at random, and in the general case, perturbation is done by assigning ``new'' values (within bounds) to the chosen variables. Default number of trials before perturbation, $N_{trial}=\min(2,0.005n)$. For all binary problems presented in this section the primal residual error was zero, i.e. the problems were solved to optimality.

As the default solver for sub-problems, RACQP uses Gurobi, but any other solver that admits mixed integer quadratic problem would suffice. The results reported in this section are based on Gurobi 7.5, and may be outdated. However, since we use Gurobi as the sub-solver, we expect RACQP to implicitly gain by the improvements made to Gurobi. Gurobi was ran using its default run-time settings (e.g. presolve option was turned on).

In \cite{stellato:2018mip} the authors present a mixed integer quadratic solver, MIQPs, which uses OSQP solver for solving sub-problems resulting from branch-and-bound strategy. Since the solver is built for small and medium size problems that occur in embedded applications, we do not include it in our current study. However, given that MIQPs showed a promising numerical performance (3x faster than Gurobi) even though being implemented in Python, it would be interesting to use it within RACQP as the external solver for MIP (Algorithm \ref{alg:RACQP}, line \ref{alg:RACQP:external}) instead of our default solver (Gurobi) and compare performance. We defer this comparison to future work. 

To solve MIP problems RACQP uses the partial Lagrangian approach, described in Section \ref{sect:partLan}, to handle bounds on variables $\X_i$, $\x_i\in\X_i$. Additionally, depending on a problem structure, equality and inequality constraints can also be moved to the local constraint set.
 Our experiments show that moving some (as it done for QAP), or all constraints (e.g. graph cut problems) to a local set is beneficial in terms of block sizes, run-time, and overall solution quality. By using local constraints we help the sub-solver (e.g. Gurobi) reduce the size of the problem and tighten its formulation (using presolve and cutting plane algorithms). 

Rather than solving the binary QP problem exactly, our goal is to find a (randomized or deterministic) algorithm that could find a better solution under a fixed solution time constraint. Our preliminary tests show that solving a large-scale problem using RAC-ADMM based approach can lead to a very good quality solution for an integer problem in a very limited time. 

The quality of solutions is given in a form of a gap between the objective value of the optimal solution $x_{\hbox{opt}}^*$ and the objective value of the solution found by a solver $S$, $x_S^*$: 
 
 \myeql{\label{eq:gap}gap_{S}=\frac{f(x_S^*) - f(x_{\hbox{opt}}^*)}{1+\hbox{abs}(f(x_{\hbox{opt}}^*))}}}

 For the instances for which the optimal solution remains unknown (e.g. QAPLIB and GSET instances), we use the best known results from the literature. Note that for maximization problems (e.g. Max-Cut, Max-Bisection) gap is the negative of (\ref{eq:gap}). All binary problems are solved with primal residual equal to zero (i.e. the solutions are feasible and integer).

\subsubsection{Randomness Helps}

We start the analysis of RACQP for binary problems with a shorth example showing that  having blocks that are randomly constructed at each iteration, as done by RAC-ADMM, is  the main feature that makes RACQP work well for combinatorial problems, without a need for any special adaptation of the algorithm for the discrete domain.

RAC-ADMM can be easily adapted to execute classical ADMM or RP-ADMM algorithms, so here we compare these three ADMM variants when applied to combinatorial problems. We use a small size problem ($n=1000$) and construct a problem using (\ref{eq:randomGen}) applied to problem of Markowitz type (\ref {eq:mark:bin}),
\myeqmodel{\label{eq:mark:bin}
 \min\limits_{\x} & \x^T V \x + \tau\m^T\x+\kappa\|\x\|^2_2\\[0.3cm]

 \mbox{s.t.}& \e^T\x= r\\[0.2cm]
    & \x\in\{0,1\}^n
 }
with $\kappa=10^{-5}$ and a positive integer number $r\in\Z_+$,  $r\in(1,n)$ that defines how many stocks from a given portfolio must be chosen.
For completeness of the comparison, we implemented distributed-ADMM (Eq. \ref{ADMM-D}) for binary problems and ran the algorithm on the same data.
\begin{figure}[h!]
    \centering
        \centering
        \includegraphics[scale=.4]{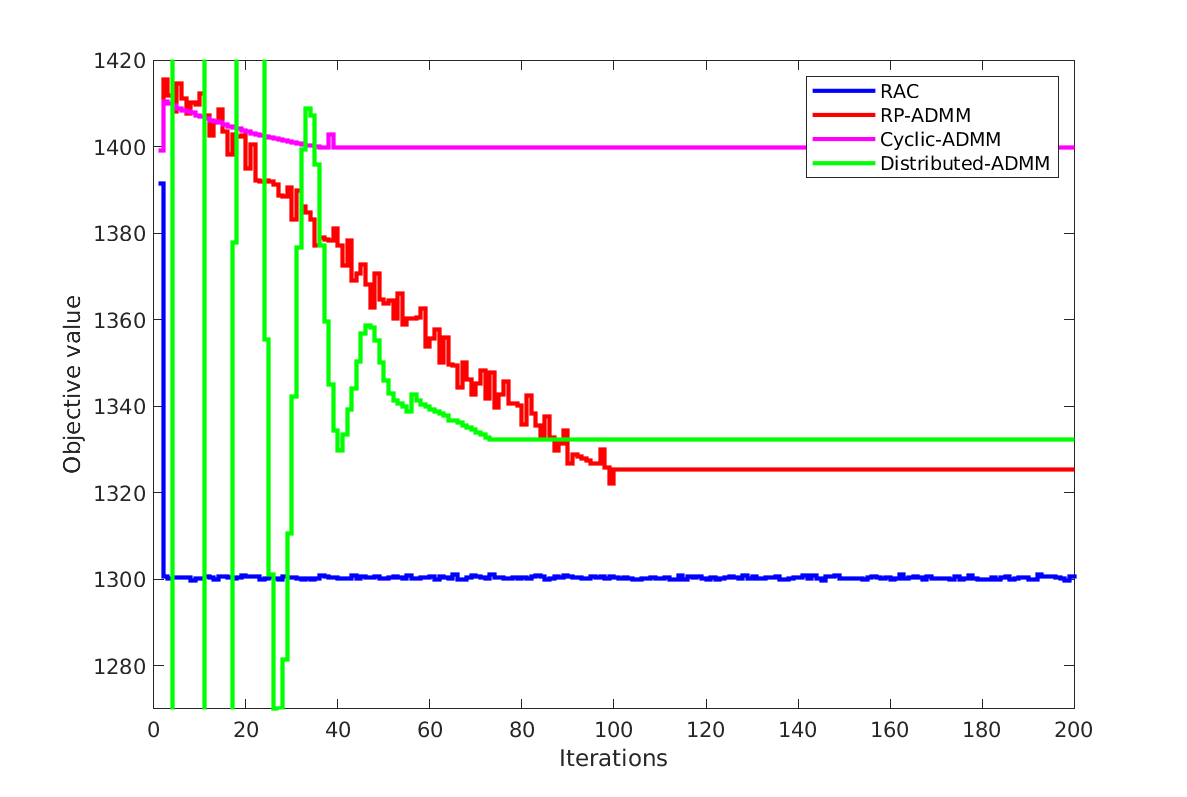}
\caption{A typical evaluation of the objective function value of (\ref{eq:mark:bin}): RAC-ADMM, RP-ADMM (\ref{RP-ADMM}), Cyclic ADMM (\ref{ADMM}) and Distributed ADMM (\ref{ADMM-D})}\label{fig:bin}
\end{figure}
Results show that RAC-ADMM is much better suited for binary optimization problems than either cyclic ADMM or RP-ADMM or distributed-ADMM, which is not surprising since more randomness is adapted into the algorithm making it more likely to escape local optima. All the algorithms are quick to find a local optimum, but besides RAC-ADMM stay at that first found point, while RAC-ADMM continues to find local optima, which could be better or worse than previously found. Because of this behavior, one can keep track the best solution found ($\x_{best}$, Algorithm \ref{alg:RACQP:mip}).  The algorithms seem robust with respect to the structure of the Hessian and choice of initial point. A typical evaluation of the algorithms is shown in Figure \ref{fig:bin}. Note that distributed-ADMM has a very low objective value in early iterations, which is due to the large feasibility errors.

\subsubsection{Markowitz Portfolio Selection}
\label{subsect:num:bin:mark}
Similarly to the section on continuous problems, we compare RACQP performance with that of Gurobi on Markowitz cardinality constrained portfolio selection problem (\ref{eq:mark:bin}) using real data coming from CRSP 2018 \cite{WRDS}.  In the experiments, we set $r=n/2$ with all other settings identical to those used in Section \ref{subsect:num:cont:markov}, including  $V$ and $\m$, estimated from CRSP 2018 data. The default perturbation RACQP settings with $\beta=0.05$, $p=100$ were used in the experiments. Gap is measured from the ``Optimal'' objective values of the solutions found by Gurobi in about 1 hour run-time after relaxing MIPGAP parameter to 0.1. 

\begin{table}[h!]
\centering
\footnotesize
\begin{threeparttable}
\begin{tabular}{ccrcccccc}  
\toprule
\multirow{3}{*}{CRSP 2018}&\multirow{3}{*}{Problem}&\multirow{3}{*}{Optimal}&\multicolumn{6}{c}{Gap}\\
\cmidrule(lr){4-9}
\multirow{2}{*}{data} &\multirow{2}{*}{size ($n$)}&\multirow{2}{*}{Obj. Val.}& \multicolumn{2}{c}{run time = 1 min}&\multicolumn{2}{c}{run time = 5 min}& \multicolumn{2}{c}{run time = 10 min}   \\
\cmidrule(lr){4-5}\cmidrule(lr){6-7}\cmidrule(lr){8-9}
& & &Gurobi$^*$ & RACQP & Gurobi & RACQP & Gurobi & RACQP \\
\midrule
\up quarterly & 7958 & 0.055 & 36.8 & 2.0$\cdot$10$^{-3}$ & 0 & 9.0$\cdot$10$^{-4}$ & 0 & 9.0$\cdot$10$^{-4}$\\
monthly & 7958 & 0.144 & 25.9 & 1.1$\cdot$10$^{-2}$ & -8.7$\cdot$10$^{-6}$ & 2.0$\cdot$10$^{-3}$ & -8.7$\cdot$10$^{-6}$ & 2.0$\cdot$10$^{-3}$\\
daily & 4628 & 1.164 &2.9 & 2.8$\cdot$10$^{-4}$ & 0 & 4.6$\cdot$10$^{-5}$ & 0 & 0\\
\bottomrule
\end{tabular}
\begin{tablenotes}
  \item[*] Root relaxation step not finished. Gurobi returned a heuristic feasible solution.
  \end{tablenotes}
\end{threeparttable}
\caption{Markowitz portfolio selection model (\ref{eq:mark:bin}). CRSP 2018 data \cite{WRDS}.}\label{tbl:markov:real_bin}
\end{table}

From the results (Table \ref{tbl:markov:real_bin}) it is noticeable that RACQP finds relatively good solutions (gap $10^{-2} - 10^{-4}$) in a very short time, in some cases even before Gurobi had time to finalize root relaxation step of its binary optimization procedure. Maximal allowed run-time of 1 min  was far too short for Gurobi to find any solution, so it returned a heuristic ones. Note that those solutions (third column of the table) are extremely weak, suggesting that a RAC-ADMM based solution could be implemented and used instead.

{\it Low-rank Markowitz portfolio selection model}

Similarly to (\ref{eq:mark:cont2}) we formulate the model for low-rank covariance matrix $V$ as

\myeqmodel{\label{eq:mark:bin2}
 \min\limits_{\x} & \|\y\|^2_2 - \tau\m^T\x + \kappa\|\x\|^2_2 \\[0.3cm]
 \mbox{s.t.}& \e^T\x= r\\[0.2cm]
    & B\x - \y = \bm{0}\\[0.2cm]
    & \x\in\{0,1\}^n
 }

and solve the model for CRSP 2018 data. We use $\beta=0.5$, $p=50$. RACQP gap was measured from the optimal solution returned by Gurobi.
 In Table \ref{tbl:markov:real_bin2} we report on the best solutions found by RACQP with max run-time limited to 60 seconds. Results are hard to compare. When Hessian is diagonal and the number of constraints are small, as the case for this data, Gurobi has a very easy time solving the problems (monthly and daily data) -- it finds good heuristic points to start with, and solves problems at a root node after a couple of hundreds of simplex iterations. On the other hand, RACQP, which does not directly benefit from diagonal Hessian, needs to execute multiple iterations of ADMM. Even though the problems are small and solved very quickly, the overhead of preparing the sub-problems and initializing Gurobi to solve sub-problems accumulates to the point of overwhelming RACQP run-time. In that light, for the rest of this section we consider problems where Hessian is a non-diagonal matrix, and address the problems that are hard to solve directly by Gurobi (and possibly other MIP QP solvers).

\begin{table}[h!]
\centering
\footnotesize
\begin{tabular}{ccccc}  
\toprule
\multirow{2}{*}{CRSP 2018}&\multirow{2}{*}{Problem}&\multirow{2}{*}{Optimal}&\multicolumn{2}{c}{Gap}\\
\cmidrule(lr){4-5}
data&size ($n$) & Obj. Val. &Gurobi & RACQP \\
\midrule
\up quarterly & 7958 & 0.015 & -2.2$\cdot$10$^{-7}$ & 1.6$\cdot$10$^{-3}$\\
monthly & 7958 & 0.104 & -4.2$\cdot$10$^{-6}$ & -1.3$\cdot$10$^{-5}$\\
daily & 4628 & 1.140 &3.3 & -1.0$\cdot$10$^{-2}$\\
\bottomrule
\end{tabular}
\caption{Low-rank reformulation Markowitz portfolio selection model (\ref{eq:mark:bin}). CRSP 2018 data \cite{WRDS}. Max run-time 1 min.}\label{tbl:markov:real_bin2}
\end{table}

\subsubsection{QAPLIB}
\label{subsect:num:bin:qap}

The binary quadratic assignment problem (QAP) is known to be NP-hard and that binary instances of larger sizes (dimension of the permutation matrix $r>40$) are considered to be intractable and cannot be solved exactly (though some instances of a large size with special structure have been solved). Currently, the only practical solutions for solving large QAP instances are heuristic methods.

For binary QAP we apply the same method for variance reduction as we did for relaxed QAP (Section \ref{subsect:qap_relaxed}). We group variables following the structure of constraints, which is dictated by the permutation matrix $X\in\{0,1\}^{r\times r}$ (see Eq. \ref{eq:qaplib:bin} for QAP problem formulation) -- we construct one super-variable, $\x_i$  for each row $i$ of $X$. Next we make the use of the partial Lagrangian, and split constraints into the local constraint set consisting of (\ref{eq:qaplib:bin}) (a) and the global constraint set consisting of (\ref{eq:qaplib:bin}) (b), so that the partial Lagrangian is
\myeqln{L_{\beta}(\x,y)=\half\x^TH\x-\y^T(A_{global}\x-\bo)+\frac{
\beta}{2}||A_{global}\x-\bo||^2. } 
At each iteration, we update the $i^{th}$ block by solving 
\myeq{\x_i^{k+1}=\argmin \{L_{\beta}(\cdot) |\, A_{local}\x_i=\bo,\allowbreak \ \x_i\in \{0,1\}^n\}.}
Next, continuing on the discussion on perturbation from the previous section, we turn the feature on and set parameters as follows: number of super-variables to perturb is drawn from truncated exponential distribution, $N_p\sim\hbox{Exp}(\lambda)$ with parameter $\lambda=0.4r$, minimum number of variables $N_{p,\hbox{min}}=2$ and maximum number of variables $N_{p,\hbox{max}}=r$. The number of trials before perturbation $N_{trial}$ is set to its default value.

Note that we do not perturb single variables ($x_{i,j}$), rather super-variables that we choose at random.  If a super-variable $\x_i$ has value of '1' at one location, and '0' on all other entries, then we  randomly swap location of '1' within the super variable (thus keeping the row-wise constraint on $X$ for row $i$ satisfied). If the super-variable is not feasible (number of '1'$\not=1$), we flip values of a random number of variables that make $\x_i$.
The initial point is a random feasible vector. The penalty parameter is a  function of the problem size, $\beta=n$, while the number of blocks depends on the permutation matrix size and it is $p=\lceil r/2\rceil$.

\begin{table}[h!]
\centering
\footnotesize
\begin{tabular}{lcc}   
\toprule
QAPLIB  \cite{Bur:97} benchmark results summary&Gurobi&  RACQP   \\
\midrule
\uph Num. instances opt/best found & 3 &18 \\
Num. instances gap $< 0.01$ (excluding opt/best) & 0 &17 \\
Num. instances gap $< 0.1$ (excluding opt/best and $<0.01$)& 3 &70 \\
\bottomrule
\end{tabular}
\caption{Number of instances = 133. Max run-time: 10 min}\label{tbl:qap:bin_summary}
\end{table}

The summary of the QAPLIB benchmark \cite{bench:qaplib} results is given in Table \ref{tbl:qap:bin_summary}. Out of 133 total instances the benchmark includes, RACQP found the optimal solution (or the best known from literature as not all instances have proven optimal solution) for 18 instances within 10 min of run-time. For the rest of the instances, RACQP returned solutions with an average gap of $\mu=0.07$.
Gurobi solved only three instances to optimality. The average gap of the unsolved instances is $\mu=12.15$, which includes heuristic solutions returned when root relaxation step was not finalized (20 instances). Removing those outliers results in the average gap of $\mu=5.57$.

\begin{table}[h!]
\centering
\footnotesize
\begin{threeparttable}
\begin{tabular}{ccccccc}   
\toprule
\multirow{3}{*}{Instance}&\multirow{3}{*}{Problem}&\multirow{3}{*}{Density}&\multirow{3}{*}{Best known}&\multicolumn{3}{c}{Gap}\\
\cmidrule(lr){5-7}
\multirow{2}{*}{name} &\multirow{2}{*}{size($n$)}&\multirow{2}{*}{($H$)}&\multirow{2}{*}{Obj val}& \multicolumn{1}{c}{Gurobi} &\multicolumn{2}{c}{RACQP}\\
\cmidrule(lr){5-5}\cmidrule(lr){6-7}
&&& &10 min & 5 min & 10 min \\
\midrule
\up lipa80a & 6400 & 0.96 & 253195 & 0.15$^*$ & 0.02 & 0.01\\ 
lipa80b & 6400 & 0.96 & 7763962 & -0.96$^*$ & 0.23 & 0.23\\ 
lipa90a & 8100 & 0.97 & 360630 & 0.22$^*$ & 0.01 & 0.01\\ 
lipa90b & 8100 & 0.97 & 12490441 & -0.96$^*$ & 0.23 & 0.23\\ 
sko81 & 6561 & 0.69 & 90998 & 1.11$^*$ & 0.02 & 0.02\\ 
sko90 & 8100 & 0.68 & 115534 & 1.17$^*$ & 0.04 & 0.03\\ 
sko100a & 10000 & 0.68 & 152002 & 1.34$^*$ & 0.05 & 0.04\\ 
sko100b & 10000 & 0.68 & 153890 & 1.38$^*$ & 0.04 & 0.03\\ 
sko100c & 10000 & 0.67 & 147862 & 1.21$^*$ & 0.04 & 0.03\\ 
sko100d & 10000 & 0.67 & 149576 & 1.21$^*$ & 0.04 & 0.04\\ 
sko100e & 10000 & 0.67 & 149150 & 1.17$^*$ & 0.05 & 0.03\\ 
sko100f & 10000 & 0.67 & 149036 & 1.18$^*$ & 0.04 & 0.03\\ 
tai80a & 6400 & 0.96 & 13499184 & -0.98$^*$ & 0.06 & 0.05\\ 
tai80b & 6400 & 0.43 & 818415043 & -1.00 & 0.26 & 0.22\\ 
tai100a & 10000 & 0.96 & 21043560 & -0.97$^*$ & 0.06 & 0.05\\ 
tai100b & 10000 & 0.43 & 1185996137 & -1.00$^*$ & 0.21 & 0.21\\ 
tai150b& 22500 & 0.44 & 498896643 & -1.00$^*$ & 0.21 & 0.20\\ 
tho40 & 1600 & 0.38 & 240516 & -0.92 & 0.04 & 0.03\\ 
tho150 & 22500 & 0.42 & 8133398 & -0.89$^*$ & 0.08 & 0.06\\ 
wil50 & 2500 & 0.86 & 48816 & 0.53 & 0.01 & 0.01\\ 
wil100 & 10000 & 0.88 & 273038 & 1.19 & 0.03 & 0.02\\ 
\bottomrule
\end{tabular}
\begin{tablenotes}
  \item[*] Root relaxation step not finished. Gurobi returned a heuristic feasible solution.
  \end{tablenotes}
\end{threeparttable}
\caption{QAPLIB, large problems. Gap between best known results  \cite{Bur:97, misevivcius:2019} and RACQP/Gurobi objective values.}\label{tbl:qap:bin_large}
\end{table}

Table \ref{tbl:qap:bin_large} gives detailed information on 21 large instances from QAPLIB data set. The most important takeaway from the table is that Gurobi can not even start solving very large problems as it can not finalize  the root relaxation step within given maximum run time, while RACQP can.

\subsubsection{Maximum Cut Problem}
\label{subsect:num:gset_mc}

The maximum-cut (Max-Cut) problem consists of finding a partition of
the nodes of a graph $G = (V,E)$, into two disjoint sets $V_1$ and $V_2$ ($V_1\cap V_2 = \emptyset$, $V_1\cup V_2 = V$) in such a way that the total weight of the edges that have one endpoint in $V_1$ and the other in $V_2$ is maximized. The problem has numerous important practical applications, and is one of Karp’s 21 NP-complete problems. A standard formulation of the problem is \myeq{ \max\limits_{\y_i\in\{-1,1\}}  \frac{1}{4}\sum_{i,j} w_{i,j}(1-y_iy_j) }, which can be re-formulated into quadratic unconstrained binary problem

\myeqmodel{\label{eq:maxcut}
 \min\limits_{\x} & \x^T H \x \\[0.2cm]
 \mbox{s.t.}& \x\in\{0,1\}^n
 }

where $h_{i,j}=w_{i,j}$ and $h_{i,i}=-\half(\sum_{j=1}^n w_{i,j} + \sum_{j=1}^n w_{j,i})$.

We use the Gset benchmark from \cite{bench:gset}, and compare the results of our experiments with the optimal solutions (found by Gurobi) and the best known solutions from the literature \cite{benlic:2013, ma:2017}. 
For perturbation we use default parameters and perform perturbation by choosing a random number of variables and negating their values, i.e. $x_i=1-x_i$.  
 The number of blocks is equal for all instances, $p=4$, and the initial point is set to zero ($\x_0=\bz$) for all the experiments. Note that as the max-cut problem is unconstrained, the enalty parameter $\beta$ is not used (and RACQP is doing a randomly assembled cyclic BCD).

\begin{table}[h!]
\centering
\setlength\tabcolsep{5pt}
\footnotesize
\begin{tabular}{cccccccccccc}  
\toprule
\multirow{3}{*}{Instance}&\multirow{3}{*}{Problem}&\multirow{3}{*}{Density}&\multirow{3}{*}{Best known}&\multicolumn{8}{c}{Gap}\\
\cmidrule(lr){5-12}
\multirow{2}{*}{name} &\multirow{2}{*}{size($n$)}&\multirow{2}{*}{($H$)}&\multirow{2}{*}{Obj val}& \multicolumn{2}{c}{run time = 5 min} &\multicolumn{2}{c}{run time = 10 min}& \multicolumn{2}{c}{run time = 30 min} & \multicolumn{2}{c}{run time = 60 min}   \\
\cmidrule(lr){5-6}\cmidrule(lr){7-8}\cmidrule(lr){9-10}\cmidrule(lr){11-12}
 &&&& Gurobi & RACQP & Gurobi & RACQP & Gurobi & RACQP & Gurobi & RACQP\\
\midrule
\up G1 & 800 &6.1$\cdot10^{-2}$& 11624 & -0.006 & -0.003 & -0.005 & -0.003 & -0.005 & -0.002 & -0.005 & -0.002\\
G6 & 800 &6.1$\cdot10^{-2}$& 2178 & -0.015 & -0.011 & -0.014 & -0.011 & -0.012 & -0.008 & -0.012 & -0.008\\
G11 & 800 &5.8$\cdot10^{-3}$& 564 & 0 & -0.004 & 0 & -0.004 & 0 & -0.004 & 0 & -0.004\\
G14 & 800 &1.5$\cdot10^{-2}$& 3064 & -0.021 & -0.001 & -0.021 & -0.001 & -0.021 & -0.001 & -0.020 & -0.001\\
G18 & 800 &1.5$\cdot10^{-2}$& 992 & -0.081 & -0.011 & -0.081 & -0.011 & -0.081 & -0.011 & -0.081 & -0.011\\
G22 & 2000 &1.1$\cdot10^{-2}$& 13359 & -0.062 & -0.008 & -0.062 & -0.008 & -0.052 & -0.008 & -0.052 & -0.007\\
G27 & 2000 &1.1$\cdot10^{-2}$& 3848 & -0.157 & -0.152 & -0.155 & -0.152 & -0.152 & -0.151 & -0.149 & -0.141\\
G32 & 2000 &2.3$\cdot10^{-3}$& 1410 & 0 & -0.014 & 0 & -0.014 & 0 & -0.014 & 0 & -0.013\\
G36 & 2000 &6.4$\cdot10^{-3}$& 7678 & -0.026 & -0.004 & -0.026 & -0.004 & -0.026 & -0.004 & -0.026 & -0.004\\
G39 & 2000 &6.3$\cdot10^{-3}$& 2408 & -0.102 & -0.011 & -0.102 & -0.010 & -0.102 & -0.006 & -0.102 & -0.006\\
G43 & 1000 &2.1$\cdot10^{-2}$& 6660 & -0.046 & -0.002 & -0.046 & -0.002 & -0.045 & -0.002 & -0.045 & -0.002\\
G50 & 3000 &1.7$\cdot10^{-3}$& 5880 & 0 & -0.001 & 0 & -0.001 & 0 & -0.001 & 0 & -0.001\\
G51 & 1000 &1.3$\cdot10^{-2}$& 3848 & -0.021 & -0.008 & -0.021 & -0.003 & -0.021 & -0.003 & -0.021 & -0.003\\
G55 & 5000 &1.2$\cdot10^{-3}$& 10299 & -0.044 & -0.007 & -0.041 & -0.005 & -0.039 & -0.005 & -0.038 & -0.005\\
G56 & 5000 &1.2$\cdot10^{-3}$& 4016 & -0.112 & -0.016 & -0.112 & -0.015 & -0.112 & -0.011 & -0.112 & -0.011\\
G58 & 5000 &2.6$\cdot10^{-3}$& 19276 & -0.054 & -0.008 & -0.040 & -0.007 & -0.039 & -0.005 & -0.039 & -0.004\\
G60 & 7000 &8.4$\cdot10^{-4}$& 14187 &-0.120 & -0.008 & -0.098 & -0.007 & -0.096 & -0.005 & -0.091 & -0.004\\
G61 & 7000 &8.1$\cdot10^{-4}$& 5796 &-0.222 & -0.015 & -0.181 & -0.013 & -0.158 & -0.013 & -0.117 & -0.012\\
G63 & 7000 &1.8$\cdot10^{-3}$& 26997 & -0.046 & -0.007 & -0.046 & -0.006 & -0.032 & -0.005 & -0.032 & -0.005\\
G67 & 10000 &4.6$\cdot10^{-4}$& 6940 & -0.003 & -0.018 & -0.002 & -0.016 & 0.001 & -0.014 & 0.001 & -0.013\\
G70 & 10000 &2.8$\cdot10^{-4}$& 9581 &-0.006 & -0.006 & -0.006 & -0.006 & -0.005 & -0.004 & -0.004 & -0.004\\
G77 & 14000 &3.3$\cdot10^{-4}$& 9926 & -0.016 & -0.017 & -0.010 & -0.017 & 0.001 & -0.013 & 0.001 & -0.012\\
G81 & 20000 &2.3$\cdot10^{-4}$& 14030 & -0.119 & -0.023 & -0.031 & -0.017 & -0.023 & -0.014 & 0.002 & -0.014\\
\up & \multicolumn{3}{c}{\hfill Average:} & -0.0635 & -0.0151 & -0.0546 & -0.0141 & -0.0505 & -0.0128 & -0.0466 & -0.0126\\
\bottomrule
\end{tabular}
\caption{Max-Cut, GSET instances. Gap between best known results and RACQP/Gurobi objective values.}\label{tbl:gset}
\end{table}

In contrast to continuous sparse problems (rule 4, Section \ref{sect:racqp_rules}), sparse binary problems benefit from using a randomized multi-block approach, as shown in Table \ref{tbl:gset}. The table compares RACQP and Gurobi results collected from experiments on Gset instances for three different maximum run-time limit settings, 10, 30 and 60 minutes. RACQP again outperforms Gurobi, overall, it finds better solutions when run-time is limited. Although Gurobi does better on a few problems, on average RACQP is better.
Note that for large(r) problems ($n\ge5000$) RACQP keeps improving, which can be explained by the difference in number of perturbations -- for smaller problems, good points have already being visited and a chance to find a better one are small. Adaptively changing perturbation parameters could help, but this topic is out of scope of this work.

\subsubsection{Maximum Bisection Problem}
\label{subsect:num:gset_bs}

The maximum bisection problem is a variant of the Max-Cut problem that involves partitioning the vertex set $V$ of a graph $G = (V,E)$ into two disjoint sets $V_1$ and $V_2$ of equal cardinality (i.e. $V_1\cap V_2 = \emptyset$, $V_1\cup V_2 = V$, $|V_1| = |V_2|$) such that the total weight of the edges whose endpoints belong to different subsets is maximized. The problem formulation follows (\ref{eq:maxcut}) with the addition of a constraint $\e^T\x= \lfloor n/2 \rfloor$, where $n$ is the graph size. 

For Max-Bisection, at each iteration we would update the $i^{th}$ block by solving 
\myeqln{\x_{\omega_i}^{k+1}=\argmin \{\x_{\omega_i}^T H_{\omega_i} \x_{\omega_i} -y(\e^T\x_{\omega_i}-b_{\omega_i}) +\bhalf\|\e^T\x_{\omega_i}-b_{\omega_i}\|^2 \ | \x_{\omega_i}\in \{0,1\}^{d_i}\}}
where $d_i$ is the size of block $i$, $\x_{\omega_i}$ is a sub-vector of $\x$ constructed of components of $\x$ with indices $\omega_i\in\Omega$, and $b_{\omega_i}=\lfloor n/2 \rfloor-\e^T\x_{-\omega_i}$ with $\x_{-\omega_i}$ being the sub-vector of $\x$ with indices not chosen by $\omega_i$. Solving the sub-problems directly has shown to be very time consuming. However, noticing that Gurobi, while solving the problem as whole, makes a good use of cuts for this type of problems (matrix Q structure), we decided to reformulate the sub-problems as follows

\myeqmodeln{\label{eq:maxbisect}
 \min\limits_{\x} & \x_{\omega_i}^T H_{\omega_i} \x_{\omega_i} -yr +\bhalf r^2 \ \\[0.2cm]
 \mbox{s.t.}&  \e^T\x_{\omega_i}-r=\bb_{\omega_i} \\[0.2cm]
  &\x_{\omega_i}\in\{0,1\}^{d_i},\ r\in\{0,1\}.
 }
Note that $r$ can be also defined as a bounded continuous or integer variable, but because the optimal value is zero and because Gurobi makes good use of binary cuts, we decided to define $r$ as binary. 

As in the previous section, we use Gset benchmark library  and compare the results of our experiments with the best known solutions for max-bisection problems found in the literature \cite{ma:2017}. The experimental setup is identical to  that of Max-Cut experiments except for the use of the penalty parameter $\beta=0.005$ and the initial point $\x_0$ which is a feasible random vector. Perturbation is done with a simple swap -- an equal number of variables with values ``1'' and ``0'' is chosen and the new value set to be the negation of the old value.

\begin{table}[h!]
\centering
\setlength\tabcolsep{5pt}
\footnotesize
\begin{tabular}{cccccccccccc}  
\toprule
\multirow{3}{*}{Instance}&\multirow{3}{*}{Problem}&\multirow{3}{*}{Density}&\multirow{3}{*}{Best known}&\multicolumn{8}{c}{Gap}\\
\cmidrule(lr){5-12}
\multirow{2}{*}{name} &\multirow{2}{*}{size($n$)}&\multirow{2}{*}{($H$)}&\multirow{2}{*}{Obj val}& \multicolumn{2}{c}{run time = 5 min} &\multicolumn{2}{c}{run time = 10 min}& \multicolumn{2}{c}{run time = 30 min} & \multicolumn{2}{c}{run time = 60 min}   \\
\cmidrule(lr){5-6}\cmidrule(lr){7-8}\cmidrule(lr){9-10}\cmidrule(lr){11-12}
 &&&& Gurobi & RACQP & Gurobi & RACQP & Gurobi & RACQP & Gurobi & RACQP\\
\midrule
\up G1 & 800 &6.1$\cdot10^{-2}$& 11624 & -0.004 & -0.005 & -0.004 & -0.005 & -0.004 & -0.005 & -0.002 & -0.001\\ 
G6 & 800 &6.1$\cdot10^{-2}$& 2177 & -0.023 & -0.004 & -0.022 & -0.003 & -0.018 & -0.003 & -0.015 & -0.003\\ 
G11 & 800 &5.8$\cdot10^{-3}$& 564 & 0 & -0.014 & 0 & -0.011 & 0 & -0.007 & m & -0.007\\ 
G14 & 800 &1.5$\cdot10^{-2}$& 3062 & -0.019 & -0.008 & -0.019 & -0.008 & -0.019 & -0.002 & -0.018 & -0.002\\ 
G18 & 800 &1.5$\cdot10^{-2}$& 992 & -0.062 & -0.004 & -0.062 & -0.004 & -0.062 & -0.001 & -0.062 & -0.001\\ 
G22 & 2000 &1.1$\cdot10^{-2}$& 13359 & -0.207 & -0.009 & -0.207 & -0.005 & -0.171 & -0.003 & -0.066 & -0.003\\ 
G27 & 2000 &1.1$\cdot10^{-2}$& 3341 & -0.050 & -0.023 & -0.050 & -0.021 & -0.043 & -0.016 & -0.042 & -0.014\\ 
G32 & 2000 &2.3$\cdot10^{-3}$& 1410 & 0 & -0.010 & 0 & -0.010 & 0 & -0.009 & m & -0.009\\ 
G36 & 2000 &6.4$\cdot10^{-3}$& 7678 & -0.021 & -0.004 & -0.021 & -0.004 & -0.021 & -0.004 & -0.021 & -0.004\\ 
G39 & 2000 &6.3$\cdot10^{-3}$& 2408 & -0.088 & -0.011 & -0.073 & -0.010 & -0.072 & -0.010 & -0.072 & -0.010\\ 
G43 & 1000 &2.1$\cdot10^{-2}$& 6659 & -0.075 & -0.004 & -0.059 & -0.004 & -0.057 & -0.001 & -0.057 & -0.001\\ 
G50 & 3000 &1.7$\cdot10^{-3}$& 5880 & -0.012 & 0.000 & -0.012 & 0.000 & -0.004 & 0.000 & -0.004 & 0.000\\ 
G51 & 1000 &1.3$\cdot10^{-2}$& 3847 & -0.017 & -0.005 & -0.015 & -0.005 & -0.015 & -0.004 & -0.014 & -0.004\\ 
G55 & 5000 &1.2$\cdot10^{-3}$& 10299 & -0.120 & -0.008 & -0.041 & -0.007 & -0.040 & -0.006 & -0.038 & -0.006\\ 
G56 & 5000 &1.2$\cdot10^{-3}$& 4016 & -0.197 & -0.019 & -0.109 & -0.018 & -0.098 & -0.017 & -0.089 & -0.017\\ 
G58 & 5000 &2.6$\cdot10^{-3}$& 19276 & -0.169 & -0.007 & -0.169 & -0.007 & -0.037 & -0.005 & -0.037 & -0.005\\ 
G60 & 7000 &8.4$\cdot10^{-4}$& 14187 & -0.166 & -0.011 & -0.136 & -0.006 & -0.074 & -0.004 & -0.074 & -0.004\\ 
G61 & 7000 &8.1$\cdot10^{-4}$& 5796 & -0.359 & -0.019 & -0.359 & -0.019 & -0.180 & -0.019 & -0.167 & -0.018\\ 
G63 & 7000 &1.8$\cdot10^{-3}$& 26988 & -0.226 & -0.007 & -0.158 & -0.006 & -0.128 & -0.004 & -0.038 & -0.003\\ 
G67 & 10000 &4.6$\cdot10^{-4}$& 6938 & -0.258 & -0.016 & -0.173 & -0.014 & -0.004 & -0.011 & -0.001 & -0.010\\ 
G70 & 10000 &2.8$\cdot10^{-4}$& 9581 & -0.009 & -0.008 & -0.009 & -0.006 & -0.008 & -0.004 & -0.004 & -0.003\\ 
G77 & 14000 &3.3$\cdot10^{-4}$& 9918 & -0.468 & -0.015 & -0.468 & -0.013 & -0.211 & -0.012 & -0.015 & -0.010\\ 
G81 & 20000 &2.3$\cdot10^{-4}$& 14030 & -0.280 & -0.017 & -0.280 & -0.015 & -0.253 & -0.014 & -0.127 & -0.012\\ 
\up & \multicolumn{3}{c}{\hfill Average:}& -0.1348 & -0.0099 & -0.1165 & -0.0087 & -0.0722 & -0.0070 & -0.0459 & -0.0064\\
\bottomrule
\end{tabular}
\caption{Max-Bisection, GSET instances. Gap between best known results and RACQP/Gurobi objective values.}\label{tbl:gset_ec}
\end{table}

The results are shown in Table \ref{tbl:gset_ec}. Compared to the unconstrained max-cut problem, RACQP seems to have less trouble solving max-bisection problem -- adding a single constraint boosted its performance by up to 2x. Gurobi performance on the other worsened. Overall, RACQP outperforms Gurobi, finding better solutions when run-time is limited.
Both Gurobi and RACQP continue gaining on solution quality (gap gets smaller) with longer time limits.

%% file: Empirical_Analysis_Subsec_ML.tex
\subsection{Selected Machine Learning Problems}
\label{sect:ML}

In this section we apply RAC method and RP method to few selected
machine learning (ML) problems related to convex quadratic optimization, namely Linear
Regression (Elastic-Net) and Support Vector Machine (SVM). To solve the former we apply a specialized implementation of RAC-ADMM (available for download at \cite{RACQP:code}), while for the latter we use RACQP solver.

\input{Empirical_Analysis_Subsec_ML_Lasso}

\input{Empirical_Analysis_Subsec_ML_SVM}

%% file: Empirical_Analysis_Subsec_ML_Lasso.tex
\subsubsection{Linear Regression using Elastic Net}
\label{sect:lasso}

For a classical linear regression model, with observed features \myeq{\mathbf{X}\in\mathbb{R}^{n\times p }}, where $n$ is number of observations and $p$ is number of features, one solves the following unconstrained optimization problem
\myeql{\min_{\beta} \ \frac{1}{2n}(y-\mathbf{X}\beta)'(y-\mathbf{X}\beta) + P_{\lambda,\alpha}(\beta)}
with  \myeq{P_{\lambda,\alpha}(\beta)=\lambda\{\frac{1-\alpha}{2}\|\beta\|_2+\alpha\|\beta\|_1\}} used for Elastic Net model. By adjusting $\alpha$ and $\lambda$, one could obtain different models: for ridge regression, $\alpha = 0$, for lasso $\alpha=1$, and for classic linear regression, $\lambda=0$. 
For the problem to be solved by ADMM, we use variable splitting and reformulate the problem as follows

\myeqdl{\label{eq:el_net}
\begin{array}{cl}
     \min\limits_{\beta} & \frac{1}{2N}(y-X\beta)^{T}(y-X\beta)+P_{\lambda,\alpha}(z) \\
     \mbox{s.t.} & \beta - \z \ = \ \bz 
\end{array}
}

Note that in (\ref{eq:el_net}) we follow the standard machine learning Elastic Net notation in which $\beta$ is the decision variable in the optimization formulation, rather than $\x$.

Let $c=-\frac{1}{n}X^Ty$, $A=\frac{X}{\sqrt{n}}$, and let $\gamma$ denote the  augmented Lagrangian penalty parameter with respect to constraint $\beta-z$, and $\xi$ be the dual with respect to constraint $\beta-z$. 
The augmented Lagrangian could then be written as 
\myeqln{
\begin{array}{ll}L_{\lambda}=&\frac{1}{2}\beta^T(A^TA+\gamma I)\beta + (c-\xi)^T\beta+(\xi-\gamma \beta)^Tz+\frac{\gamma}{2} z^Tz+P_{\lambda,\alpha}(z)
\end{array}}
We apply RAC-ADMM algorithm by partitioning $\beta$ into multi-blocks, but solve $z$ as one block. For any given $\beta_{k+1}$, optimizer $z^*_{k+1}$ has the closed form solution. 
\myeqln{z^*_{k+1}(i)(\beta_{k+1}(i), \xi_k(i))=\frac{S(\xi_k(i)-\gamma\beta_{k+1}(i),\lambda
\alpha)}{(1-\alpha)\lambda+\gamma},}
where $\xi_i$ is the dual variable with respect to constraint $\beta_i-z_i=0$, and $S(a,b)$ is soft-threshold operation \cite{friedman:2007}.
\myeqln{S(a,b)=\begin{cases}
-(a-b),  &\textup{if} \  b<|a|, \ a>0 \\
-(a+b),  &\textup{if} \ b<|a|, \ a\leq 0 \\
0,     & \textup{if}\ b\geq |a|  \\
\end{cases}}

In order to solve classic linear regression directly,  $\mathbf{X}^T\mathbf{X}$ must be positive definite which can not be satisfied for $p>n$. However, RAC-ADMM only requires  each sub-block  $\mathbf{X}_{sub}^T\mathbf{X}_{sub}$ to be positive definite, so, as long as block size $s<n$, RAC-ADMM can be used to solve the classic linear regression.

We compare our solver with glmnet \cite{friedman:2010, friedman:2011}  and Matlab lasso implementation on synthetic data (sparse and dense problems) and benchmark regression data from LIBSVM \cite{Chang:2011}.

{\it Synthetic Data}

The data set for dense problems $\mathbf{X}$ is generated uniform randomly with $n=10,000$, $p=50,000$, with zero sparsity, while for the ground truth $\beta^{*}$ we use  standard Gaussian and set sparsity of $\beta^*$ to $0.1$.  Due to the nature of the problem, estimation requires lower feasibility precision, so we fix number of iterations to $10$ and $20$.
Glmnet solver benefits from having a diminishing sequence of $\lambda$, but given that many applications (e.g. see  \cite{mohsen:2015}) require a fixed $\lambda$ value , we decided to use fixed $\lambda$ for all solvers. 
Note that the computation time of RAC-ADMM solver is invariant regardless of whether $\lambda$ is decreasing or fixed.

\begin{table}[h!]
  \centering
  \footnotesize
  \begin{tabular}{ccrrrrrrrr}
    \toprule
\multirow{2}{*}{$\lambda$} & \multirow{2}{*}{Num.} & \multicolumn{4}{c}{Absolute L2 loss} & \multicolumn{4}{c}{Total time [s]} \\ 
\cmidrule(lr){3-6}\cmidrule(lr){7-10}
\uph & iterations & RAC & RP & glmnet & Matlab & RAC & RP & glmnet & Matlab\\ 
\midrule
\multirow{2}{*}{0.01} & 10 & 204.8 & 204.6  & 213.9 & 249.1 & 396.5 &  227.6  & 2465.9 & 1215.2\\ 
 & 20 & 208.1 & 230.2 & 213.9 & 237.1 & 735.2 & 343.9 & 3857.9 & 2218.2\\ 
\uph \multirow{2}{*}{0.1} & 10 & 217.8 &215.6  & 220.5 & 213.1 & 388.7 & 212.5 & 4444.3 &2125.9  \\ 
 & 20 & 272.6  &202.4  &    220.5 & 212.4 & 739.7 & 337.2 &    4452.4 & 2434.6 \\ 

\uph\multirow{2}{*}{1} & 10 & 213.6& 209.0 & 203.1 & 210.5 & 415.3&213.6  & 3021.1 & 1138.9\\ 
 & 20 & 213.8& 212.4 & 210.5 & 203.1 & 686.3& 392.1 & 5295.5 & 1495.6\\
    \bottomrule
  \end{tabular}
  \caption{Comparison on solver performance, dense elastic net model.  Dense problem, $n=10,000,\ p=50,000$}
  \label{tbl:en_synt}
\end{table}

Table \ref{tbl:en_synt} reports on the average cross-validation run-time and the average absolute $L2$ loss for  all possible pairs $(\alpha, \lambda)$ with parameters chosen from  
$\alpha=\{0, \ 0.1, \ 0.2,\dots,1\}$ and $\lambda = \{1,\ 0.01\}$. 
Without specifying, RAC-ADMM solver run-time parameters were identical across the experiments, with augmented Lagrangian penalty parameter $\gamma = 0.1\lambda$ for sparsity $< 0.995$, $\gamma = \lambda$ for sparsity $> 0.995$,  and block size $s=\!=\!100$.

Large scale sparse data set $\mathbf{X}$ is generated uniform randomly with $\{n=40,000, \ p=4,000,000\}$, using sparsity $=0.998$. For ground truth $\beta^*$, the standard Gaussian with sparsity $\beta^*=0.5$ and fixed $\lambda$. Noticing from the previous experiment that increasing a step size from $10$ to $20$ didn't significantly improve prediction error, we fix number of iteration to $10$.

\begin{table}[h!]
  \caption{Comparison on solver performance, elastic net model}
  \label{tbl:en_synt_1}
  \centering
  \footnotesize
  \begin{tabular}{ccrrrrrrrrr}
    \toprule
\multirow{2}{*}{$\lambda$} & \multirow{2}{*}{Num.} & \multicolumn{3}{c}{Avg Absolute L2 loss} & \multicolumn{3}{c}{Best Absolute L2 loss}& \multicolumn{3}{c}{Total time [s]} \\ 
\cmidrule(lr){3-5}\cmidrule(lr){6-8} \cmidrule(lr){9-11}
\uph & iterations & RAC & RP & glmnet & RAC &RP & glmnet & RAC & RP & glmnet \\ 
\midrule
\multirow{1}{*}{0.01} & 10 &   1293.3&  1356.7 & 8180.3 & 745.2& 703.52  &   4780.2& 4116.1&2944.8   &17564.2  \\
\uph\multirow{1}{*}{0.1} & 10 &   777.31& 717.92  & 4050.4 &  613.9& 611.79 &3125.6 & 3756.3 &2989.1  & 12953.7\\
\uph\multirow{1}{*}{1} & 10 &   676.17&  671.23  &3124.5 &  615.7&614.79 &1538.9&3697.8&  3003.8 & 8290.5\\ 
    \bottomrule
  \end{tabular}
  \caption*{  Sparse problem, $n=40,000,\ p=4,000,000$}
\end{table}

Table \ref{tbl:en_synt_1}, report on the average cross-validation run time and the average absolute $L2$ loss for all possible pairs $(\alpha, \lambda)$ with parameters chosen from  
$\alpha=\{0, \ 0.1, \ 0.2,\dots,1\}$ and $\lambda = \{1,\ 0.01\}$.
The table  also shows the best $L2$ loss for each solver.
Because it took more than $10,000$ seconds for Matlab lasso to solve even one estimation, the table reports only comparison between glmnet and RAC.

Experimental results on synthetic data show that RAC-ADMM solver outperforms significantly all other solvers in total time  while being competitive in absolute $L2$ loss. Further RAC-ADMM speedups could be accomplished by fixing block-structure (RP-ADMM).
In terms of run-time, for dense problem, RAC-ADMM is 3 times faster compared with Matlab lasso and 7 times faster compared with glmnet. RP-ADMM is 6 times faster compared with Matlab lasso, and 14 times faster compared with glmnet. For sparse problem, RAC-ADMM is more than 30 times faster compared with Matlab lasso, and 3 times faster compared with glmnet. RP-ADMM is 4 times faster compared with glmnet.

Following Corollary \ref{col:rac_speed} RP-ADMM is slower that RAC-ADMM when convergence is measured in number of iterations, and experimental evidence (Table \ref{tbl:RAC_Markov_comp}) show that it also suffers from slow convergence to a high precision level on L1-norm of equality constraints. However, the benefit of RP-ADMM is that it could store pre-factorized sub-block matrices, as block structure is fixed at each iteration,in contrast to RAC-ADMM which requires  reformulation of sub-blocks at each iteration, what it turn makes each iteration more time-wise costly . In many machine learning problems, including regression, due to the nature of problem, a less precision level is required. This makes RP-ADMM an attractive approach, as  it could converge within fewer steps and potentially be faster than RAC-ADMM. 
In addition, while performing simulations we observed 
that increasing number of iteration does not significantly improve performance of prediction. In fact, absolute $L2$ loss remains similar even when number of iteration is increased to 100. This further gives an advantage to RP-ADMM, as it benefits the most when number of iteration is relatively small.

{\it Benchmark instances, LIBSVM \cite{Chang:2011}}

LIBSVM regression data E2006-tfidf feature size is $150,360$ with number of training and testing data points of  $16,087$ and $3,308$ respectively. The null training error of test set is $221.8758$. Following the findings from the section on synthetic problems and noticing that this dataset is sparse (density=$0.991$), this setup uses fixed number of iterations to $10$, and vary $\lambda = \{1, \ 0.01\}$ and $\alpha = \{0, \ 0.1, \ 0.2,\dots,1\}$. The training set is used to predict $\beta^*$, and the model error (ME) of test set is compared across different solvers.

Table \ref{tbl:en_ben_l} shows the performance of OSQP and Matlab lasso for $\alpha=1$ and $\lambda = 0.01$, and Table \ref{tbl:en_ben} compares compare RAC-ADMM with glmnet. The reason for splitting the results in two tables is related to inefficiency of factorizing a big matrix by OSQP solver and Matlab lasso implementation. Each solver requires more than than $1000$ seconds to solve the problem for even $10$ iterations, making them impractical to use.
On the other hand, glmnet, which uses a cyclic coordinate descent algorithm on each variable, performs significantly faster than OSQP and Matlab lasso. However, glmnet can still be inefficient, as a complete cycle through all $p$ variables requires $O(pN)$ operations \cite{friedman:2010}.

\begin{table}[h!]
  \centering
  \footnotesize
  \begin{tabular}{ccc}
    \toprule
  Solver &   Training ME &Total time [s]\\ 
\midrule
OSQP& 64.0 & 1482.5 \\
Matlab  & 61.1 & 3946.6\\
\bottomrule
\end{tabular}
\caption{E2006-tfidf Performance Summary for Lasso problem ($\alpha =1$, $\lambda = 0.01$)}
\label{tbl:en_ben_l}
\end{table}

\begin{table}[h!]
  \centering
  \footnotesize
  \begin{tabular}{ccccccc}
    \toprule
\multirow{2}{*}{$\lambda$}   &  \multicolumn{3}{c}{Training ME} & \multicolumn{3}{c}{ Total time [s]}  \\ 
\cmidrule(lr){2-4}\cmidrule(lr){5-7}
 \uph &  RAC  & RP &  glmnet  &  RAC &RP &  glmnet\\ 
\midrule
0.01  &  22.4& 22.4 &  29.9  &  106.5& 50.9 &  653.2\\ 
0.1  &  22.1&  22.1&  22.7  &   100.5& 51.9  &  269.3\\
1  &  25.7&25.7  &  23.5  &  102.5  &54.2 & 282.9 \\
    \bottomrule
  \end{tabular}
  \caption{E2006-tfidf performance summary }
  \label{tbl:en_ben}
\end{table}

Results given in Table  \ref{tbl:en_ben} are the averages over run-time and training error collected from experiments with $\alpha=\{0,0.1,\dots,1\}$. The results show that RAC-ADMM is faster than glmnet for all different parameters and that it achieves the best training model error, $22.0954$, among all the solvers.
In terms of run-time, RAC-ADMM is 14 times faster than OSQP, 38 times faster than  Matlab lasso, and 4 times faster than glmnet. RP-ADMM is 28, 18 and 8  times faster than OSQP, Matlab lasso and glmnet, respectively.

For log1pE2006 benchmark , feature size is $4,272,227$, number of training data is $16,087$ and number of testing data is $3,308$. The null training error of test set is $221.8758$ and sparsity of  data is $0.998$.
Similarly to the previous benchmark, the performance results are split into two tables. Table \ref{tbl:en_ben_log_1} shows the performance of OSQP and Matlab lasso, while Table \ref{tbl:en_ben_log} compares RAC-ADMM and glmnet.

\begin{table}[h!]
  \centering
  \footnotesize
  \begin{tabular}{ccc}
    \toprule
  Solver &   Training ME &Total time [s]\\ 
\midrule
OSQP&  66.6  &  11437.4\\
Matlab  & - & $>$3 days\\
\bottomrule
\end{tabular}
\caption{log1pE2006 Performance Summary for Lasso problem ($\alpha =1$, $\lambda = 0.01$)}
\label{tbl:en_ben_log_1}
\end{table}

\begin{table}[h!]
  \centering
  \footnotesize
  \begin{tabular}{ccccccc}
    \toprule
\multirow{2}{*}{$\lambda$}   &  \multicolumn{3}{c}{Training ME} & \multicolumn{3}{c}{ Total time [s]}  \\ 
\cmidrule(lr){2-4}\cmidrule(lr){5-7}
 \uph &  RAC& RP  &  glmnet  &  RAC& RP  &  glmnet\\ 
\midrule
0.01  & 43.0 &  41.8& 22.0  &  962.2&   722.5 &  7639.6\\ 
0.1  &  30.8&   31.8&  22.5  &  978.7&  721.4 &  4945.2\\ 
1  &    32.1& 35.5 &  29.3  &  958.5&  749.2  &  1889.5\\
    \bottomrule
  \end{tabular}
  \caption{log1pE2006 performance summary }
  \label{tbl:en_ben_log}
\end{table}

The results show that RAC-ADMM and RP-ADMM are still competitive and are of same level as glmnet with respect to model error, and all outperform OSQP and Matlab. In terms of run-time, RAC-ADMM is 12 times faster than OSQP, and 5 times faster than glmnet. RP-ADMM is 16 and 7 times faster than OSQP and glmnet, respectively.

%% file: Empirical_Analysis_Subsec_ML_SVM.tex
\subsubsection{Support Vector Machine}
\label{sect:svm}

A Support Vector Machine (SVM) is a machine learning method for classification, regression, and other learning tasks. The method learns a mapping between the features $\x_i\in\R^r$, $i=1,\dots n$ and the target label $y_i\in\{-1,1\}$ of a set of data points using a {\it training set} and constructs a hyperplane \myeq{\w^T\phi(\x)+b} that separates the data set. 
This hyperplane is then used to predict the class of further data points.  
The objective uses Structural Risk Minimization principle which aims to minimize the empirical risk (i.e. misclassification  error) while maximizing the confidence interval (by maximizing the separation margin) \cite{vapnik:1998, vapnik:1995}.

Training an SVM is a convex optimization problem, with multiple formulations, such as C-support vector classification (C-SVC), $\upsilon$-support vector
classification ($\upsilon$-SVC), $\epsilon-$support vector regression ($\epsilon-$SVR), and many more. As our goal is to compare RACQP, a general QP solver, with specialized SVM software and not to compare SVM methods themselves, we decided on using C-SVC (\cite{boser:1992, cortes:1995}), with the dual problem formulated as
\myeqdl{\label{svm}
\begin{array}{cc}
     \min\limits_{\z} & \frac{1}{2}\z^T Q \z\ -\ \e^T\z  \\
     \mbox{s.t.} & \y^T\z \ = \ 0  \\
          & \z\in[0,C]
\end{array}
}
 
with $Q\!\in\!\R^{n\times n}$, $Q\!\succeq\!0$, \myeq{q_{i,j}=y_iy_jK(\x_i,\x_j)}, where \myeq{K(\x_i,\x_j):=\phi(\x_i)^T\phi(\x_j)} is a kernel function, and regularization parameter $C\!>\!0$. The optimal $\w$ satisfies \myeq{\w=\sum_{i=1}^ny_i\z_i\phi(\x_i)},  and 
the bias term $b$ is calculated using the support vectors that lie on the margins (i.e. $0<\z_i<C$) as \myeq{b_i= \w^T\phi(\x_i) - y_i }.  To avoid numerical stability issues, $b$ is then found by averaging over $b_i$.  The decision function is defined with \myeq{f(\x) =\sgn(\w^T\phi(\x)+b)}.

We compare RACQP with LIBSVM \cite{Chang:2011}, due its popularity, and with Matlab-SVM , due to its ease of use. These methods implement specialized approaches to address the SVM problem (e.g. LIBSVM uses a Sequential Minimal Optimization, SMO, type decomposition method  \cite{fan:2005, bottou:2007}), while our approach solves the optimization problem (\ref{svm}) directly. 

The LIBSVM benchmark library provides a large set of instances for SVM, and we selected a representative subset: training data sets with sizes ranging from 20,000 to 580,000; number of features from eight to 1.3 million. We use the test data sets when provided, otherwise, we create test data by randomly choosing 30\% of testing data and report cross-validation accuracy results.

In Table \ref{tbl:svm} we report on model training run-time and accuracy, defined as (num. correctly predicted data)/(total testing data size)$\times$100\%. RAC-ADMM parameters were as follows: max block size $s=100, 500,$ and $1000$ for small, medium and large instances, respectively and augmented Lagrangian penalty $\beta=0.1p$, where $p$ is the number of blocks, which in this case is found to be $p=\lceil n/s \rceil$ with $n$ being the size of training data set. 
In the experiments we use Gaussian  kernel, \myeq{K(\x_i,\x_j)\!=\!\exp(-\frac{1}{2\sigma^2}\|\x_i-\x_j\|^2)}. 
Kernel parameters $\sigma$ and $C$ were estimated by running a {\it grid-check} on  cross-validation. We tried different pairs $(C,\sigma)$ and picked those that returned the best cross-validation accuracy (done using randomly choose 30\% of train data) when instances were solved using RAC-ADMM. Those pairs were then used to solve the instances with LIBSVM and Matlab. The pairs were chosen from a relatively coarse grid, $\sigma,C\in\{0.1, 1, 10\}$ because the goal of this experiment is to compare RAC-ADMM with heuristic implementations rather than to find the best classifier.
Termination criteria were either primal/dual residual tolerance ($\epsilon_p = 10^{-1}$ and $\epsilon_d = 10^{-0}$) or maximum number of iterations, $k=10$, whichever occurs the first. Dual residual was set to such a low value because empirical observations showed that restricting the dual residual does not significantly increase accuracy of the classification but effects run-time disproportionately. Maximum run-time was limited to 10 hours for mid-size problems, and unlimited for the large ones. Run-time is shown in seconds, unless noted otherwise.

\begin{table}[h!]
  \centering
  \footnotesize
\begin{threeparttable}
  \begin{tabular}{lllrrrrrrr}
    \toprule
 \multirow{2}{*}{Instance} & \multirow{2}{*}{Training} 
&\multirow{2}{*}{Testing}   
& \multirow{2}{*}{Num.} & \multicolumn{3}{c}{Accuracy [\%]}  &\multicolumn{3}{c}{Training run-time [s]}  \\ 
\cmidrule(r){5-7}\cmidrule(r){8-10}
\uph name &  set size  
&  set size 
&  features  &  RAC  &  LIBSVM  &  Matlab  &  RAC  &  LIBSVM  &  Matlab\\ 
\midrule
\up
a8a & 22696 & 9865 & 122 & 76.3 & 78.1 & 78.1 & 91 & 250 & 2653\\ 
w7a & 24692 & 25057 & 300 & 97.1 & 97.3 & 97.3 & 83 & 133 & 2155\\ 
rcv1.binary & 20242 & 135480 & 47236 & 73.6 & 52.6 & -- & 78 & 363 & 10+h\\ 
news20.binary$^*$ & 19996 & 5998 & 1355191 & 99.9 & 99.9 & -- & 144 & 3251 & NA\\ 
\up
a9a & 32561 & 16281 & 122 & 76.7 & 78.3 & 78.3 & 211 & 485 & 5502\\ 
w8a & 49749 & 14951 & 300 & 97.2 & 99.5 & 99.5 & 307 & 817 & 20372\\ 
ijcnn1 & 49990 & 91701 & 22 & 91.6 & 91.3 & 91.3 & 505 & 423 & 0\\ 
cod\_rna & 59535 & 271617 & 8 & 79.1 & 73.0 & 73.0 & 381 & 331 & 218\\  
real\_sim$^*$ & 72309 & 21692 & 20958 & 69.5 & 69.5 & -- & 1046 & 9297 & 10+h\\ 
\up
skin\_nonskin$^*$ & 245057 & 73517 & 3 & 99.9 & 99.9 & -- & 2.6h & 0.5h & NA\\  
webspam\_uni$^*$ & 350000 & 105000 & 254 & 64.3 & 99.9 & -- & 13.8h & 11.8h & NA\\ 
covtype.binary$^*$ & 581012 & 174304 & 54 & 91.3 & 99.9 & -- & 16.2h & 45.3h & NA\\
    \bottomrule
  \end{tabular}
\begin{tablenotes}
  \item[*]  No test set provided, using 30\% of randomly chosen data from the training set. Reporting cross-validation accuracy results.
  \end{tablenotes}
\end{threeparttable}
  \caption{Model training performance comparison for SVM }
  \label{tbl:svm}
\end{table}

The results show that RACQP produces classification models of competitive quality as models produced by specialized software implementations in a much shorter time.  RACQP is in general faster than LIBSVM (up to 27x) except for instances where ratio of number of observations $n$ with respect to number of features $r$ is very large. It is noticeable that while producing (almost) identical results as LIBSVM, the Matlab implementation is significantly slower.

For small and mid-size instances (training test size $<$ 100K) we tried, the difference in accuracy prediction is less than 2\%, except for problems where test data sets are much larger than the training sets. In the case of ``rcv1.binary'' instance test data set is 5x larger than the training set, and for ``cod\_rna'' instance is 4x larger. In both cases RACQP outperforms LIBSVM (and Matlab) in accuracy, by 20\%  and  9\%, respectively. 

All instances except for ``news20.binary'' have $n>>r$ and the choice of the Gaussian kernel is the correct one. For instances where the number of features is larger than the number of observations, linear kernel is usually the better choice as the separability of the model can be exploited \cite{woodsend:2011} and problem solved to similar accuracy in a fraction of time required to solve it with the non-linear kernel. 
The reason we used the Gaussian kernel on ``news20.binary' instance is that we wanted to show that RACQP is only mildly affected by the feature set size. Instances of similar sizes but different number of features  are all solved by RACQP in approximately the same time, which is 
in contrast with  LIBSVM and Matlab that are both affected by the feature space size. LIBSVM slows down significantly while Matlab, in addition to slowing down could not solve ''news.binary`` -- the implementation of fitcsvm() function that invokes Matlab-SVM algorithm requires full matrices to be provided as the input which in the case of ''news.binary`` requires 141.3GB of main memory.

``Skin\_nonskin'' benchmark instance ``marks'' a point where our direct approach starts showing weaknesses -- LIBSVM is 5x faster than RACQP because of the fine-tuned heuristics which exploit very small feature space (with respect to number of observations). The largest instance we addressed is ``covtype.binary'', with more than half of million observations and the (relatively) small feature size ($p=54$). For this instance, RACQP continued slowing down proportionately  to the increase in problem size, while LIBSVM experienced a large hit in run-time performance, requiring almost two days to solve the full size problem. This indicates that the algorithms employed by LIBSVM are put to the limit and specialized algorithms (and implementations) are needed to handle large-scale SVM problems. 
RACQP accuracy is lower than that of LIBSVM, but can be improved by tightining residual tolerances under the cost of increased run-time.

For large-size problems RACQP performance degraded, but the success with the mid-size problems suggests that a specialized  ``RAC-SVM'' algorithm could be developed to address very large problems.
Such a solution could merge RAC-ADMM algorithm 
with heuristic techniques to (temporarily) reduce the size of the problem (e.g. \cite{joachims:1998}),  smart 
kernel approximation techniques, probabilistic approach(es) to shrinking the support vector set (e.g. \cite{rudi:2017}), and similar.

%% file: Summary.tex
\section{Summary}
\label{sec:summary}

In this paper, we introduced a novel randomized algorithm, randomly assembled multi-block and cyclic alternating direction method of multipliers (RAC-ADMM), for solving continuous and binary convex quadratic problems.  We provided a theoretical proof of the performance of our algorithm for solving linear-equality constrained continuous convex quadratic programming, including the expected convergence of the algorithm and sufficient condition for almost surely convergence of the algorithm.  We further provided open source code of our solver, RACQP, and numerical results on demonstrating  the efficiency of our algorithm. 

We conducted multiple numerical tests on solving synthetic, real-world, and bench-mark quadratic optimization problems, which include continuous and binary problems. We compare RACQP with Gurobi, Mosek and OSQP for cases that do not require high accuracy, but a strictly improved solution in shortest possible run-time. Computational results show that RACQP, except for a couple of instances with a special structure, finds solutions of a very good quality in a much shorter time than the compared solvers. 

In addition to general linearly constrained quadratic problems we applied RACQP to few selected machine learning problems, Linear Regression, LASSO, Elastic-Net, and SVM. Our solver matches the performance of the best tailored methods such as Glmnet and LIBSVM, and often gives better results than that of tailored methods. In addition, our solver uses much less computation memory space than other ADMM based method do, so that it is suitable in real applications with big data.

The following is a quick summary of the pros and cons of RACQP, implementation of RAC-ADMM, for solving quadratic problems, and suggests the future research.

\begin{itemize}
\item RACQP is remarkably effective for solving continuous and binary convex QP problems when the Hessian is non-diagonal, the constraint matrix are unstructured, or the number of constraints are small. These findings are demonstrated by solving Markowitz portfolio problems with real or random data, and randomly generated sparse convex QP problems. 

\item RACQP, coupled with smart-grouping and a partial augmented Lagrangian, is equally effective when the structure of the constraints is known. This finding is supported by solving continuous and binary bench-mark Quadratic Assignment, Max-Cut, and Max-Bisection problems.
However, efficiently deciding on grouping strategy is also challenging. We plan to build an ``automatic-smart-grouping'' method as a pre-solver for unknown structured problem data. 

\item Computational studies done on binary problems show that RAC-ADMM approach to solving problems offers an advantage over the traditional direct approach (solving the problem as whole) when finding a good quality solution for a large-scale integer problem in a very limited time. However, exact binary QP solvers, such as Gurobi, are needed, because our binary RACQP relies on solving many small or medium sized binary sub-problems. Of course, we plan to explore more high efficiency solvers for medium-sized binary problems for RACQP.

\item The ADMM-based approach, either RACQP or OSQP, is less competitive when the Hessian of the convex quadratic objective is diagonal and the constraints are sparse but structured such as a network-flow type. We believe in this case both Gurobi and Mosek can utilize more efficient Cholesky factorization is that is commonly used by interior-point algorithms for solving linear programs; see more details in Section \ref{subsect:solver:cont}.
 In contrary, RACQP has considerable overhead cost of preparing block data and initialization time of the sub-problem solver, and the time spent on solving diagonal sub-problems was an order of magnitude shorter than time needed to prepare data. 
This, together with the divergence problem of multi-block ADMM, hints that there must be something connected to the problem structure that makes such instances hard for the ADMM-based approach. We plan on conducting additional research to identify problem instances that are well-suited and those that are unsuitable for ADMM.

\item There are still many other open questions regarding RAC-ADMM. For example, there is little work on how to optimally choose run-time parameters to work with RAC-ADMM, including penalty parameter $\beta$, number of blocks, and so for.
\end{itemize}